\documentclass[final,leqno,onefignum,onetabnum]{siamltex}

\usepackage{fullpage}
\usepackage{upgreek}
\usepackage[mathscr]{eucal}
\usepackage{color}
\usepackage[usenames,dvipsnames,svgnames,table]{xcolor}
\usepackage{amsmath,amssymb,amsopn,mathtools}
\usepackage{graphicx}
\usepackage{hyperref}       % internal hyperlink to contents
\hypersetup{
	colorlinks=true,   % boxed links
	citecolor=blue,     % citation links (bibliography)
	filecolor=blue,     % color of file links
	linkcolor=red,    % internal links (sections, pages, etc.)
	urlcolor=blue}      % color of URL links (mail, web)

\usepackage{relsize}
\numberwithin{equation}{section}
\usepackage{enumitem}
\usepackage[normalem]{ulem}
\usepackage{verbatim}
\usepackage{footmisc}
\usepackage{soul}

\newlist{Assumption}{enumerate}{1}
\setlist[Assumption]{label=A\arabic*}

\newlist{steps}{enumerate}{1}
\setlist[steps, 1]{label = Step \arabic*:}

\definecolor{Blue}{rgb}{0,0,1}
\definecolor{Red}{rgb}{1,0,0}
\definecolor{Green}{rgb}{0,1,0}
\definecolor{Bronze}{rgb}{0.8,0.5,0.2}
\definecolor{Violet}{rgb}{0.54,0.17,0.89}

\newcommand{\fontDiscrete}{\mathcal}
\newcommand{\knArg}[1]{k(t^{#1})}
\newcommand{\kn}{\knArg{n}}

\newcommand{\bds}[1]{{\boldsymbol{#1}}}

\newcommand{\defeq}{:=}

\newcommand{\stateContinuousDummyEntryNo}{\boldsymbol{u}}

\newcommand{\stateContinuousDummy}[1]{\bds{\stateContinuousDummyEntryNo}} 
\newcommand{\stateTwoContinuousDummyEntryNo}{\boldsymbol{v}}

\newcommand{\stateTwoContinuousDummy}[1]{\bds{\stateTwoContinuousDummyEntryNo}}

\newcommand{\ndof}{\nspacedof}

\newcommand{\unitvec}{\bds{e}}
\newcommand{\unitvecArg}[1]{\unitvec_{#1}}

\newcommand{\argmin}[1]{\underset{#1}{\text{argmin}}}

\newcommand{\orthogonalmat}{\boldsymbol{Q}}

\newcommand{\righttrianglemat}{\boldsymbol R}

\newcommand{\ndofPorts}[1]{n^p}

\newcounter{remctr}
\newcounter{propctr}
\newcounter{proposctr}
\newcommand{\RR}[1]{\ensuremath{\mathbb{R}^{ #1 }}}
\newcommand{\NN}{\mathbb{N}}
\newcommand{\RRplus}[1]{\ensuremath{\mathbb{R}_+^{ #1 }}}
\newcommand{\natNo}{\NN}
\newcommand{\nat}[1]{\natNo(#1)}

\newcommand{\innat}[1]{\in\nat{#1}}

\newcommand{\Span}[1]{\mathrm{span}\{#1\}}
\newcommand{\FE}{\mathrm{FE}}
\newcommand{\BE}{\mathrm{BE}}
\newcommand{\BDF}{\mathrm{BDF}}
\newcommand{\ADM}{\mathrm{AM}}
\newcommand{\ADB}{\mathrm{AB}}

\newcommand{\identity}[1]{\boldsymbol I_{#1}}

\newcommand{\resSymb}{r}
\newcommand{\resRedSymb}{{\mathsf r}}
\newcommand{\fluxSymb}{f}
\newcommand{\paramSymb}{\mu}
\newcommand{\solSymb}{x}
\newcommand{\obliqueprojector}{\mathcal P}

\newcommand{\spaceSymb}{s}

\newcommand{\timeSymb}{t}

\newcommand{\weightmatSymb}{A}

\newcommand{\basismatspaceSymb}{\Phi}
\newcommand{\basisvecspaceSymb}{\phi}

\newcommand{\samplematSymb}{Z}
\newcommand{\samplevecSymb}{z}
\newcommand{\dummySymb}{v}
\newcommand{\snapshotSymb}{X}

\newcommand{\nsmall}{n}

\newcommand{\nparam}{\nsmall_\mu}
\newcommand{\ntrain}{\nsmall_\text{train}}

\newcommand{\nressample}{\nsmall_{\samplevecSymb}}

\newcommand{\nbasisflux}{\nsmall_\fluxSymb}
\newcommand{\nbasisres}{\nsmall_\resSymb}

\newcommand{\nbasisspace}{{\nsmall_\spaceSymb}}

\newcommand{\nbig}{N}
\newcommand{\nspacedof}{\nbig_\spaceSymb}
\newcommand{\ntimedof}{{\nbig_\timeSymb}}

\newcommand{\paramDomain}{\mathcal D}

\newcommand{\totaltime}{T}

\newcommand{\spatialSubspace}{\fontDiscrete S}

\newcommand{\res}{\boldsymbol \resSymb}

\newcommand{\resArg}[1]{\boldsymbol \resSymb^{#1}}

\newcommand{\resRedApprox}{\tilde{\boldsymbol \resRedSymb}}
\newcommand{\resRedApproxArg}[1]{\resRedApprox^{#1}}

\newcommand{\resn}{\resArg{n}}

\newcommand{\param}{\boldsymbol \paramSymb}

\newcommand{\paramDummy}{\boldsymbol \nu}
\newcommand{\dt}{\Delta \timeSymb}

\newcommand{\sol}{\boldsymbol \solSymb}
\newcommand{\vel}{\dot{\boldsymbol \solSymb}}

\newcommand{\flux}{\boldsymbol \fluxSymb}
\newcommand{\fluxArg}[1]{\flux_{#1}}

\newcommand{\solDummy}{\boldsymbol w}

\newcommand{\timeArg}[1]{t^{#1}}
\newcommand{\timeDummy}{\tau}
\newcommand{\timeDomain}{[0,\totaltime]}

\newcommand{\solArg}[1]{\sol_{#1}}
\newcommand{\solFunc}{\sol}

\newcommand{\solFuncArg}[1]{\solFunc(t^{#1};\param)}

\newcommand{\solapproxFunc}{\solapprox}
\newcommand{\solapproxFuncArg}[1]{\solapproxFunc(t^{#1};\param)}

\newcommand{\solapproxFuncOnlyParam}{\solapproxFunc(\cdot;\param)}

\newcommand{\redsolapproxFunc}{\redsolapprox}
\newcommand{\redsolapproxFuncArg}[1]{\redsolapproxFunc(t^{#1};\param)}

\newcommand{\testbasisArg}[2]{\boldsymbol \xi_{ij}}

\newcommand{\solapproxArg}[1]{\solapprox_{#1}}

\newcommand{\maxError}{\text{maximum relative error}}
\newcommand{\projError}{\text{projection error}}

\newcommand{\solapprox}{\tilde\sol}
\newcommand{\velapprox}{\dot{\tilde\sol}}

\newcommand{\redsolapprox}{\hat\sol}

\newcommand{\redsolapproxArg}[1]{\redsolapprox_{#1}}

\newcommand{\redsol}{\hat\sol}
\newcommand{\redvel}{\dot{\hat\sol}}

\newcommand{\dummy}{\boldsymbol {\dummySymb}}
\newcommand{\reddummy}{\hat{\dummy}}

\newcommand{\redres}{\hat{\res}}

\newcommand{\basismatspace}{\boldsymbol{\basismatspaceSymb}}

\newcommand{\basisvecspace}{\boldsymbol{\basisvecspaceSymb}}

\newcommand{\basismatres}{\basismatspace_\resSymb}

\newcommand{\basisresvecArg}[1]{\boldsymbol \phi_{r,#1}}

\newcommand{\weightmat}{\boldsymbol{\weightmatSymb}}

\newcommand{\snapshots}{\boldsymbol \snapshotSymb}

\newcommand{\samplematNT}{\boldsymbol \samplematSymb}
\newcommand{\samplemat}{\samplematNT^T}

\newcommand{\unscaledAE}{\boldsymbol{A}}
\newcommand{\unscaledEncoder}{\boldsymbol{E}}
\newcommand{\unscaledDecoder}{\boldsymbol{D}}
\newcommand{\normalEncoder}{\mathbf{en}}

\newcommand{\scaledEncoder}{\boldsymbol{h}}
\newcommand{\scaledDecoder}{\boldsymbol{g}}
\newcommand{\nnweight}{\boldsymbol{W}}
\newcommand{\nnbias}{\boldsymbol{b}}
\newcommand{\activation}{\boldsymbol{\sigma}}

\newcommand{\nny}{\boldsymbol{y}}
\newcommand{\mask}{\boldsymbol{S}}
\newcommand{\jacobian}{\boldsymbol{J}}
\newcommand{\NMLSPGHR}{\text{NM-LSPG-HR}}
\newcommand{\NMGalerkinHR}{\text{NM-Galerkin-HR}}

\newcommand{\spaceDomain}{\Omega}
\newcommand{\discreteU}{\boldsymbol{U}}
\newcommand{\discreteV}{\boldsymbol{V}}
\newcommand{\nx}{n_x}
\newcommand{\ny}{n_y}
\newcommand{\nt}{n_t}
\newcommand{\nxy}{n_{xy}}
\newcommand{\RN}[1]{%
  \textup{\uppercase\expandafter{\romannumeral#1}}%
}
\newcommand{\subnetSymbol}{sn}

\newcommand{\bmat}[1]{\begin{bmatrix}#1\end{bmatrix}} % Need \usepackage{amsmath}
\def\Abold{\boldsymbol{A}}

\def\Ubold{\boldsymbol{U}}
\def\Vbold{\boldsymbol{V}}

\def\ubold{\boldsymbol{u}}

\def\Sigmabold{\boldsymbol{\Sigma}}

\def\zerobold{{\bf 0}}

\newtheorem{remark}{Remark}[section]

\usepackage{subfigure, algorithmic, algorithm}% subcaptions for subfigures
\usepackage{amsmath,amssymb}% subcaptions for subfigures
\usepackage{epsfig, algorithmic, algorithm}
\usepackage{multirow, booktabs}
\usepackage{siunitx}
\usepackage{xr}
\mathtoolsset{showonlyrefs=true}
\usepackage{xr}

\title{A fast and accurate physics-informed neural network reduced order model
with shallow masked autoencoder}

\author{
    Youngkyu Kim\thanks{Mechanical Engineering, 
      University of California, Berkeley, CA 94720
      (youngkyu$\_$kim@berkeley.edu, zohdi@berkeley.edu)}\and
    Youngsoo Choi\thanks{Center for Applied Scientific Computing,
      Lawrence Livermore National Laboratory, Livermore,
      CA 94550 (choi15@llnl.gov)}\and
    David Widemann\thanks{Computational Engineering Division,
      Lawrence Livermore National Laboratory,
      Livermore, CA 94550 (widemann1@llnl.gov)}\and
    Tarek Zohdi\footnotemark[1]
}

\begin{document}
\setlength{\abovedisplayskip}{3pt}
\setlength{\belowdisplayskip}{3pt} 
\setlength{\abovedisplayshortskip}{3pt} 
\setlength{\belowdisplayshortskip}{3pt}

\maketitle

\begin{abstract}
Traditional linear subspace reduced order models (LS-ROMs) are able to
  accelerate physical simulations in which the intrinsic solution space falls
  into a subspace with a small dimension, i.e., the solution space has a small
  Kolmogorov $n$-width. However, for physical phenomena not of this type, e.g.,
  any advection-dominated flow phenomena such as in traffic flow, atmospheric
  flows, and air flow over vehicles, a low-dimensional linear subspace poorly
  approximates the solution. To address cases such as these, we have
  developed a fast and accurate physics-informed neural network ROM, namely
  nonlinear manifold ROM (NM-ROM), which can better approximate high-fidelity
  model solutions with a smaller latent space dimension than the LS-ROMs. Our
  method takes advantage of the existing numerical methods that are used to
  solve the corresponding full order models.  The efficiency is achieved by
  developing a hyper-reduction technique in the context of the NM-ROM.
  Numerical results show that neural networks can learn a more efficient latent
  space representation on advection-dominated data from 1D and 2D Burgers'
  equations.  A speedup of up to $2.6$ for 1D Burgers' and a speedup of $11.7$ for
  2D Burgers' equations are achieved with an appropriate treatment of the
  nonlinear terms through a hyper-reduction technique.  Finally, a posteriori
  error bounds for the NM-ROMs are derived that take account of the
  hyper-reduced operators.
\end{abstract}


%\begin{AMS}
%  15A23,35K05,35N20,35L65,65D25,65D30,65F15,65L05,65L06,65L60,65M22
%\end{AMS}

%\pagestyle{myheadings}
%\thispagestyle{plain}
%\markboth{Y. CHOI AND D. COOMBS}{Nonlinear term basis}

\section{Introduction}\label{sec:intro}
  Physical simulations are influencing developments in science, engineering, and
  technology more rapidly than ever before. However, high-fidelity, forward
  physical simulations are computationally expensive and, thus, make intractable
  any decision-making applications, such as design optimization, inverse
  problems, optimal controls, and uncertainty quantification, for which many
  forward simulations are required to explore the parameter space in the outer
  loop.

  To compensate for the computational expense issue, the projection-based
  reduced order models (ROMs) take advantage of both the known governing
  equation and the data. ROMs generate the solution data from the corresponding
  physical simulations and then compress the data to find an intrinsic solution
  subspace, which is represented by a linear combination of basis vectors, i.e.,
  LS-ROMs.  This condensed solution representation is plugged back into the
  (semi-)discretized governing equation to reduce the number of unknowns,
  resulting in an over-determined system, i.e., more equations than unknowns.
  Note that the full governing equations are used to constrain the LS-ROM
  through this substitution. Therefore, this can be considered as a
  physics-informed surrogate model. Additionally, the existing numerical methods
  for the corresponding full order model (FOM) is utilized in the LS-ROM
  solution process.  Therefore, the LS-ROM fully respects the original
  discretization of the governing equations that describe/approximate the
  underlying physical laws, unlike black-box approaches. 

  The LS-ROM approach has been successfully applied to many problems and
  applications, including, but not limited to, rocket nozzle shape design
  \cite{amsallem2015design}, flutter avoidance wing shape optimization
  \cite{choi2020gradient}, topology optimization of wind turbine blades
  \cite{choi2019accelerating}, porous media flow/reservoir simulations
  \cite{ghasemi2015localized, jiang2019implementation, yang2016fast},
  computational electro-cardiology \cite{yang2017efficient}, inverse problems
  \cite{fu2018pod}, shallow water equations \cite{zhao2014pod,
  cstefuanescu2013pod}, computing electromyography \cite{mordhorst2017pod},
  spatio-temporal dynamics of a predator–-prey systems
  \cite{dimitriu2013application}, and acoustic wave-driven microfluidic biochips
  \cite{antil2012reduced}. A survey paper for the projection-based LS-ROM
  techniques can be found in \cite{benner2015survey}. 

  In spite of its successes, the linear subspace solution representation suffers
  from not being able to represent certain physical simulation solutions with a
  small basis dimension, such as advection-dominated or sharp gradient
  solutions. This is because LS-ROMs work only for physical problems in which
  the intrinsic solution space falls into a subspace with a small dimension,
  i.e., the solution space has a small Kolmogorov $n$-width.  Unfortunately,
  even though problems that are advection-dominated or have sharp gradient
  solutions are important, they do not have small Kolmogorov $n$-width. Such
  physical simulations include, but are not limited to, the hyperbolic equations
  with high Reynolds number, the Boltzmann transport equations, and the traffic
  flow simulations.

  Therefore, there have been many attempts to build efficient ROMs for the
  advection-dominated or sharp gradient problems. The attempts can be divided
  mainly into two categories: the first one is to enhance the solution
  representability of the linear subspace by introducing some special treatments
  and adaptive schemes and the second one is to replace the linear subspace
  solution representation with the nonlinear manifold. 

  The effort of enhancing the solution representability of the linear subspace
  includes the artificial viscosity, the Petrov--Galerkin projection applied to
  the computational fluid dynamics problems \cite{carlberg2013gnat,
  carlberg2018conservative, choi2020sns}, the residual discrete empirical
  interpolation approach to handle the Navier-–Stokes equations with a large
  Reynolds number \cite{xiao2014non}, and the space--time ROM
  \cite{choi2020space, choi2019space, taddei2020space} where the temporal as
  well as spatial dimensions were reduced to maximize the compressibility even
  with the advection-dominated problems.  A dictionary-based model reduction
  method was developed in \cite{abgrall2016robust} where $\ell_1$ minimization
  is used to project onto the reduced linear subspace.  A fail-safe $h$-adaptive
  algorithm was developed in \cite{carlberg2015adaptive} where the reduced
  linear subspace basis vectors are broken algebraically to enrich the solution
  subspace.  The shifted proper orthogonal decomposition (POD) was introduced to
  address the issue that arises from the advection-dominated problems
  \cite{reiss2018shifted} where a transport operator is incorporated within the
  POD process.  The drawback with this approach is that the speed of the
  transport operator must be known a priori.  In a similar spirit, the transport
  reversal was introduced in \cite{rim2018transport}, which was inspired by the
  template fitting \cite{kirby1992reconstructing}.  The windowed least-squares
  Petrov--Galerkin model reduction for dynamical systems with implicit time
  integrators is introduced in \cite{parish2019windowed}, which can overcome the
  challenges arising from the advection-dominated problems by representing only
  a small time window with a local ROM.  In order to capture the sharp gradient
  accurately,  many approaches use localization strategies. The examples of such
  methods include the online adaptive bases and sampling approach in
  \cite{peherstorfer2018model} and \cite{constantine2012reduced}.  Transformed
  snapshot interpolation method was also developed in
  \cite{welper2020transformed} to capture a sharp gradient in the solution, by
  introducing a new transform discretization near singularities.

  Even though all the approaches mentioned above do show some remedies of
  overcoming the challenges that arise from the advection-dominated problems,
  the solution representability of the linear subspace is still limited in a
  sense that the treatments introduced in the methods above are problem-specific
  and require some a priori knowledge, such as advection direction.  In order to
  maximize the representability and make the methodology as general as possible,
  it seems unavoidable to transition from the linear subspace to a nonlinear
  manifold solution representation. 

  There are many works available in the current literautre that looked into the
  nonlinear manifold solution represenation in physical simulations. Many of
  them treat the weights and biases of a neural network (NN) to be unknowns in
  the solution process.  For example, Lagaris, et al., used a single output NN
  as an argument for trial functions and minimized the partial/ordinary
  differential equation (PDE/ODE) residual norm \cite{lagaris1998artificial},
  where the weights of the NNs are used as optimization variables.  Dissanayake
  and Phan-Thien used the universal approximator of NNs as a solution
  representation for solving PDEs. They also used the weights of the NNs as
  parameters as in the work by Lagaris, et al.  \cite{dissanayake1994neural}. A
  similar method was also applied to a plasma equilibrium solver
  \cite{van1995neural}.  Meade and Fernandez used hard limit transfer functions
  for linear ordinary differential equations \cite{meade1994numerical}. However,
  these approaches can introduce too many unknowns because all the wieghts and
  biases need to be found during the PDE/ODE solution process.

  Recently, similar attempts have been made to incorporate physical laws into
  NN-based surrogate models —- so called physics-informed surrogate models,
  where the weights and biases of the NN are determined in the training phase.
  Such models include, but are not limited to, attempts to mimic temporal
  evolution by incorporating a time integrator in a loss function
  \cite{raissi2019physics, chen2018neural, khoo2017solving, long2018pde,
  beck2019machine} and to represent the solution with a trained NN
  \cite{raissi2019physics, zhu2019physics, berg2018unified}, the deep Galerkin
  method \cite{sirignano2018dgm}, approximating spatial gradient functions with
  a multilayer feedforward NN \cite{han2018solving}, DeepONet
  \cite{lu2019deeponet}, DeepXDE \cite{lu2019deepxde}, fractional
  physics-informed NNs (fPINNs) \cite{pang2019fpinns}, PINNs with uncertainty
  quantification \cite{zhang2019quantifying}, and Deep Ritz method of minimizing
  the energy functional with trial functions of NNs \cite{weinan2018deep,
  he2018relu}.  However, inclusions of NNs in the governing equations of the
  underlying physical laws, such as those above, do not take advantage of the
  existing numerical methods for high-fidelity physical simulations.

  Recently, a neural network-based ROM is developed in \cite{lee2020model},
  where the weights and biases are determined in the training phase and the
  existing numerical methods are utilized in their models. The same technique is
  extended to preserve the conserved quantities in the physical conservation
  laws \cite{lee2019deep}.  However, their approaches do not achieve any
  speed-up with respect to the corresponding FOM because the nonlinear terms that still scale with the FOM size need to be updated every time step or
  Newton step. 

  Two interesting papers were written by Rim, et al., recently. First of all,
  manifold approximations via transported subspaces in \cite{rim2019manifold}
  introduced a nonlinear solution representation by explicitly composing global
  transport dynamics with locally linear approximations of the solution
  manifolds.  However, their approach is only applicable to 1D problem for now. 
  The other work by Rim, et al., is the depth separation for reduced deep
  networks in nonlinear model reduction \cite{rim2020depth}, where they applied
  a compression technique on weight matrices and bias vectors to achieve the
  reduced deep networks. 

  We present a fast and accurate physics-informed neural network ROM with a
  nonlinear manifold solution representation, i.e., the nonlinear manifold ROM
  (NM-ROM). We train a shallow masked autoencoder with solution data from the
  corresponding FOM simulations and use the decoder as the nonlinear manifold
  solution representation. Our NM-ROM is different from the aformentioned
  physics-informed neural networks in that we take advantage of the existing
  numerical methods of solving PDE/ODEs in our approach. Furthermore, our NM-ROM
  is different from the neural network-based ROM of \cite{lee2020model} in a
  sense that we use a shallow masked autoencoder, while they used a deep
  convolutional autoencoder. The choice of the shallow masked NN over the deep
  convolutional NN is determined by the efficiency of the hyper-reduction
  technique we have developed.

\subsection{Nomenclature}\label{sec:nomenclature}
\begin{figure}[!t]
  \centering
  \includegraphics[width=0.7\textwidth]{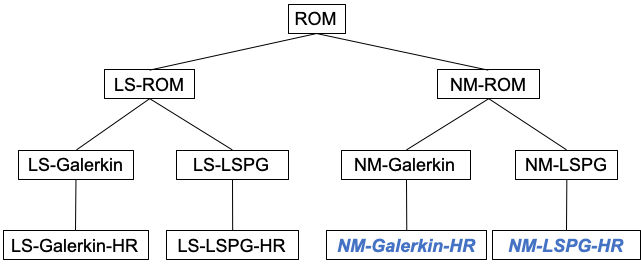}
  \caption{The figure shows the hierarchy of several ROMs. If the governing
  equation is nonlinear, then a hyper-reduction is required to achieve both
  accuracy and speed-up with respect to the corresponding FOM.  This paper
  contributes to the development of NM-LSPG-HR and NM-Galerkin-HR that achieve
  both speedup and accuracy with the NM-ROM. Throughout the paper, we will
  compare the performance of LS-ROMs and NM-ROMs.}
  \label{fg:rom_hierarchy}
\end{figure}

We use the following nomenclature/abbreviation for various ROMs throughout the paper:
\begin{itemize}
  \item FOM: full order model
  \item LS-ROM: linear subspace reduced order model
  \item LS-Galerkin: linear subspace Galerkin
  \item LS-LSPG: linear subpsace least-squares Petrov--Galerkin
  \item LS-Galerkin-HR: linear subspace Galerkin hyper-reduction
  \item LS-LSPG-HR: linear subspace least-squares Petrov--Galerkin hyper-reduction
  \item NM-ROM: nonlinear manifold reduced order model
  \item NM-Galerkin: nonlinear manifold Galerkin
  \item NM-LSPG: nonlinear manifold least-squares Petrov--Galerkin
  \item NM-Galerkin-HR: nonlinear manifold Galerkin hyper-reduction
  \item NM-LSPG-HR: nonlinear manifold least-squares Petrov--Galerkin hyper-reduction
\end{itemize}
These ROMs form a hierarchy that is depicted in Fig.~\ref{fg:rom_hierarchy}.
 
\subsection{Organization of the paper} 
We organize the subsequent sections by starting to discuss some background
materials in Section~\ref{sec:background}, where the FOM is stated in
Section~\ref{sec:FOM} and two LS-ROMs, i.e., LS-Galerkin and LS-LSPG, are
described in Sections~\ref{sec:LS-Galerkin} and \ref{sec:LS-LSPG},
respectively. Our NM-ROM is introduced in Section~\ref{sec:NM-ROM}, where the
nonlinear manifold solution representation is explained in
Section~\ref{sec:NM}. The shallow masked autoencoder that is used for the
solution representation is described in Section~\ref{sec:NN}.  The NM-Galerkin
is explained in Section~\ref{sec:NM-Galerkin} and the NM-LSPG is descirbed in
Section~\ref{sec:NM-LSPG}. The hyper-reduction technique that enables the
NM-ROM to achieve a speed-up is elaborated in Section~\ref{sec:HR}. The error
analysis is presented in Section~\ref{sec:error_analysis}. Finally, the
performance of our NM-ROM is demonstrated in two numerical experiments in
Section~\ref{sec:numericresults}. Finally, the paper is concluded with
summary and discussion in Section~\ref{sec:discussion-conclusion}.

\section{Background}\label{sec:background}
\subsection{Full order model}\label{sec:FOM}
A parameterized nonlinear dynamical system is considered, characterized by
a system of nonlinear ordinary differential equations (ODEs), which can be
considered as a resultant system from semi-discretization of Partial
Differential Equations (PDEs) in space domains
 \begin{equation} \label{eq:fom}
  \frac{d\sol}{dt} = \flux(\sol,t; \param),\quad\quad
  \sol(0;\param) = \solArg{0}(\param),
 \end{equation}     
where $t\in \timeDomain$ denotes time with the final time
$\totaltime\in\RRplus{}$, and $\sol(t;\param)$ denotes the time-dependent,
parameterized state implicitly defined as the solution to
problem~\eqref{eq:fom} with $\sol:\timeDomain\times \paramDomain\rightarrow
\RR{\nspacedof}$.  Further, $\flux: \RR{\nspacedof} \times \timeDomain
\times \paramDomain \rightarrow \RR{\nspacedof}$ with
$(\solDummy,\timeDummy;\paramDummy)\mapsto\flux(\solDummy,\timeDummy;\paramDummy)
$ denotes the velocity of $\sol$, which we assume to be nonlinear
in at least its first argument.  The initial state is denoted by
$\solArg{0}:\paramDomain\rightarrow \RR{\nspacedof}$, and $\param \in
\paramDomain$ denotes parameters in the domain
$\paramDomain\subseteq\RR{\nparam}$. 

A uniform time discretization is assumed throughout the paper, characterized
by time step $\dt\in\RRplus{}$ and time instances $\timeArg{n} =
\timeArg{n-1} + \dt$ for $n\in\nat{\ntimedof}$ with $\timeArg{0} = 0$,
$\ntimedof\in\natNo$, and $\nat{N}\defeq\{1,\ldots,N\}$.  To avoid
notational clutter, we introduce the following time discretization-related
notations: $\solArg{n} \defeq \solFuncArg{n}$, $\solapproxArg{n} \defeq
\solapproxFuncArg{n}$, $\redsolapproxArg{n} \defeq \redsolapproxFuncArg{n}$,
and $\fluxArg{n} \defeq \flux(\solFuncArg{n},t^{n}; \param)$, where
$\solapprox$, $\redsolapprox$ and $\redsolapproxFunc$ are defined in
Section~\ref{sec:LSROM}. 

Implicit time integrators are considered as time discretization methods. To illustrate this, we mainly consider the backward Euler time integrator for
an implicit scheme. Several other time integrators are shown in
Appendix~\ref{sec:appendixTimeIntegrators}. 

\begin{comment}
The explicit Forward Euler (FE) method numerically solves
Eq.~\eqref{eq:fom}, by time-marching with the following update: 
  \begin{equation} \label{eq:forwardEuler}
    \solArg{n} - \solArg{n-1} =\dt\fluxArg{n-1}.
  \end{equation}
Eq.~\eqref{eq:forwardEuler} implies the following subspace inclusion:
  \begin{equation}\label{eq:forwardEuler_partial_inclusion}
    \Span{\fluxArg{n-1}} \subseteq \Span{\solArg{n-1},\solArg{n}}.    
  \end{equation}
By induction, we conclude the following subspace inclusion relation:
  \begin{equation}\label{eq:forwardEuler_total_inclusion}
    \Span{\fluxArg{0},\dots,\fluxArg{\ntimedof-1}} \subseteq \Span{
      \solArg{0},\ldots, \solArg{\ntimedof}},
  \end{equation}
which shows that the span of nonlinear term snapshots is included in the
span of solution snapshots. The residual function with
the forward Euler time integrator is defined as
  \begin{align}\label{eq:residual_FE} 
  \begin{split}
    \resn_{\FE}(\solArg{n};\solArg{n-1},\param) &\defeq 
    \solArg{n} - \solArg{n-1} -\dt\fluxArg{n-1} 
  \end{split}
  \end{align} 
\end{comment}

The implicit Backward Euler (BE) method numerically solves Eq.~\eqref{eq:fom},
by solving the following nonlinear system of equations for $\solArg{n}$ at
$n$-th time step:
  \begin{equation} \label{eq:backwardEuler}
    \solArg{n} - \solArg{n-1} = \dt\fluxArg{n}.
  \end{equation} 
Eq.~\eqref{eq:backwardEuler} implies the following subspace inclusion:
  \begin{equation}\label{eq:backwardEuler_partial_inclusion}
    \Span{\fluxArg{n}} \subseteq \Span{\solArg{n-1},\solArg{n}}.    
  \end{equation}
By induction, we conclude the following subspace inclusion relation:
  \begin{equation}\label{eq:backwardEuler_total_inclusion}
    \Span{\fluxArg{1},\dots,\fluxArg{\ntimedof}} \subseteq \Span{\solArg{0},\ldots,
    \solArg{\ntimedof}},
  \end{equation}
which shows that the span of nonlinear term snapshots is included in the
span of solution snapshots. The residual function with
the backward Euler time integrator is defined as
  \begin{align}\label{eq:residual_BE} 
  \begin{split}
    \resn_{\BE}(\solArg{n};\solArg{n-1},\param) &\defeq 
    \solArg{n} - \solArg{n-1} -\dt\fluxArg{n}.
  \end{split}
  \end{align}

\subsection{Linear subspace reduced order model (LS-ROM)}\label{sec:LSROM}
Many projection-based reduced order models with linear subspace solution
representation can be considered for nonlinear dynamical systems. We consider
Galerkin and least-squares Petrov-Galerkin projection methods, which are the
most relevant to our proposed method, i.e., NM-ROM.

\subsubsection{Linear subspace solution representation}\label{sec:LS}
The linear subspace reduced order model approach applies spatial projection
using a subspace $\spatialSubspace \defeq
\Span{\basisvecspace_i}_{i=1}^\nbasisspace \subseteq \RR{\nspacedof}$ with
$\dim(\spatialSubspace)=\nbasisspace\ll\nspacedof$. Using this subspace, it
approximates the solution as $ \sol\approx\solapprox
\in\solArg{ref}+\spatialSubspace $ (i.e., in a trial subspace) or equivalently 
  \begin{equation}\label{eq:spatialLSROMsolution} 
    \sol \approx \solapprox= \solArg{ref} + \basismatspace\redsolapprox 
  \end{equation} 
and the time derivative of the solution as
  \begin{equation} \label{eq:spatialLSROMvelocity}
    \frac{d\sol}{dt}\approx\frac{d\solapprox}{dt} =
    \basismatspace\frac{d\redsolapprox}{dt}
  \end{equation}
where $\solArg{ref}\in\RR{\nspacedof}$ denotes a reference solution and
$\basismatspace \defeq [\basisvecspace_1 \cdots
\basisvecspace_{\nbasisspace}] \in\RR{\nspacedof\times\nbasisspace}$ denotes a
basis matrix and $\redsolapprox\in\RR{\nbasisspace}$  denotes the generalized
coordinates.  The initial condition for the generalized coordinate,
$\redsolapproxArg{0}\in\RR{\nbasisspace}$, is given by
$\redsolapproxArg{0}=\basismatspace^T\left(\solArg{0}-\solArg{ref}\right)$.

For constructing $\basismatspace$, Proper Orthogonal Decomposition (POD) is
commonly used.  POD \cite{berkooz1993proper} obtains
$\basismatspace$ from a truncated Singular Value Decomposition (SVD)
approximation to a FOM solution snapshot matrix. It is related to  principal
component analysis in statistical analysis \cite{hotelling1933analysis} and 
Karhunen--Lo\`{e}ve expansion \cite{loeve1955} in stochastic analysis.
POD forms a solution snapshot matrix,
$\snapshots\defeq\bmat{\solArg{0}^{\param_1} - \solArg{ref} &
\cdots &
\solArg{\ntimedof}^{\param_{\nparam}} -
\solArg{ref}}\in\RR{\ndof\times\nparam(\ntimedof+1)}$,
where $\solArg{n}^{\param_k}$ is a solution state at $n$-th time step with
parameter $\param_k$ for $n\innat{\ntimedof}$ and $k\innat{\nparam}$.
Then, POD computes its thin SVD: 
     \begin{equation}\label{eq:SVD} 
       \snapshots = \Ubold\Sigmabold\Vbold^T,
     \end{equation} 
where $\Ubold\in\RR{\ndof\times\nparam(\ntimedof+1)}$ and
$\Vbold\in\RR{\nparam(\ntimedof+1)\times\nparam(\ntimedof+1)}$ are orthogonal
matrices and $\Sigmabold\in\RR{\nparam(\ntimedof+1)\times\nparam(\ntimedof+1)}
$ is a diagonal matrix with singular values on its diagonals.  Then POD chooses
the leading $\nbasisspace$ columns of $\Ubold$ to set $\basismatspace$ (i.e.,
$\basismatspace = \bmat{\ubold_1 & \cdots & \ubold_{\nbasisspace}}$, where
$\ubold_k$ is $k$-th column vector of $\Ubold$).  The POD basis minimizes
$\|\snapshots - \basismatspace\basismatspace^T\snapshots \|_F^2$ over all
$\basismatspace\in\RR{\ndof\times\nbasisspace}$ with orthonormal columns, where
$\|\Abold\|_F$ denotes the Frobenius norm of a matrix $\Abold\in\RR{I\times J}$,
defined as $\|\Abold\|_F = \sqrt{\sum_{i=1}^{I}\sum_{j=1}^{J} a_{ij}^2}$ with
$a_{ij}$ being an $(i,j)$-th element of $\Abold$.  Since the objective function
does not change if $\basismatspace$ is post-multiplied by an arbitrary
$\nbasisspace\times\nbasisspace$ orthogonal matrix, the POD procedure seeks the
optimal $\nbasisspace$–-dimensional subspace that captures the snapshots in the
least-squares sense.  For more details on POD, we refer to
\cite{hinze2005proper,kunisch2002galerkin}.

\subsubsection{Linear subspace Galerkin projection}\label{sec:LS-Galerkin}
We derive LS-Galerkin using
time continuous residual minimization. First, we rewrite FOM ODE
Eq.~\eqref{eq:fom} as
\begin{equation}\label{eq:resFOM}
    \res(\vel,\sol,t;\param):= \vel - \flux(\sol,t; \param) = 0, \quad
    \sol(0;\param)=\solArg{0}(\param) 
\end{equation}
where $\res:\RR{\nspacedof}\times\RR{\nspacedof}\times\RRplus{}\times\paramDomain \rightarrow\RR{\nspacedof}$ with $(\dot{\solDummy},\solDummy,\timeDummy;\paramDummy) \mapsto \res(\dot{\solDummy},\solDummy,\timeDummy;\paramDummy)$ denotes the time continuous residual. Here, we denote $\dot{(\cdot)}$ as time derivative of $(\cdot)$ for notational simplicity. Replacing $\sol$ with $\solapprox$ given by Eq.~\eqref{eq:spatialLSROMsolution} and $\vel$ with $\velapprox$ given by Eq.~\eqref{eq:spatialLSROMvelocity} leads to the following residual function with the reduced number of unknowns
\begin{align}\label{eq:LSGalerkinresApprox}
  \resRedApproxArg{}(\redvel,\redsol,t;\param) :=
  \res(\basismatspace\redvel,\sol_{ref}+\basismatspace \redsol,t;\param),
\end{align}
where
$\resRedApproxArg{} :
\RR{\nbasisspace}\times\RR{\nbasisspace}\times\RRplus{}\times\paramDomain
\rightarrow\RR{\nspacedof}$ with $(\dot{\hat{\solDummy}},\hat{\solDummy},
\timeDummy;\paramDummy) \mapsto
\resRedApproxArg{}(\dot{\hat{\solDummy}},\hat{\solDummy},\timeDummy;\paramDummy)$ denotes the time
continuous residual.  Note that $\resRedApproxArg{}(\redvel,\redsol,t;\param) =
\zerobold$ is an over-determined system. Therefore, it is likely that no
solution exists. To close the system, we minimize the squared norm of the
residual vector function:
\begin{equation}\label{eq:LSGalerkinRes}
  \redvel = \argmin{\reddummy\in\RR{\nbasisspace}}
  \|\resRedApproxArg{}(\reddummy,\redsol,t;\param)\|_2^2
\end{equation}
with $\redsol(0;\param)=\redsolapproxArg{0}(\param)$
$=\basismatspace^T\left(\solArg{0}(\param)-\solArg{ref}\right)$. The solution to
Eq.~\eqref{eq:LSGalerkinRes} leads to the LS-Galerkin
\begin{equation}\label{eq:LSGalerkinResSol}
  \redvel=\basismatspace^{T}\flux(\sol_{ref}+\basismatspace \redsol,t;\param),
  \quad \redsol(0;\param)=\redsolapproxArg{0}(\param).
\end{equation}

% Replacing
% $\sol$ with $\solapprox$ in Eq.~\eqref{eq:fom} leads to the following system of
% equations with reduced number of unknowns:
%     \begin{equation} \label{eq:spatialLSROMvelocity}
%         \frac{d\solapprox}{dt}=\basismatspace\frac{d\redsolapprox}{dt} =
%         \flux(\solArg{ref}+\basismatspace\redsolapprox,t; \param).
%     \end{equation}
    
% Note that Eq.~\eqref{eq:spatialLSROMvelocity} has more equations than unknowns i.e., it is an
% over-determined system.  It is likely that there is no solution satisfying
% Eq.~\eqref{eq:spatialLSROMvelocity}.  In order to close the system, the Galerkin projection
% solves the following reduced system:
%     \begin{equation} \label{eq:lsrom_galerkin}
%         \frac{d\redsolapprox}{dt} =
%         \basismatspace^T\flux(\solArg{ref}+\basismatspace\redsolapprox,t;
%         \param).
%     \end{equation}

Applying a time integrator to Eq.~\eqref{eq:LSGalerkinResSol} leads to a fully
discretized reduced system, denoted as the reduced O$\Delta$E.  
Note that the reduced O$\Delta$E has $\nbasisspace$ unknowns and
$\nbasisspace$ equations. If an implicit time integrator is applied, a
Newton--type method can be applied to solve for unknown generalized coordinates
each time step. If an explicit time integrator is applied, time marching
updates will solve the system. However, we cannot expect any speed-up because
the size of the nonlinear term and its Jacobian, which need to be updated for
every Newton step, scales with the FOM size. In order to handle this issue, the
hyper-reduction will be applied (see Section~\ref{sec:LS-Galerkin-HR})

\subsubsection{Linear subspace least-squares Petrov--Galerkin projection}\label{sec:LS-LSPG}
    The Least-Squares Petrov--Galerkin (LSPG) method projects a fully
    discretized solution space onto a trial subspace. That is, it discretizes
    Eq.~\eqref{eq:fom} in time domain and replaces $\solArg{n}$ with
    $\solapproxArg{n} \defeq \solArg{ref}+\basismatspace\redsolapproxArg{n}$
    for $n\innat{\ntimedof}$ in residual functions defined in
    Section~\ref{sec:FOM} and Appendix \ref{sec:appendixTimeIntegrators}.  Here, we consider
    only implicit time integrators because the LSPG projection is equivalent to
    the Galerkin projection when an explicit time integrator is used as shown
    in Section 5.1 in \cite{carlberg2017galerkin}.  The residual functions for
    implicit time integrators are defined in \eqref{eq:residual_BE},
    \eqref{eq:residual_AM2}, and \eqref{eq:residual_BDF2} for various time
    integrators.  For example, the residual function with the backward Euler
    time integrator\footnote{\label{fn:timeIntegrators}Although the backward
    Euler time integrator is used extensively in the paper for illustrative purposes, many other time integrators introduced in
    Appendix~\ref{sec:appendixTimeIntegrators} can be applied to all the ROM methods dealt in
    the paper in a straight forward way.} after the trial subspace projection
    becomes 
    \begin{align}\label{eq:trialsub_residual_BE} 
    \begin{split}
      \resRedApproxArg{n}_{\BE}(\redsolapproxArg{n};\redsolapproxArg{n-1},\param)
&\defeq 
      \resn_{\BE}(\solArg{ref}+\basismatspace\redsolapproxArg{n};
      \solArg{ref}+\basismatspace\redsolapproxArg{n-1},\param) 
\\ &= \basismatspace(\redsolapproxArg{n} - \redsolapproxArg{n-1})
      -\dt\flux(\solArg{ref}+\basismatspace\redsolapproxArg{n},t_n;\param).
    \end{split}
    \end{align}
 The basis matrix $\basismatspace$ can be found by the POD as in the Galerkin
 approach.  Note that Eq.~\eqref{eq:trialsub_residual_BE} is an over-determined
 system.  To close the system and solve for the unknown generalized
 coordinates, $\redsolapproxArg{n}$, the LSPG takes the squared norm of the
 residual vector function and minimize it at every time step:
\begin{align} \label{eq:spOpt1}
  \begin{split} 
    \redsolapproxArg{n} = \argmin{\reddummy\in\RR{\nbasisspace}} \quad&
    \frac{1}{2} \left \|\resRedApproxArg{n}_{\BE}(\reddummy;\redsolapproxArg{n-1},\param)
    \right \|_2^2.
  \end{split} 
\end{align}
The Gauss--Newton method with the starting point $\redsolapproxArg{n-1}$ is
applied to solve the minimization problem~\eqref{eq:spOpt1} in LSPG.  However,
as in the Galerkin approach, a hyper-reduction, which will be discussed in
Section \ref{sec:LS-LSPG-HR}, is required for a speed-up due to the presence of
the nonlinear residual vector function that scales with the full order model
size.

\section{Nonlinear manifold reduced order model (NM-ROM)}\label{sec:NM-ROM}
A projection-based reduced order model with nonlinear manifold solution representation is introduced in this section.  The ROM formulation with
nonlinear manifold solution representation is introduced in
Section~\ref{sec:NM}. Section~\ref{sec:NN} describes how we construct the neural
network that is used as a nonlinear manifold solution representation.
As in the LS-ROMs of Section~\ref{sec:LSROM}, Galerkin and
least-squares Petrov--Galerkin projections will be applied in
Sections~\ref{sec:NM-Galerkin} and \ref{sec:NM-LSPG}. Finally, the
hyper-reduction for the NM-ROM is described in Section~\ref{sec:HR}.

\subsection{Nonlinear manifold solution representation}\label{sec:NM}
The NM-ROM applies solution representation using a nonlinear manifold
$\spatialSubspace \defeq \{\scaledDecoder\left(\reddummy\right)|\reddummy \in
\RR{\nbasisspace}\}$, where $\scaledDecoder: \RR{\nbasisspace} \rightarrow
\RR{\nspacedof}$ with $\nbasisspace\ll\nspacedof$ denotes a nonlinear function
that maps a latent space of dimension $\nbasisspace$ to the full order model
space of dimension, $\nspacedof$. That is, the NM-ROM approximates the solution
in a trial manifold as 
\begin{equation}\label{eq:spatialNMROMsolution} 
  \sol \approx \solapprox= \solArg{ref} + \scaledDecoder \left(\redsolapprox
  \right) 
\end{equation} 
and the time derivative of the solution as
\begin{equation}\label{eq:spatialNMROMvelocity}
  \frac{d\sol}{dt} \approx \frac{d\solapprox}{dt} =
  \jacobian_{g}\left(\redsolapprox\right)\frac{d\redsolapprox}{dt}
\end{equation}
where $\redsolapprox\in\RR{\nbasisspace}$ denotes the
generalized coordinates. The initial condition for the generalized coordinate,
$\redsolapproxArg{0}\in\RR{\nbasisspace}$, is given by
$\redsolapproxArg{0}=\scaledEncoder\left(\solArg{0}-\solArg{ref}\right)$, where
$\scaledEncoder \approx \scaledDecoder^{-1}$ (i.e., $\sol-\solArg{ref} \approx
\scaledDecoder\left(\scaledEncoder\left(\sol-\solArg{ref}\right)\right)$). The
details about the nonlinear functions, $\scaledEncoder$ and $\scaledDecoder$,
are presented in Section \ref{sec:NN}. 

\subsection{Shallow masked autoencoder}\label{sec:NN}
In this section, we present the approach for constructing a nonlinear manifold.
Here, we use an autoencoder, $\unscaledAE$, in the form of a feedforward neural
network, that is trained to reconstruct its input. The autoencoder architecture
is composed of an encoder, $\unscaledEncoder$ and a decoder, $\boldsymbol{D}$.
The encoder maps a high dimensional input, $\sol \in \RR{\nspacedof}$ to a
low-dimensional latent vector, $\redsolapprox \in \RR{\nbasisspace}$, i.e.,
$\unscaledEncoder(\sol)=\redsolapprox$, and the decoder then
maps the latent vector to $\solapprox \in \RR{\nspacedof}$, i.e.,
$\unscaledDecoder(\redsolapprox) = \solapprox$, where
$\nbasisspace \ll \nspacedof$. Therefore, we have
\begin{equation}\label{eq:AED}
    \sol \approx \solapprox = \unscaledAE(\sol) =
    \unscaledDecoder(\unscaledEncoder(\sol)).
\end{equation}

\begin{figure}[!htbp]
  \centering
  \includegraphics[width=0.5\textwidth]{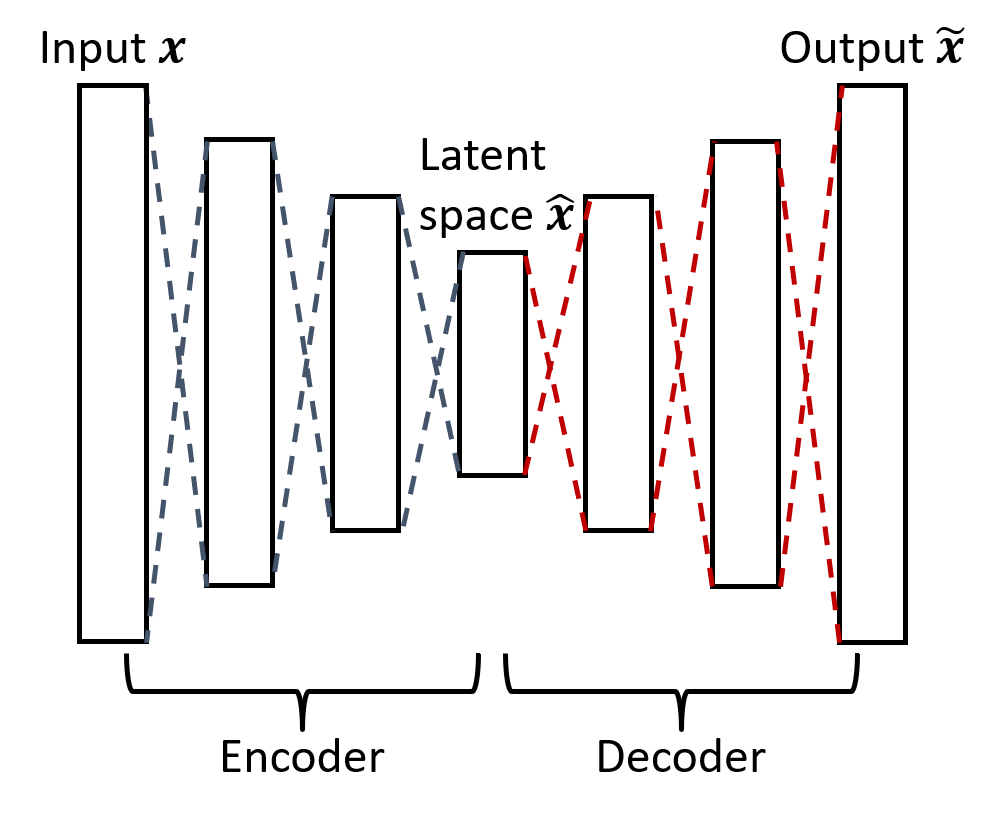}
  \caption{General description of an autoencoder: $\sol$ being encoded to a latent
  vector, $\redsolapprox$, by the encoder and decoded by the decoder, to
  $\solapprox$. The mean square error between $\sol$ and $\solapprox$ is
  minimized to update neural network weights and bias. }
  \label{fg:AE}
\end{figure}

The main idea behind an autoencoder is that it forces the model to learn salient
features by compressing the input into a low-dimensional space
and then reconstructing the input. 

\begin{figure}[!htbp]
  \centering
  \subfigure[Without masking]{
  \includegraphics[width=0.45\textwidth]{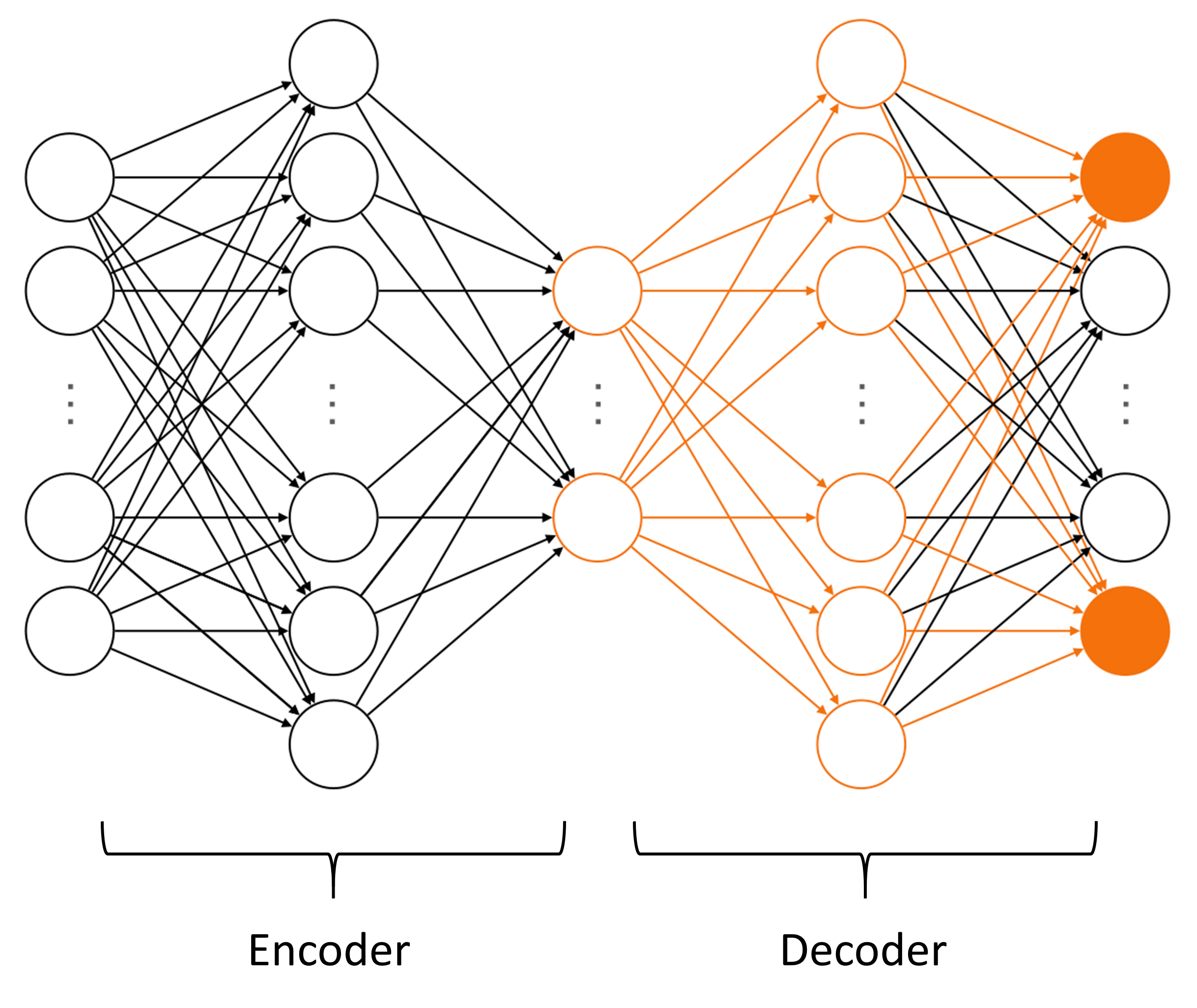}
  }
  \subfigure[With masking]{
  \includegraphics[width=0.45\textwidth]{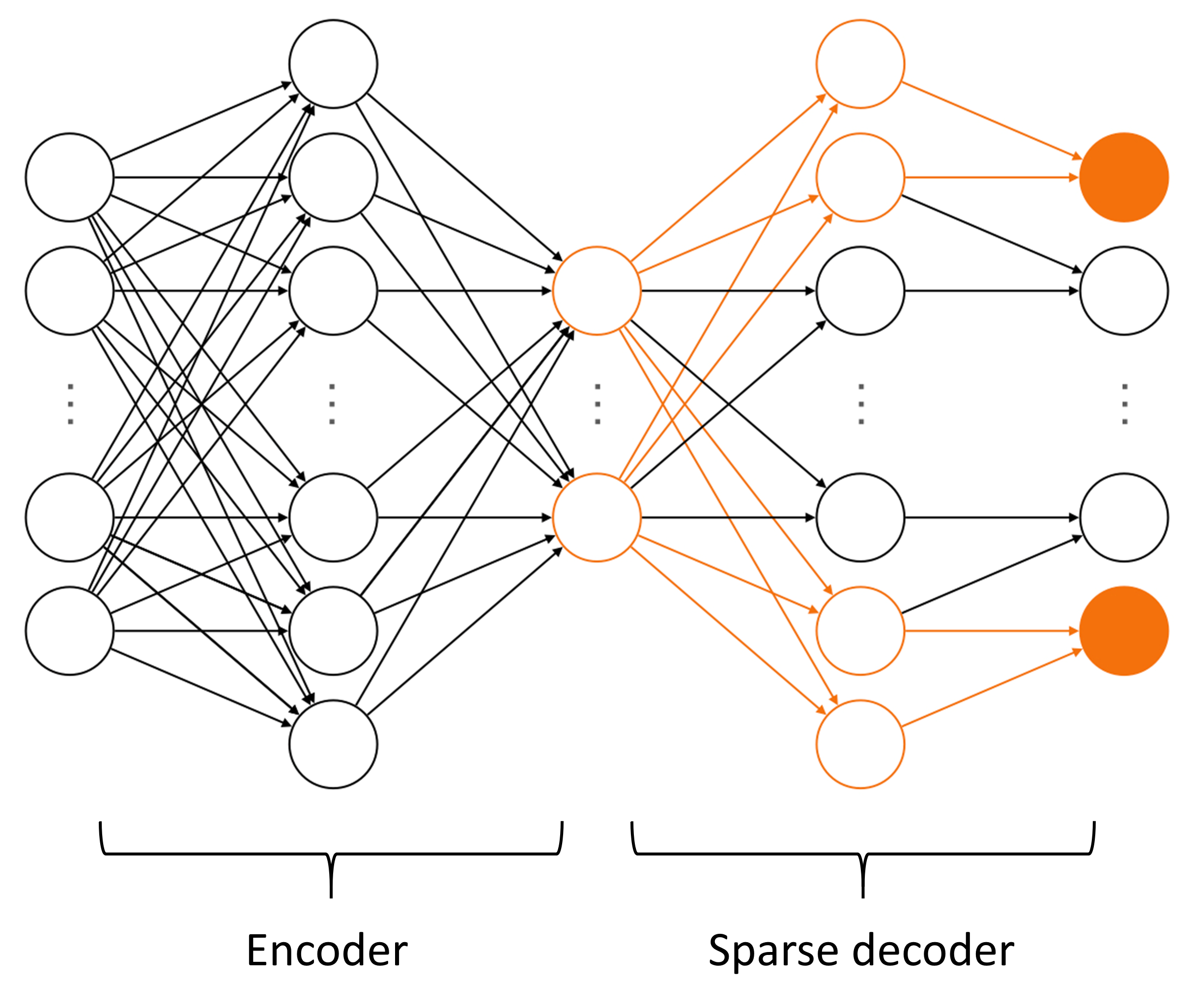}
  }
  \caption{Three layer autoencoder architecture: (a) unmasked and (b) masked
  shallow neural neural network. Nodes and edges in orange color represent
  active path that stems from the sampled outputs that are marked as the orange
  disks. Note that the masked shallow neural network has a sparser structure
  than the unmasked one.  }
  \label{fg:threeLayerAE}
\end{figure}

The universal approximation theorem
\cite{cybenko1989mathematics,pinkus1999approximation}, proves that
functions of the form,
\begin{equation}\label{eq:univeralapprox}
    v_k =
    \sum_{j=1}^{N_2}w_{jk2} \sigma
    \left(\sum_{i=1}^{N_1}w_{ij1}u_i+\theta_j\right) \hspace{12pt}
    \text{for}\hspace{2pt} \: k \innat{N_3},
\end{equation}
where $w_{ij1}, w_{jk2} \in \RR{}$ are weights, $\theta_j \in \RR{}$ is a bias, $\sigma$ is a
non-polynomial activation function, $u_i$ is an input and $v_k$ is an output,
can approximate any continuous, real-valued function arbitrarily well. Eq.
\eqref{eq:univeralapprox} is a simple, single hidden layer neural network
with a non-polynomial activation function. Its input dimension is $N_1$, width
of the hidden layer is $N_2$, and output dimension is $N_3$. We construct two single hidden layer neural networks, one is the 
encoder, $\unscaledEncoder$, and the other is the decoder, $\unscaledDecoder$. For
non-polynomial activation functions, a sigmoidal function given by
\begin{equation}\label{eq:sigmoid}
    \sigma(x) =  \frac{1}{1+\exp{(-x)}}
\end{equation}
or a swish function given by
\begin{equation}\label{eq:swish}
    \sigma(x) =  \frac{x}{1+\exp{(-x)}}
\end{equation}
are used. We use a non-deep neural network for the decoder because the decoder
and its Jacobian are computed many times during the ROM computation. In order
for this computation to be on par with POD methods, it is necessary to limit the
depth of the decoder network. The dimension of the encoder input and the decoder
output is $\nspacedof$ and the dimension of the encoder $\unscaledEncoder$
output and the decoder $\unscaledDecoder$ input is $\nbasisspace$. The width of
the encoder $\unscaledEncoder$ and decoder $\unscaledDecoder$ are
hyper-parameters. The first layers of the encoder $\unscaledEncoder$ and decoder
$\unscaledDecoder$ are fully-connected layers, where the nonlinear activation
functions are applied and the last layer of the encoder $\unscaledEncoder$ is fully-connected layer with no activation functions. The last layer of the decoder $\unscaledDecoder$ is either fully-connected layer or
sparsely-connected layer with no activation functions. These network
architectures are shown in Fig.~\ref{fg:threeLayerAE}.

Then, combining the encoder and the decoder yields the autoencoder which can be
trained to learn the identity mapping in an unsupervised manner because the 
desired output is the input. During the training phase, the error measured by 
\begin{equation}
    \|\snapshots-\Tilde{\snapshots}\|_F^2,
\end{equation}
where $\snapshots$ is solution snapshot matrix and $\Tilde{\snapshots}$ is a
reconstructed solution snapshot matrix, is minimized by optimizing learnable
parameters (i.e., weights and bias) in the two networks. The error is
back-propagated through the networks and the gradient with respect to the
learnable parameters are computed by using the chain rule
\cite{rumelhart1986learning,parker1985learnins,werbos1974beyond}. Then, the
parameters are updated in the steepest descent direction with respect to the
gradient.  Here, ADAM \cite{kingma2014adam}, a variant of stochastic gradient descent (SGD),
is used to approximate the gradient with a few data samples to make training
process faster. Stochastic gradient noise helps the neural network avoiding
over-fitting \cite{bottou2008tradeoffs}. Furthermore, graphics processing units
(GPUs) are utilized to parallelize the autonencoder's training by
simultaneously approximating multiple snapshots  \cite{raina2009large}.
In practice, a dataset is usually normalized before the training
process. 
Here, we normalize the dataset (i.e., solution snapshots) in the following way:
\begin{equation}\label{eq:normalize}
    \sol_{normal}=\sol_{scale}\odot \left(\sol - \sol_{ref} \right)
\end{equation}
where $\sol$ is a column vector of the dataset matrix $\snapshots$  and
$\odot$ denotes the element-wise product.  $\sol_{scale}$ and $\sol_{ref}$ are
directly computed from the dataset along each feature direction such that
$\sol_{normal}$ ranges either $[-1,1]$ or $[0,1]$.

After data normalization, an autoencoder can be trained to learn the identity
mapping with the normalized dataset.
Now, a normalized encoder maps from a high dimensional
normalized input $\sol_{normal}\in\RR{\ndof}$ to a low dimensional latent
vector $\redsolapprox\in\RR{\nbasisspace}$ in the form:
\begin{equation}\label{eq:en}
    \redsolapprox = \normalEncoder \left( \sol_{normal} \right)
\end{equation}
and a normalized decoder maps from the low dimensional latent vector $\redsolapprox\in\RR{\nbasisspace}$ to a
reconstructed normalized input $\solapprox_{normal}\in\RR{\ndof}$ in the form:
\begin{equation}\label{eq:de}
    \solapprox_{normal} = \mathbf{de} \left( \redsolapprox \right).
\end{equation}
Next, the encoder $\unscaledEncoder$ and the decoder $\unscaledDecoder$ can be
written by
\begin{align}\label{eq:ende}
    \unscaledEncoder\left(\sol\right) &= \mathbf{en} \left( \sol_{scale} \odot \left( \sol -\sol_{ref} \right) \right) \\
    \unscaledDecoder\left(\redsolapprox\right) &= \sol_{ref} + \mathbf{de}\left( \redsolapprox \right) \oslash \sol_{scale}
\end{align}
where $\odot$ and $\oslash$ denote the element-wise product and division,
respectively. Moreover, the row-wise product of $\sol_{scale}$ and the first
layer weight matrix of $\mathbf{en}$ yields the scaled encoder
$\scaledEncoder$.  Likewise, the row-wise division of $\sol_{scale}$ and the
last layer weight matrix of $\mathbf{de}$ gives us the scaled decoder
$\scaledDecoder$. 
Finally, the encoder $\unscaledEncoder$ and the decoder
$\unscaledDecoder$ are given by
\begin{align}\label{eq:hg}
    \unscaledEncoder\left(\sol\right)&=\boldsymbol{h}\left(\sol-\sol_{ref}\right) \\
    \unscaledDecoder\left(\redsolapprox\right)&=\sol_{ref}+\scaledDecoder\left(\redsolapprox\right).
\end{align}
We set the decoder
$\unscaledDecoder\left(\redsolapprox\right)=\sol_{ref}+\scaledDecoder\left(\redsolapprox\right)$
as the nonlinear manifold solution representation discussed in Section \ref{sec:NM-ROM}.

\begin{figure}[htbp]
  \centering
        \subfigure[Mask matrix for 1D Burgers equation]{
  \includegraphics[width=0.4\textwidth]{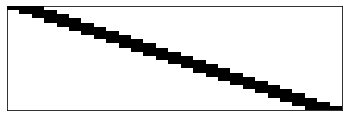}
  }
        \subfigure[Mask matrix for 2D Burgers equation]{
  \includegraphics[width=0.4\textwidth]{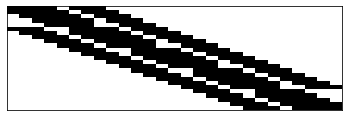}
  }
  \caption{Mask matrix. Note that the mask matrices have the analogical
  structure to the ones of Mass matrix that arises from a numerical
  discretization, such as the finite element or difference method, with 1D or 2D
  diffusion equations.}
  \label{fg:mask}
\end{figure}

The scaled decoder $\scaledDecoder$ can be written in the form
\begin{align}\label{eq:algebraicDecoder}
    \scaledDecoder\left(\redsolapprox\right)=\nnweight_2\activation\left(\nnweight_1\redsolapprox+\nnbias_1\right)
\end{align}
where $\nnweight_1$ and $\nnweight_2$ are weight matrices, $\nnbias_1$ is a
bias vector, and $\activation$ is an element-wise activation function. The
decoder can have more than two hidden layers (i.e., deep network). However, we
use the single layer decoder (i.e., shallow network) because the Jacobian
computation of the multiple hidden layer decoder involves multiple
matrix–-matrix multiplications. The output layer of the decoder
$\scaledDecoder$ is fully-connected as depicted in Fig.~\ref{fg:threeLayerAE}
(a) (i.e., $\nnweight_2$ is a dense matrix), which means all nodes in the
previous layer are required to compute even one element of the output vector. We apply a sparsity mask on the output layer of the decoder. Then, sampling a subset of the output vector doesn't need all nodes in the
previous layer as depicted in Fig.~\ref{fg:threeLayerAE} (b). Thus, more
speed-up can be achieved by a hyper-reduction technique that is described in
Section~\ref{sec:HR}. For example, the orange color nodes in
Fig.~\ref{fg:threeLayerAE} show the required nodes when the first and the last
elements of the output are selected, which are represented as solid orange
disks. To create a sparsely connected layer, we use a mask matrix $\mask$
which contains either zero or one as shown in Fig.~\ref{fg:mask}. By
element-wise product $\mask \odot \nnweight_2$, a sparse weight matrix is
obtained. The mask matrix $\mask$ is constructed to reflect local connectivity
as in the Laplacian operator approximated by the central difference scheme in
Finite Difference Method. The autoencoder composed of the encoder and the
sparse decoder is trained by using custom pruning in PyTorch \cite{paszke2019pytorch} pruning module.

In the autoencoder, the number of learnable parameters (i.e., weights and bias)
is determined by the number of nodes in the hidden layers in the encoder and
the decoder, dimension of latent vector, and the sparsity in the mask matrix.
The sparsity is determined by how many nodes in the hidden layer are used
to compute one element of the output and how many nodes in the hidden layer are
shared for neighboring elements of the output. To generate a mask matrix for 1D
problem, we use two variables $b$ and $\delta b$, where $b$ denotes the number
of nodes in the hidden layer to compute one output element (width of the block
in each row in Fig.~\ref{fg:mask}(a)) and $\delta b$ denotes the amount by
which the block shifts. For example, at the $i$th row,
$j\in\{(i-1)db,(i-1)db+1,\cdots,(i-1)db+b\}$th column is one and the others are
zero. For a mask of the 2D problem, we create a building matrix in the same way
as the mask matrix for 1D problem. Then, we add all rows neighboring $i$th row
(e.g., 5-point stencil for 2D and 7-point stencil for 3D) to $i$th row and
change nonzero values to one.  Note that the mask matrix for 2D problem as in
Fig.~\ref{fg:mask}(b) looks similar to 2D finite difference Laplacian operator.

There is no way to determine these hyper-parameters \textit{a priori}. If the
number of learnable parameters is not enough, the decoder is not able to
represent the nonlinear manifold well.  On the other hand, too many learnable
parameters may result in over-fitting, so the decoder is not able to generalize
well, which means the trained decoder can't be used for problems whose data is
unseen, i.e., the predictive case. To avoid over-fitting, there are two
options to consider. In the first option, one first divides the data into two
sets, i.e., train and test sets.  Then, the autoencoder is trained using the
train set only and is tested for the generalization ability using the test set.
If the mean squared error on the test and train sets are very
different, the over-fitting occurs and we should reduce the number of learnable
parameters \cite{kramer1991nonlinear}. 

The second option of avoiding the overfitting is to use Akaike's information
criteria (AIC) which is given by
\begin{align}\label{eq:aic}
    \text{AIC}=ln(e)+2\frac{N_w}{N}
\end{align}
where $e=\frac{\|\snapshots-\Tilde{\snapshots}\|_F^2}{2N}$, $N_w$ is the total
number of learnable parameters, and $N$ is the number of elements in the data
set matrix (i.e., $\snapshots$). If one minimizes only the first term of AIC,
then an over-fit network will be obtained. On the other hand, if one
minimizes only the second term of AIC, i.e., $N_w = 0$, then the network will not fit the training distribution.  Therefore, the minimum of AIC helps train a model that is not over-fit and generalizes well. 
\cite{ljung1999system,kramer1991nonlinear}. However, finding the minimum of
AIC requires a lot more training processes than the first option above. Because
of randomness in training, e.g., the random initialization of weights and
bias in neural networks and SGD optimization method, $e$ will be different for
every training process even with the same $N_w$ and the dataset matrix.
Therefore, AIC needs to be averaged over several training for each $N_w$ to
find the minimum of AIC. 

\begin{figure}[!htbp]
  \centering
  \subfigure[1D Burgers equation]{
  \includegraphics[width=0.30\textwidth]{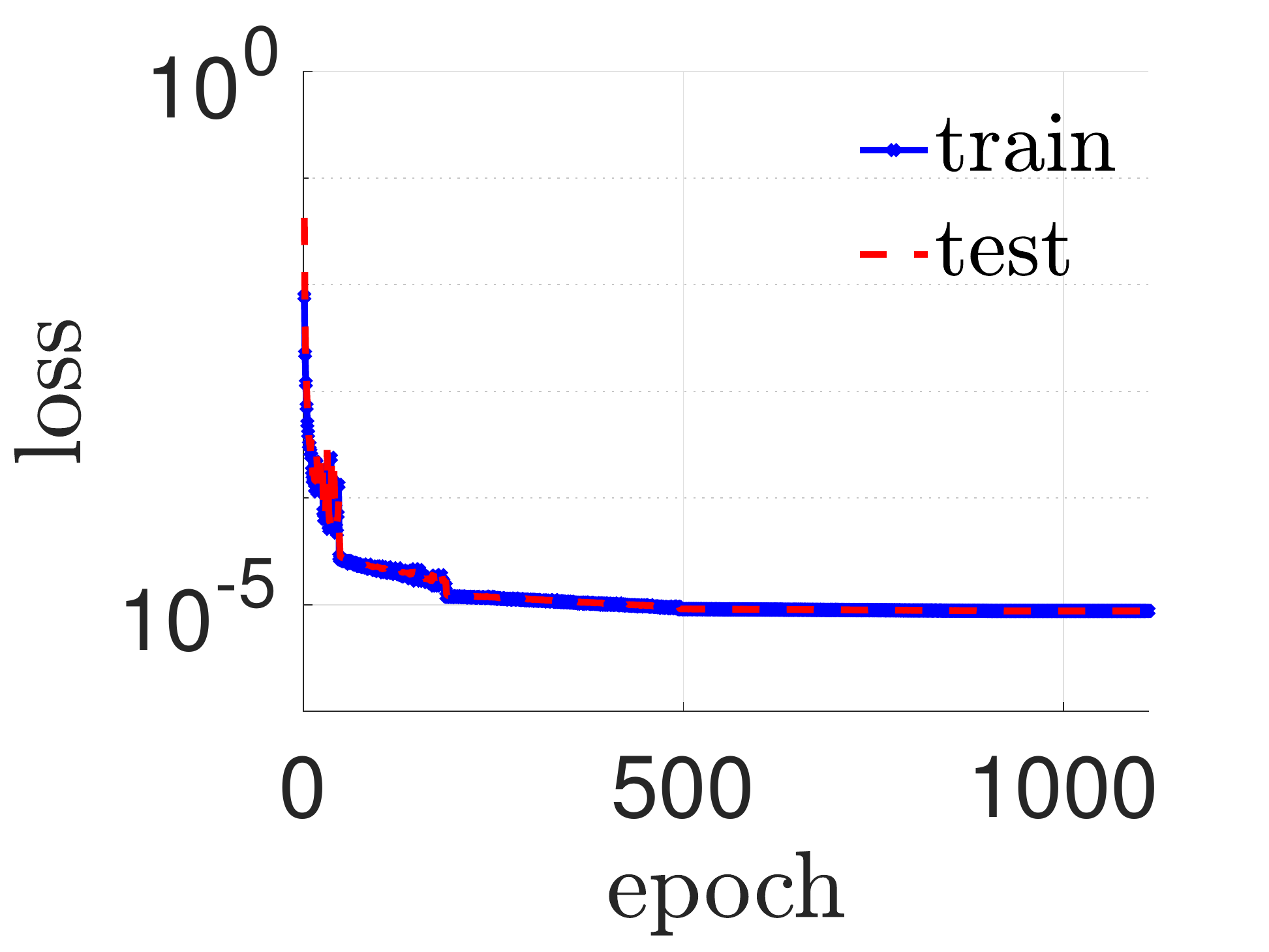}
  }~~~~~~~
  \subfigure[2D Burgers equation, $u$]{
  \includegraphics[width=0.30\textwidth]{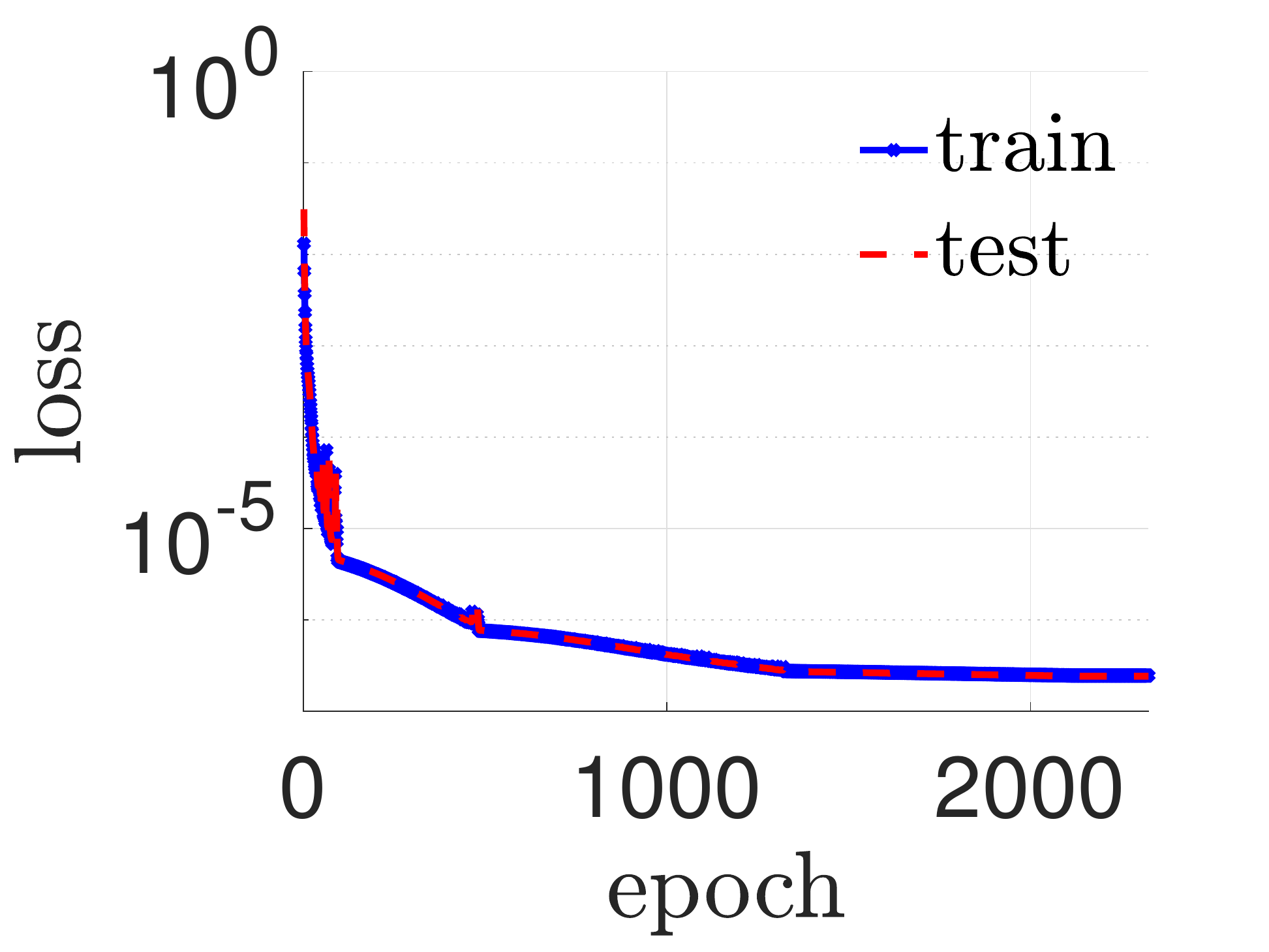}
  }~~~~~~~
  \subfigure[2D Burgers equation, $v$]{
  \includegraphics[width=0.30\textwidth]{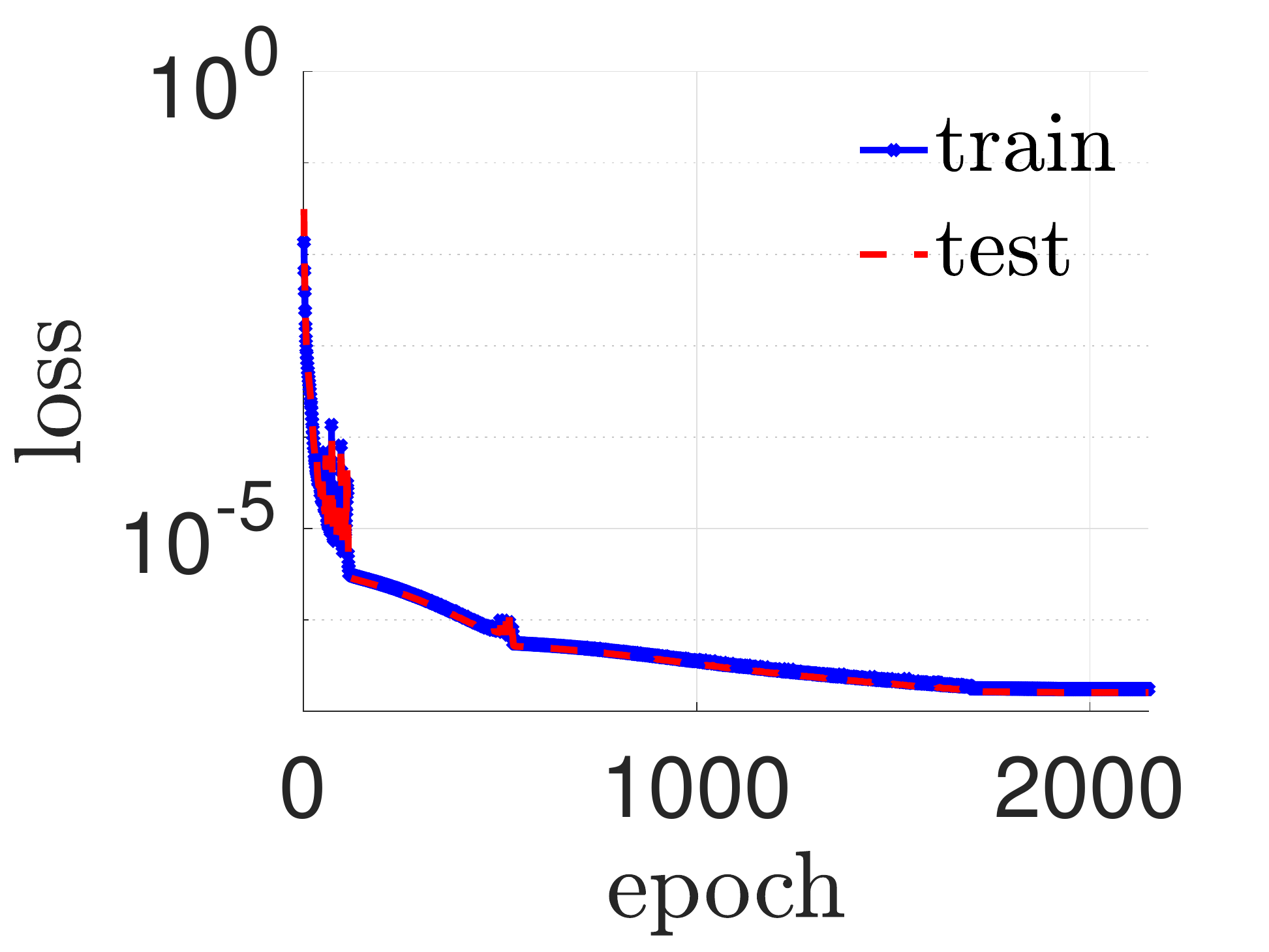}
  } 
  \caption{Loss history of decoders for various problems; all three figures show
  good agreement between train and test loss history, which is a sign for good
  balance between overfitting and accuracy.}
  \label{fg:lossHistory}
\end{figure}

Because of the practicality of the first option of avoiding the overfitting
over the second option, we use the first option in our numerical experiments.
For example, as shown in Fig.~\ref{fg:lossHistory}, the mean squares error on
the test and train data sets are very close. This implies that the trained
autoencoder is not over-fit.

\subsection{Nonlinear manifold Galerkin projection}\label{sec:NM-Galerkin}
We derive NM-Galerkin using time continuous residual minimization. Replacing $\sol$ with $\solapprox$ given by Eq.~\eqref{eq:spatialNMROMsolution} and $\vel$ with $\velapprox$ given by Eq.~\eqref{eq:spatialNMROMvelocity} in Eq.~\eqref{eq:resFOM} leads to the following residual function with the reduced number of unknowns
\begin{equation}\label{eq:NMGalerkinresApprox}
    \resRedApproxArg{}(\redvel,\redsol,t;\param):=\res(\jacobian_{g}(\redsol)\redvel,\sol_{ref}+\scaledDecoder(\redsol),t;\param).
\end{equation}
Note that Eq.~\eqref{eq:NMGalerkinresApprox} is an over-determined system. Therefore, it is likely that no solution exists. To close the system, we minimize the squared norm of the residual vector function:
\begin{equation}\label{eq:resNMGalerkin}
    \redvel = \argmin{\reddummy\in\RR{\nbasisspace}}\|\resRedApproxArg{}(\reddummy,\redsol,t;\param)\|_2^2
\end{equation}
with $\redsol(0;\param)=\redsolapproxArg{0}(\param)$
$=\scaledEncoder\left(\solArg{0}(\param)-\solArg{ref}\right)$. The solution to
Eq.~\eqref{eq:resNMGalerkin} leads to the NM-Galerkin
\begin{equation}\label{eq:resNMGalerkinSol}
    \redvel=\jacobian_{g}(\redsol)^{\dagger}\flux(\sol_{ref}+\scaledDecoder(\redsol),t;\param), \quad \redsol(0;\param)=\redsolapproxArg{0}(\param)
\end{equation}
where the Moore--Penrose inverse of a matrix $\weightmat \in \RR{\nspacedof\times \nbasisspace}$ with full column rank is defined as $\weightmat^{\dagger}:= (\weightmat^T\weightmat)^{-1}\weightmat^T$.

Applying a time integrator to Eq.~\eqref{eq:resNMGalerkinSol} leads to a fully
discretized reduced system, denoted as the reduced O$\Delta$E.  
Note that the reduced O$\Delta$E has $\nbasisspace$ unknowns and
$\nbasisspace$ equations. If an implicit time integrator is applied, a
Newton--type method can be applied to solve for unknown generalized coordinates
each time step. 
If an explicit time integrator is applied, time marching
updates will solve the system. 
However, we cannot expect any speed-up because the size of the
nonlinear terms and their Jacobians, which need to be updated for every
Newton step, scales with the FOM size. In order to handle this issue, the
hyper-reduction will be applied (see Section~\ref{sec:NM-Galerkin-HR}).

\subsection{Nonlinear manifold least-squares Petrov--Galerkin projection}\label{sec:NM-LSPG}
Alternatively, the nonlinear manifold least-squares Petrov--Galerkin (NM-LSPG) approach projects a fully discretized
solution space onto a trial manifold. That is, it discretizes Eq.~\eqref{eq:fom}
in time domain and replaces $\solArg{n}$ with $\solapproxArg{n} \defeq
\solArg{ref}+\scaledDecoder\left(\redsolapproxArg{n}\right)$ for
$n\innat{\ntimedof}$ in residual functions defined in Section~\ref{sec:FOM} and
Appendix \ref{sec:appendixTimeIntegrators}.  Here, we
consider only implicit time integrators for simplicity. See Ref.
\cite{leveque2007finite} for other types of time integrators. The residual
functions for several implicit time integrators are defined in \eqref{eq:residual_BE},
\eqref{eq:residual_AM2}, and \eqref{eq:residual_BDF2}. For example, the
residual function with the backward Euler time
integrator\footref{fn:timeIntegrators} after the trial manifold projection
becomes 
\begin{align}\label{eq:trialManifold_residual_BE} 
\begin{split}
 \resRedApproxArg{n}_{\BE}(\redsolapproxArg{n};\redsolapproxArg{n-1},\param)
  &\defeq
  \resn_{\BE}(\solArg{ref}+\scaledDecoder\left(\redsolapproxArg{n}\right);
  \solArg{ref}+\scaledDecoder\left(\redsolapproxArg{n-1}\right),\param) \\ &=
  \scaledDecoder\left(\redsolapproxArg{n}\right)-\scaledDecoder\left(\redsolapproxArg{n-1}\right)
  -\dt\flux(\solArg{ref}+\scaledDecoder\left(\redsolapproxArg{n}\right),t_n;\param).
\end{split}
\end{align}
 The nonlinear manifold $\scaledDecoder$ can be found by training the
 autoencoder as described in Section \ref{sec:NN}. Note that
 Eq.~\eqref{eq:trialManifold_residual_BE} is an over-determined system.
 Therefore, it is likely that no solution exists. To
 close the system and solve for the unknown generalized coordinates,
 $\redsolapproxArg{n}$, the NM-LSPG takes the squared norm of the residual
 vector function and minimizes it at every time step:
\begin{align} \label{eq:manifoldOpt1}
  \begin{split} 
    \redsolapproxArg{n} = \argmin{\reddummy\in\RR{\nbasisspace}} \quad&
    \frac{1}{2} \left \|\resRedApproxArg{n}_{\BE}(\reddummy;\redsolapproxArg{n-1},\param)
    \right \|_2^2.
  \end{split} 
\end{align}
The Gauss--Newton method with the starting point $\redsolapproxArg{n-1}$ is
applied to solve the minimization problem~\eqref{eq:manifoldOpt1}. However, as
in the Galerkin approach, a hyper-reduction which will be discussed in Section
\ref{sec:NM-LSPG-HR} is required for a speed-up due to the presence of the nonlinear
residual vector function that scales with the full order model size.
More specifically, $\scaledDecoder\left(\redsolapproxArg{n}\right)$,
$\flux(\solArg{ref}+\scaledDecoder\left(\redsolapproxArg{n}\right),t;\param)$,
and their Jacobians are needed to be updated whenever $\redsolapproxArg{n}$
chagnes if the backward Euler time integrator is used.

\section{Hyper-reduction}\label{sec:HR}
As mentioned in Section~\ref{sec:LSROM} and \ref{sec:NM-ROM}, we cannot expect
speed-up even though the dimension of unknowns in ROMs is small, i.e.,
$\nbasisspace \ll \nspacedof$, because the nonlinear term still scales with the
full order model size. To overcome this issue, there are several
hyper-reduction techniques available, e.g., \cite{chaturantabut2010nonlinear,
drmac2016new, drmac2018discrete, carlberg2013gnat, choi2020sns} for LS-ROMs.
These hyper-reduction techniques share a common feature and it plays an
important role in the development of the hyper-reduction technique in the
NM-ROMs, so we will go over one of the hyper-reduction technique that is
commonly used in the LS-ROMs.

\subsection{Nonlinear residual approximation}\label{sec:nonlinearTermApproximation}
We follow the DEIM-SNS and GNAT-SNS approaches introduced in \cite{choi2020sns}
where the solution snapshots, whose span includes a span of nonlinear term
snapshots, are taken to build a nonlinear term basis.  Then, it selects a
subset of each nonlinear term basis vector to either interpolate or data-fit in
a least-squares sense. In this way, it reduces the computational complexity of
updating nonlinear terms in an iterative solver for nonlinear problems.

In more details, the GNAT-SNS method approximates the nonlinear residual term
with gappy POD \cite{everson1995karhunen} as  
 \begin{equation}\label{eq:ResApprox}
   \resRedApproxArg{} \approx \basismatres\redres,
 \end{equation}
 where $\basismatres \defeq
 [\basisresvecArg{1},\ldots,\basisresvecArg{\nbasisres} ] \in
 \RR{\ndof\times\nbasisres}$, $\nbasisspace \leq \nbasisres \ll \ndof$, denotes
 the residual basis matrix and $\redres \in \RR{\nbasisres}$ denotes the
 generalized coordinates of the nonlinear residual term.  Here,
 $\resRedApproxArg{}$ represents a residual vector function, e.g., the backward
 Euler residual, $\resRedApproxArg{n}_{BE}$, defined in
 Eq.~\eqref{eq:trialsub_residual_BE}.  The GNAT-SNS method uses the SVD
 of the FOM solution snapshot matrix to construct $\basismatres$, which reduces
 computational cost by avoiding another POD to a
 nonlinear residual term snapshots. 
 The hyper-reduction method solves the following least-squares
 problem to obtain the generalized coordinates $\redres$: 
\begin{align} \label{eq:ResApprox_least-squares}
  \begin{split} 
    \redres := \argmin{\reddummy\in\RR{\nbasisres}} \quad&
    \frac{1}{2} \left \|\samplemat(\resRedApproxArg{} - \basismatres\reddummy)
    \right \|_2^2.
  \end{split} 
\end{align}
 where $\samplemat\defeq[\unitvecArg{p_1},\ldots,\unitvecArg{p_{\nressample}}]^T
 \in\RR{\nressample\times\ndof}$, $\nbasisspace \leq \nbasisres \leq \nressample
 \ll \ndof$, is the sampling matrix and $\unitvecArg{p_i}$ is
 the $p_i$th column of the identity matrix
 $\identity{\ndof}\in\RR{\ndof\times\ndof}$.  The solution to
 Eq.~\eqref{eq:ResApprox_least-squares} is given as
 \begin{equation}\label{eq:HR-generalizedcoordinates}
   \redres = (\samplemat\basismatres)^\dagger\samplemat\resRedApproxArg{},
 \end{equation}
 where the Moore--Penrose inverse of a matrix $\weightmat \in \RR{\nressample
 \times \nbasisres}$ with full column rank is defined as $\weightmat^{\dagger}
 := (\weightmat^T\weightmat)^{-1}\weightmat^T$. Therefore,
 Eq.~\eqref{eq:ResApprox} becomes
 \begin{equation}\label{eq:HR}
   \resRedApproxArg{} \approx \obliqueprojector \resRedApproxArg{},
 \end{equation}
 where $\obliqueprojector:= \basismatres
 (\samplemat\basismatres)^\dagger\samplemat$ is the oblique
 projection matrix. The projection matrix has a pseudo-inverse
 instead of the inverse because it allows the oversampling, i.e.,
 $\nbasisres<\nressample$.  The hyper-reduction method does not construct the
 sampling matrix $\samplematNT$.  {\it Instead, it maintains the sampling indices
 $\{p_1,\ldots,p_{\nbasisflux}\}$ and corresponding rows of $\basismatres$ and
 $\resRedApproxArg{}$.} This enables hyper-reduced ROMs to achieve a speed-up
 when it is applied to nonlinear problems.

 The sampling indices (i.e., $\samplematNT$) can be determined by Algorithm 3 of
 \cite{carlberg2013gnat} for computational fluid dynamics problems and Algorithm
 5 of \cite{carlberg2011efficient} for other problems. 
 These two algorithms take greedy procedure to minimize the error 
 in the gappy reconstruction of the POD basis vectors $\basismatres$. 
 These sampling algorithms for the hyper-reduction method allows
 oversampling (i.e., $\nressample > \nbasisres$), 
 resulting in solving least-squares problems in
 the greedy procedure. These selection algorithms can be viewed as the extension
 of Algorithm 1 in \cite{chaturantabut2010nonlinear} (i.e., a row pivoted LU
 decomposition) to the oversampling case. The nonlinear residual term projection
 error associated with these sampling algorithms is presented in Appendix D of
 \cite{carlberg2013gnat}. That is,
 \begin{equation}\label{eq:traditionalNMGalerkinHRprojection_error}
   \|\resRedApproxArg{} - \obliqueprojector\resRedApproxArg{}\|_2
   \leq \| \righttrianglemat^{-1} \|_2 
   \|\resRedApproxArg{} - \basismatres\basismatres^T\resRedApproxArg{} \|_2
 \end{equation}
 where $\righttrianglemat$ is the triangular factor from the QR factorization of
 $\samplemat\basismatres$ 
 (i.e.,
 $\samplemat\basismatres=\orthogonalmat\righttrianglemat$).
 For more details, please refer to
 \cite{choi2020sns}\footnote{In this paper, GNAT-SNS in \cite{choi2020sns} is
 re-named as LS-LSPG-HR to emphasize the difference between the LS-ROMs and
 NM-ROMs.}. 
 
\subsection{Hyper-reduction for LS-ROM}\label{sec:HRforLS-ROM}
We present formulations of LS-Galerkin-HR and LS-LSPG-HR. For numerical
examples, LS-LSPG-HR is only implemented.

\subsubsection{LS-Galerkin-HR}\label{sec:LS-Galerkin-HR}
We denote the hyper-reduced linear subspace Galerkin as LS-Galerkin-HR.  The
LS-Galerkin-HR method approximates the nonlinear residual term with the gappy
POD procedure as in Section~\ref{sec:nonlinearTermApproximation}. Therefore, the
LS-Galerkin-HR method replaces the residual in \eqref{eq:LSGalerkinRes} with
$\obliqueprojector \resRedApproxArg{}(\reddummy,\redsol,t;\param)$ given by
Eq.~\eqref{eq:HR}. Thus, it minimizes the following least-squares problem:
\begin{equation}\label{eq:resLSGalerkinHR}
    \redvel = \argmin{\reddummy\in\RR{\nbasisspace}}
    \|(\samplemat\basismatres)^\dagger \samplemat
    \resRedApproxArg{}(\reddummy,\redsol,t;\param)\|_2^2 
\end{equation}
with $\redsol(0;\param)=\redsolapproxArg{0}(\param)$. The solution to
Eq.~\eqref{eq:resLSGalerkinHR} leads to the following reduced ODE:
\begin{align}\label{eq:resLSGalerkinHRSol}
    \redvel =
    ((\samplemat\basismatres)^\dagger\samplemat\basismatspace)^{\dagger}
    ((\samplemat\basismatres)^T\samplemat\basismatres)^{-1}
    (\samplemat\basismatres)^T\samplemat
    \flux(\sol_{ref}+\basismatspace\redsol,t;\param), \quad
    \redsol(0;\param)=\redsolapproxArg{0}(\param).
\end{align}

Applying a time integrator to Eq.~\eqref{eq:resLSGalerkinHRSol} leads to a fully
discretized reduced system, denoted as the reduced O$\Delta$E.  
Note that the reduced O$\Delta$E has $\nbasisspace$ unknowns and
$\nbasisspace$ equations. If an implicit time integrator is applied, a
Newton--type method can be applied to solve for unknown generalized coordinates
each time step. If an explicit time integrator is applied, time marching
updates can be applied.

Note that the operator
$((\samplemat\basismatres)^\dagger\samplemat\basismatspace)^{\dagger}
((\samplemat\basismatres)^T\samplemat\basismatres)^{-1}
(\samplemat\basismatres)^T$ can be pre-computed once for all. We avoid
constructing the sampling matrix $\samplematNT$.  For example, the operator
$\samplemat\basismatres$ can be computed simply by extracting only the selected
rows of $\basismatres$. For the term, $\samplemat\flux$, only the nonlinear term
elements that are selected by the sampling matrix need to be computed. This
implies that we have to keep track of the rows of $\basismatspace$ that are
needed to compute the selected nonlinear term elements, which is usually a
larger set than the rows that are selected solely by the sampling matrix, i.e.,
$\samplemat\basismatspace$, as in the 5-point stencil or 7-point stencil in the
finite difference method. 

\subsubsection{LS-LSPG-HR}\label{sec:LS-LSPG-HR}
We denote the hyper-reduced linear subspace LSPG as LS-LSPG-HR.  The
LS-LSPG-HR method approximates the nonlinear residual term with the gappy
POD procedure as in Section~\ref{sec:nonlinearTermApproximation}. Therefore, the
LS-LSPG-HR method replaces the residual in \eqref{eq:spOpt1} with
$\obliqueprojector \resRedApproxArg{}(\reddummy,\redsol,t;\param)$ given by
Eq.~\eqref{eq:HR}. Thus, it minimizes the following least-squares problem:
  \begin{align} \label{eq:GNAT-spOpt1}
    \redsolapproxArg{n} = \argmin{\reddummy\in\RR{\nbasisspace}} \quad&
    \frac{1}{2} \left \|\ (\samplemat\basismatres)^\dagger\samplemat
    \resRedApproxArg{n}_{\BE}(\reddummy;\redsolapproxArg{n-1},\param) \right
    \|_2^2,
  \end{align}
with $\redsol(0;\param)=\redsolapproxArg{0}(\param)$. Note that the
pseudo-inverse $(\samplemat\basismatres)^\dagger$ can be pre-computed once for
all.  Due to the definition of $\resRedApproxArg{n}_{\BE}$ in
Eq.~\eqref{eq:trialsub_residual_BE}, the sampling matrix $\samplematNT$ needs to
be applied to the following terms: $\basismatspace(\redsolapproxArg{n} -
\redsolapproxArg{n-1})$ and
$\flux(\solArg{ref}+\basismatspace\redsolapproxArg{n},t;\param)$ at every time
step. The first term $\samplemat\basismatspace$ can be precomputed by extracting
the selected rows of the basis matrix. For the second term, only the nonlinear
term elements that are selected by the sampling matrix need to be computed. This
implies that we have to keep track of the rows of $\basismatspace$ that are
needed to compute the selected nonlinear term elements, which is usually a
larger set than the rows that are selected solely by the sampling matrix, i.e.,
$\samplemat\basismatspace$, as in the 5-point stencil or 7-point stencil in the
finite difference method.

\subsection{Hyper-reduction for NM-ROM}\label{sec:HRforNM-ROM}
There are two layers of nonlinear terms in the NM-ROM: (i) the nonlinear term in
the original governing equations, i.e., $\flux$ in Eq.~\eqref{eq:fom}, and (ii)
the decoder, which is nonlinear function of the generalized coordinates, i.e.,
$\scaledDecoder$ in Eq.~\eqref{eq:spatialNMROMsolution} and appears in the
definition of residuals both for Galerkin and Petrov--Galerkin cases.  The first
layer nonlinear term can be treated in the same way as the LS-ROMs (see
Sections~\ref{sec:LS-Galerkin-HR} and \ref{sec:LS-LSPG-HR}). Now, it is the
second layer nonlinear term that requires a special attention.  For example, the
Jacobian of the decoder needs to be evaluated at every solver iteration.
Because the cost of computing the Jacobian scales with the number of learnable
parameters in the decoder, we cannot expect much speed-up.  
As we did in the hyper-reduction process of the LS-ROMs, we have to avoid
computing all the entries of the decoder or its Jacobian because they scale
with the full order model size. This will be achieved by constructing a subnet
that computes only the relevant outputs, which is discussed in
Section~\ref{sec:EHRDC}. First, we state the hyper-reduced NM-ROMs, i.e., the
NM-Galerkin-HR in Section~\ref{sec:NM-Galerkin-HR} and the NM-LSPG-HR in
Section~\ref{sec:NM-LSPG-HR}. At last, the flop count estimate comparison
between non-hyper-reduced and hyper-reduced models are shown at the end of
Section~\ref{sec:EHRDC} and their derivations are shown in
Appendix~\ref{sec:appendixComputationalCosts}. 

\subsubsection{NM-Galerkin-HR}\label{sec:NM-Galerkin-HR}
Now, we apply the hyper-reduction to the NM-Galerkin method. We denote the
hyper-reduced nonlinear manifold Galerkin as NM-Galerkin-HR.  The NM-Galerkin-HR
method approximates the nonlinear residual term with the gappy POD procedure as
in Section~\ref{sec:nonlinearTermApproximation}. Therefore, the NM-Galerkin-HR
method replaces the residual in \eqref{eq:resNMGalerkin} with $\obliqueprojector
\resRedApproxArg{}(\reddummy,\redsol,t;\param)$ given by Eq.~\eqref{eq:HR}.
Thus, it minimizes the following least-squares problem:
\begin{equation}\label{eq:resNMGalerkinHR}
  \redvel =
  \argmin{\reddummy\in\RR{\nbasisspace}}\|(\samplemat\basismatres)^\dagger
  \samplemat\resRedApproxArg{}(\reddummy,\redsol,t;\param)\|_2^2
\end{equation}
with $\redsol(0;\param)=\redsolapproxArg{0}(\param)$. The solution to
Eq.~\eqref{eq:resNMGalerkinHR} leads to the NM-Galerkin-HR
\begin{equation}\label{eq:resNMGalerkinHRSol}
    \redvel = ((\samplemat\basismatres)^\dagger
    \samplemat\jacobian_{g}(\redsol))^{\dagger}(\samplemat\basismatres)^\dagger
    \samplemat\flux(\sol_{ref}+\scaledDecoder(\redsol),t;\param), \quad
    \redsol(0;\param)=\redsolapproxArg{0}(\param).
\end{equation}

Applying a time integrator to Eq.~\eqref{eq:resNMGalerkinHRSol} leads to a fully
discretized reduced system, denoted as the reduced O$\Delta$E.  
Note that the reduced O$\Delta$E has $\nbasisspace$ unknowns and
$\nbasisspace$ equations. If an implicit time integrator is applied, a
Newton--type method can be applied to solve for unknown generalized coordinates
each time step. 
If an explicit time integrator is applied, time marching
updates will solve the system. 

Note that the pseudo inverse, $(\samplemat\basismatres)^\dagger$, can be
pre-computed once for all by extracting only the selected rows of
$\basismatres$.  However, the term, $\samplemat\jacobian_{g}(\redsol)$, cannot
be precomputed because $\jacobian_{g}$ needs to be updated every time $\redsol$
is updated.  Fortunately, we need to compute only the selected rows of
$\jacobian_{g}$ by the sampling matrix $\samplematNT$. Similarly, 
for the term, $\samplemat\flux$, only the nonlinear term elements that are
selected by the sampling matrix need to be computed. This implies that we have
to keep track of the outputs of $\scaledDecoder$ that are needed to compute the
selected nonlinear term elements, which is usually a larger set than the outputs that are selected solely by the sampling matrix, i.e.,
$\samplemat\scaledDecoder$, as in the 5-point stencil or 7-point stencil in the
finite difference method. 

\subsubsection{NM-LSPG-HR}\label{sec:NM-LSPG-HR}
We apply the hyper-reduction to the NM-LSPG method discussed in Section
\ref{sec:NM-LSPG}. The hyper-reduction procedure for the nonlinear residual
function after the trial manifold projection is the same as the one in Section
\ref{sec:nonlinearTermApproximation}, i.e., we replace the residual defined in
\eqref{eq:trialManifold_residual_BE} with
$\obliqueprojector\resRedApproxArg{n}_{BE}$ and plug it into the minimization
problem in Eq.~\eqref{eq:manifoldOpt1}. Then, the minimization problem becomes
  \begin{align} \label{eq:NM-LSPG-HR-spOpt1}
    \begin{split} 
      \redsolapproxArg{n} = \argmin{\reddummy\in\RR{\nbasisspace}} \quad&
      \frac{1}{2} \left \|\ (\samplemat\basismatres)^\dagger\samplemat
      \resRedApproxArg{n}_{\BE}(\reddummy;\redsolapproxArg{n-1},\param)
      \right \|_2^2.
    \end{split} 
  \end{align}
Note that the pseudo-inverse $(\samplemat\basismatres)^\dagger$ can be
pre-computed once for all.  Due to the definition of $\resRedApproxArg{n}_{\BE}$
in Eq.~\eqref{eq:trialManifold_residual_BE}, the sampling matrix $\samplematNT$
needs to be applied the following two terms:
$\scaledDecoder\left(\redsolapproxArg{n}\right) -
\scaledDecoder\left(\redsolapproxArg{n-1}\right)$ and
$\flux(\solArg{ref}+\scaledDecoder\left(\redsolapproxArg{n}\right),t;\param)$ at
every time step. The first term,
$\samplemat(\scaledDecoder\left(\redsolapproxArg{n}\right) -
\scaledDecoder\left(\redsolapproxArg{n-1}\right))$, requires to compute only the
selected outputs of the decoder. Furthermore, for the second term, only the
nonlinear term elements that are selected by the sampling matrix need to be
computed. This implies that we have to keep track of the outputs of $\scaledDecoder$ that are needed to compute the selected nonlinear term elements by the
sampling matrix, which is usually a larger set than the outputs that are
selected solely by the sampling matrix, i.e., $\samplemat\scaledDecoder$, as in the 5-point stencil or 7-point stencil in the finite difference method. Therefore,
we build a subnet that computes only the outputs of the decoder that is required
to compute the elements of the nonlinear term, $\flux$. Then, with the same
subnet, the outputs required for the first term,
$\samplemat(\scaledDecoder\left(\redsolapproxArg{n}\right) -
\scaledDecoder\left(\redsolapproxArg{n-1}\right))$, can be extracted from the
same subnet. The construction of the subnet is explained in
Section~\ref{sec:subnet}.

\subsection{Efficient Hyper-Reduction Decoder Computation}\label{sec:EHRDC}
\begin{figure}[!htbp]
    \centering
    \includegraphics[width=0.6\textwidth]{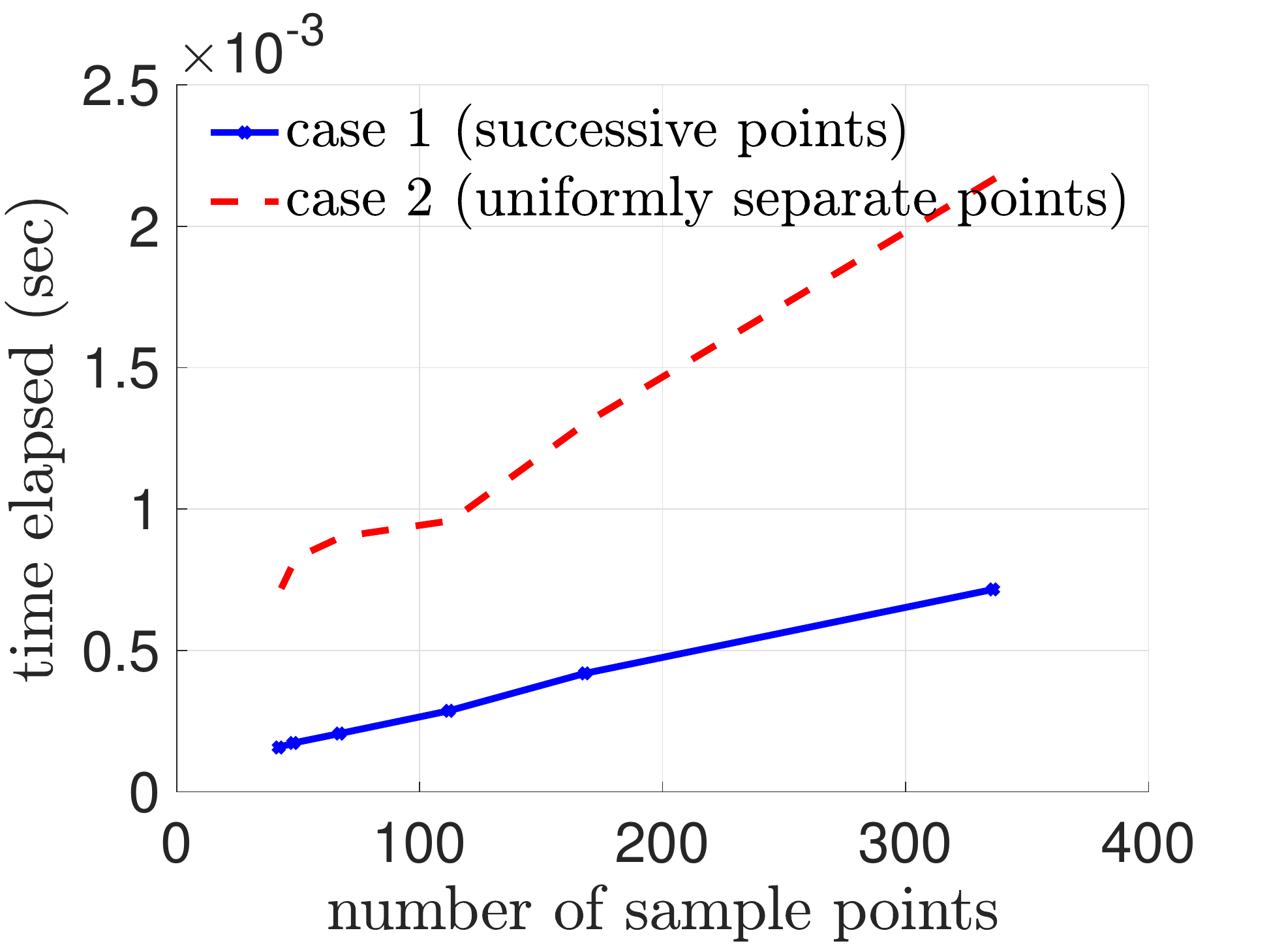}
    \caption{Illustration of the computational time elapsed for the evaluation
    of decoder and its Jacobian vs the number of sample points from 2D Burgers'
    equation in Section~\ref{sec:2dburgers}. The total number
    of points is $3364$. }
    \label{fg:DecoderJacobianTimeVSnsample}
\end{figure}

In the NM-LSPG-HR method, the residual is evaluated at the sampling points given
by the hyper-reduction. We use ``sample points" and ``hyper-reduction indices"
interchangeably throughout the paper.  Evaluating the decoder and its Jacobian
can be done efficiently by restricting the computation to the active paths of
the outputs required to compute the selected residual elements. For example,
active paths of the sparse decoder are shown in orange color in
Fig.~\ref{fg:threeLayerAE}(b). The costs of computing the decoder and its
Jacobian scale piecewise-linearly with the number of sample points because the
slopes of the costs of computing the decoder and its Jacobian vs the number of
sample points are different depending on how many nodes in hidden layer are
shared for each sample point (see Fig.~\ref{fg:DecoderJacobianTimeVSnsample}).
The distribution of the hyper-reduction indices determines the number of
overlapping nodes in hidden layer of decoder. The more the overlapping nodes in
hidden layer implies the more efficient computation of the hyper-reduced
decoder. If successive points are selected, overlapping of nodes in hidden layer
are maximized. If the selected points are uniformly  apart, then the overlapping
of nodes in hidden layer is minimized. In the case of random distribution, if
the number of selected points is small, the possibility of overlapping is low.
Our required outputs to compute the selected residual elements after the
hyper-reduction are neither successive nor uniformly separated.  Thus, the cost
of computing the decoder and its Jacobian would be between case 1 (successive
points) and case 2 (uniformly separated points) in
Fig.~\ref{fg:DecoderJacobianTimeVSnsample}. By restricting our computation to
active paths, we only compute along the subnet of the decoder network that is
needed for our required outputs.

\subsubsection{Construction of a subnet}\label{sec:subnet}
To determine the sparse decoder's active paths for given hyper-reduction
indices together with additional indices required to compute the hyper-reduced
residual (i.e., the neighbor indices that are used to approximate the
derivatives at the sample point), denoted as $\mathcal{H}$, we follow the steps
below:  
\begin{steps}
  \item Set nonlinear activation functions to be identity functions.  
  \item Replace nonzero elements of the weight matrices, $\nnweight_1,
    \nnweight_2 \odot \mask$ and the bias vector, $\nnbias_1$ with one and then
    denote each of them as $\tilde{\nnweight}_1$, $\tilde{\nnweight}_2$, and
    $\tilde{\nnbias}_1$, respectively. Zero elements of $\tilde{\nnweight}_1$,
    $\tilde{\nnweight}_2$, and $\tilde{\nnbias}_1$ represent non-connected edges
    between layers.
  \item A new decoder model, $\tilde{\scaledDecoder}(\nny_0)$, is created in the
    form
    \begin{equation}
      \tilde{\scaledDecoder}(\nny_0) =
      \tilde{\nnweight}_2(\tilde{\nnweight}_1\nny_0+\tilde{\nnbias}_1)
    \end{equation}
    or for each layer, we can write
    \begin{align}
      \nny_1&=\tilde{\nnweight}_1\nny_0+\tilde{\nnbias}_1 \\
      \nny_2&=\tilde{\nnweight}_2\nny_1.
    \end{align}
  \item Set $\nny_0=(1,\cdots,1)^T \in \RR{\nbasisspace}$ as an input. By
    construction, $\tilde{\scaledDecoder}(\nny_0)$ must be all positive.
  \item Define the target vector as $\nny_*=\nny_2-\mathbf{e} \in
    \RR{\nspacedof}$, where $i$th component of the error vector, $\mathbf{e}$,
    is defined as $e_i=\delta_{ij},j\in \mathcal{H}$.  Then the loss function,
    $L$, is defined as
    \begin{equation}
      L=\frac{1}{2}\|\nny_2-\nny_*\|_2^2
    \end{equation}
    and $\frac{\partial L}{\partial \nny_2}$ is given by
    \begin{equation}
      \frac{\partial L}{\partial \nny_2}=\nny_2-\nny_*=\mathbf{e}.
    \end{equation}
  \item Compute $\nabla_{\tilde{\nnweight}_2}L$,
    $\nabla_{\tilde{\nnweight}_1}L$, and $\nabla_{\tilde{\nnbias}_1}L$ using the
    chain rule
    \begin{align}
      \nabla_{\tilde{\nnweight}_2}L &= \frac{\partial L}{\partial
      \tilde{\nnweight}_2} = \frac{\partial L}{\partial \nny_2}\frac{\partial
      \nny_2}{\partial \tilde{\nnweight}_2} = \left(\frac{\partial L}{\partial
      \nny_2}\nny_1^T\right) \odot s(\tilde{\nnweight}_2) \\
      \nabla_{\tilde{\nnweight}_1}L &= \frac{\partial L}{\partial
      \tilde{\nnweight}_1} = \frac{\partial L}{\partial \nny_1}\frac{\partial
      \nny_1}{\partial \tilde{\nnweight}_1} = \frac{\partial L}{\partial
      \nny_2}\frac{\partial \nny_2}{\partial \nny_1}\frac{\partial
      \nny_1}{\partial \tilde{\nnweight}_1} =
      \left(\tilde{\nnweight}_2^T\frac{\partial L}{\partial \nny_2}\nny_0^T
      \right) \odot s(\tilde{\nnweight}_1) \\
      \nabla_{\tilde{\nnbias}_1}L &= \frac{\partial L}{\partial
      \tilde{\nnbias}_1} = \frac{\partial L}{\partial \nny_1}\frac{\partial
      \nny_1}{\partial \tilde{\nnbias}_1} = \frac{\partial L}{\partial
      \nny_2}\frac{\partial \nny_2}{\partial \nny_1}\frac{\partial
      \nny_1}{\partial \tilde{\nnbias}_1} =
      \left(\tilde{\nnweight}_2^T\frac{\partial L}{\partial \nny_2} \right)
      \odot s(\tilde{\nnbias}_1) 
    \end{align}
  where
  \begin{equation}
      s(x) := \begin{cases}
      0 \quad &if \quad x\leq0 \\
      1 \quad &otherwise
      \end{cases}
  \end{equation}
  is the element-wise function. Here, we make derivatives of $L$ with respect to
  non-connected edges (i.e., zero elements of $\tilde{\nnweight}_2$,
  $\tilde{\nnweight}_1$, and $\tilde{\nnbias}_1$) zeros by element-wise
  multiplication with $s(\tilde{\nnweight}_2)$, $s(\tilde{\nnweight}_1)$, and
  $s(\tilde{\nnbias}_1)$ because we do not consider non-connected edges as
  variables.
  \item Using the fact that the weights and bias that are not in the active
    paths do not contribute to computing $L$, we deduce that the derivatives of
    $L$ with respect to such weights and bias are zero. On the other hand, the
    derivatives of $L$ with respect to the weights and bias that are in the
    support of indices in $\mathcal{H}$ must be strictly positive because the special
    structure of $\tilde{\scaledDecoder}$ (i.e., the same structure as
    the sparse decoder, $\scaledDecoder$, except all the nonzero weights and
    bias are one and the nonlinear activation functions are identity functions), 
    choosing the all-ones vector as input vector, and defining the target vector
    as above should induce the positive gradient to reduce the $L$.  Thus,
    active path weights and bias are obtained by
    \begin{align}
      \nnweight_2^a &= (\nnweight_2 \odot \mask) \odot
      s(\nabla_{\tilde{\nnweight}_2}L) \\
      \nnweight_1^a &= \nnweight_1 \odot s(\nabla_{\tilde{\nnweight}_1}L) \\
      \nnbias_1^a &= \nnbias_1 \odot s(\nabla_{\tilde{\nnbias}_1}L).
    \end{align}
  \item Removing zero rows and zero columns of the active path weights and bias,
    $\nnweight_2^a$, $\nnweight_1^a$, and $\nnbias_1^a$ yields the subnet
    weights and bias, which are denoted as $\nnweight_2^{\subnetSymbol}$,
    $\nnweight_1^{\subnetSymbol}$, and $\nnbias_1^{\subnetSymbol}$,
    respectively. Then the subnet, $\scaledDecoder^{\subnetSymbol}$ is given by
  \begin{equation}
      \scaledDecoder^{\subnetSymbol}(\redsol) =
      \nnweight_2^{\subnetSymbol}\activation(\nnweight_1^{\subnetSymbol}\redsol
      + \nnbias_1^{\subnetSymbol}).
  \end{equation}
\end{steps}

\begin{figure}[!htbp]
    \centering
    \subfigure[Shallow network ]{
    \includegraphics[width=0.48\textwidth]{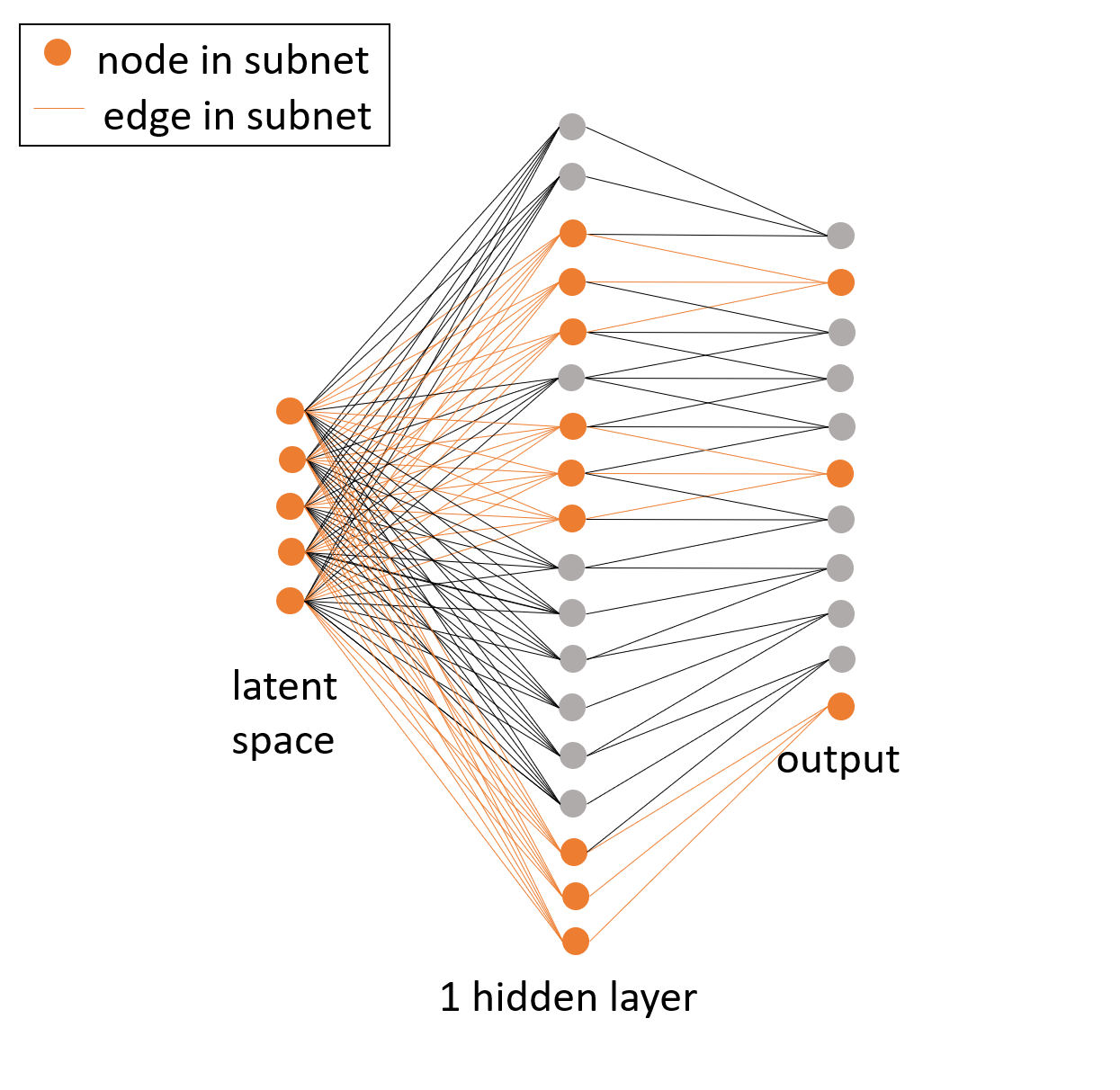}
    }
    \subfigure[Deep network]{
    \includegraphics[width=0.48\textwidth]{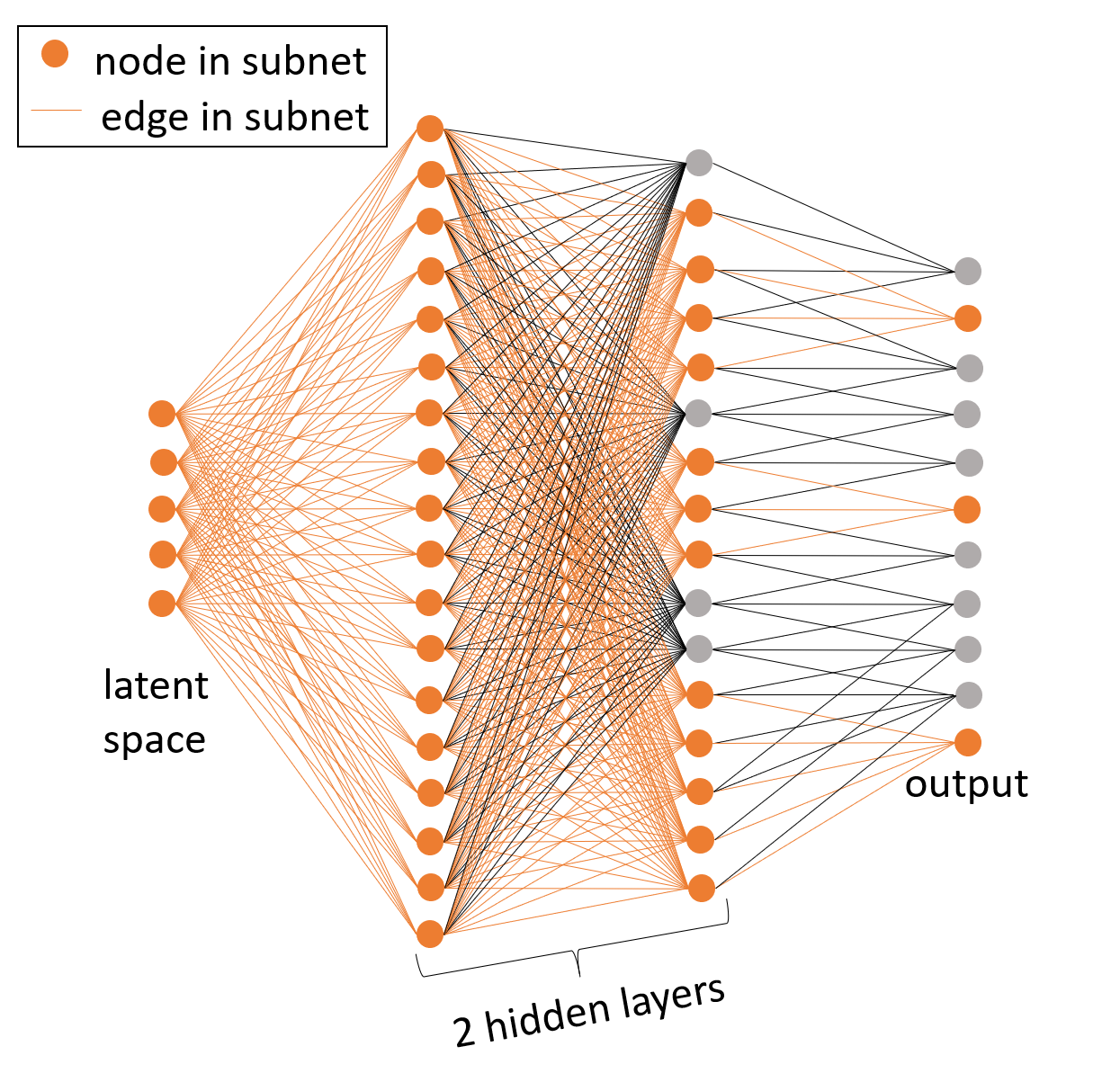}
    }
    \caption{Illustration of the effect on the sparsity of the active path for
    shallow network vs deep network. The shallow network provides a sparser
    network than the deep network in the subnet. Therefore, the shallow network
    is expected to achieve a higher speed-up than the deep network.}
    \label{fg:shallowVSdeep}
\end{figure}

This subnet strategy works for neural networks of arbitrary depth. However, we
have illustrated it in the context of the neural network with one hidden layer.
It is because that is what we use to achieve enough speed up. Please see
Fig.~\ref{fg:shallowVSdeep} for an argument of a shallow over a deep network. 

{\remark\label{remark:flop} To count flops of NM-LSPG and NM-LSPG-HR, let $m$,
$f$, $z$, and $b$ denote FOM size, ROM size, the size of subnet output, and the
number of nodes in the hidden layer to compute one output element of the sparse
decoder, respectively. Then, the flop counts of NM-LSPG is $\mathcal{O}(mbf)$
and the flop counts of NM-LSPG-HR is $\mathcal{O}(zbf)+\mathcal{O}(fz^2)$. Thus,
if $z$ is small e.g., $z^2 < m$, speed-up can be achieved. For full details, see
Appendix~\ref{sec:appendixComputationalCosts}.}

\section{Error analysis}\label{sec:error_analysis}
We present error analysis of the NM-LSPG-HR method. The error analysis is based
on \cite{lee2020model} and we added an oblique projection matrix for a
hyper-reduction method. \textit{A posteriori} discrete-time error bounds for
NM-Galerkin and NM-LSPG without any hyper-reduction are derived in
\cite{lee2020model}. Here, we apply a linear multi-step method described in
Appendix~\ref{sec:appendixTimeIntegrators}. 

{\theorem\label{theorem:errorbound} 
Let $\samplematNT \in \RR{\ndof \times \nressample}$ with $\nressample \ll
\ndof$ denote a sampling matrix, $\obliqueprojector \in \RR{\ndof \times \ndof}$
be an oblique projection matrix used in NM-Galerkin-HR and NM-LSPG-HR, i.e.,
$\obliqueprojector = \obliqueprojector_{\NMGalerkinHR}$ for NM-Galerkin-HR and
$\obliqueprojector = \obliqueprojector_{\NMLSPGHR}$ for NM-LSPG-HR, and
$\resRedApproxArg{n} \in \RR{\ndof}$ denote the nonlinear residual term,
which is defined by replacing $\solArg{n}$ with $\solapproxArg{n} \defeq
\solArg{ref}+\scaledDecoder\left(\redsolapproxArg{n}\right)$ for
$n\innat{\ntimedof}$ in residual functions defined in Section~\ref{sec:FOM} and
Appendix \ref{sec:appendixTimeIntegrators}, e.g., the residual with the
backward Euler time integrator is defined in Sections~\ref{sec:NM-LSPG} and
~\ref{sec:NM-LSPG-HR}.  Then, if the velocity $\flux$ is Lipschitz continuous
with the Lipschitz constant $L$ and the time step size $\Delta t$ is
sufficiently small such that $\Delta t < \frac{\gamma_1
|\alpha_0|}{\gamma_2|\beta_0|L}$, we have the following error bound 
\begin{equation}\label{eq:nmgalerkinerrorbound}
\begin{split}
    \| \solArg{n}-\solArg{ref}-\scaledDecoder(\redsolapproxArg{n}) \|_2 
    &\leq \frac{1}{\| \obliqueprojector \|_2 (\gamma_1-\gamma_2\frac{|\beta_0|\Delta t L}{|\alpha_0|})|\alpha_0|}\left \|\obliqueprojector \resRedApproxArg{n}(\solArg{ref}+\scaledDecoder(\redsolapproxArg{n}))
    \right \|_2 \\
    &+\sum_{j=1}^k \frac{|\alpha_j|+|\beta_j|\Delta t L}{(\gamma_1-\gamma_2\frac{|\beta_0|\Delta t L}{|\alpha_0|})|\alpha_0|}
    \left \|
    \solArg{n-j}-\solArg{ref}-\scaledDecoder(\redsolapproxArg{n-j})
    \right \|_2
\end{split}
\end{equation}
for NM-Galerkin-HR and
\begin{equation}\label{eq:nmlspgerrorbound}
\begin{split}
    \| \solArg{n}-\solArg{ref}-\scaledDecoder(\redsolapproxArg{n}) \|_2 
    &\leq \frac{1}{\| \obliqueprojector \|_2 (\gamma_1|\alpha_0|-\gamma_2|\beta_0|\Delta t L)} \min_{\reddummy\in\RR{\nbasisspace}}
    \left \|
    \obliqueprojector \resRedApproxArg{n}(\reddummy;\redsolapproxArg{n-1},\cdots,\redsolapproxArg{n-k},\param)
    \right \|_2 \\
    &+\sum_{j=1}^k \frac{|\alpha_j|+|\beta_j|\Delta t L}{(\gamma_1|\alpha_0|-\gamma_2|\beta_0|\Delta t L)}
    \left \|
    \solArg{n-j}-\solArg{ref}-\scaledDecoder(\redsolapproxArg{n-j})
    \right \|_2
\end{split}
\end{equation}
for NM-LSPG-HR, where coefficients $\alpha_j, \text{ }\beta_j \in \RR{}, \text{ }j=0, \cdots, k$ define a particular linear multi-step scheme and $\gamma_1, \gamma_2 \in \RR{}$ are $0 < \gamma_1, \gamma_2\leq 1$.

\proof 
We have
\begin{align}
    \resn(\solArg{n}) &= \sum_{j=0}^k \alpha_j \solArg{n-j} - \Delta t \sum_{j=0}^k \beta_j \flux(\solArg{n-j}) = 0 \label{eq:LMMres}, \\
    \obliqueprojector \resRedApproxArg{n}(\solArg{ref}+\scaledDecoder(\redsolapproxArg{n})) &= \obliqueprojector \left( \sum_{j=0}^k \alpha_j \left( \solArg{ref}+\scaledDecoder(\redsolapproxArg{n-j}) \right) - \Delta t \sum_{j=0}^k \beta_j \flux \left(\solArg{ref}+\scaledDecoder \left(\redsolapproxArg{n-j} \right)\right) \right) \label{eq:LMMresRedApprox}
\end{align}
where $\solArg{n} \in \RR{\ndof}$ denotes FOM solution and $\solArg{ref}+\scaledDecoder(\redsolapproxArg{n}), \text{ }\redsolapproxArg{n} \in \RR{\nbasisspace}$ is approximate solution.

Subtracting Eq.~\eqref{eq:LMMres} from Eq.~\eqref{eq:LMMresRedApprox} gives
\begin{equation}
    \begin{split}
    -\obliqueprojector \resRedApproxArg{n} \left( \solArg{ref} + \scaledDecoder(\redsolapproxArg{n}) \right)  &=  \obliqueprojector \left( \alpha_0 \left( \solArg{n}-\solArg{ref}-\scaledDecoder(\redsolapproxArg{n}) \right) -\Delta t \beta_0 \left( \flux(\solArg{n})-\flux(\solArg{ref}+\scaledDecoder(\redsolapproxArg{n})) \right) \vphantom{\sum_{j=1}^k} \right. \\
    \left. \right. &+ \left. \sum_{j=1}^k \alpha_j \left( \solArg{n-j} - \solArg{ref}-\scaledDecoder(\redsolapproxArg{n-j}) \right) - \Delta t \sum_{j=1}^k \beta_j \left( \flux \left(\solArg{n-j}\right) - \flux \left(\solArg{ref}+\scaledDecoder(\redsolapproxArg{n-j})\right) \right) \right).
    \end{split}
\end{equation}
We can re-write this in the following form
\begin{multline}
    \underbrace{
    \obliqueprojector \left( \solArg{n}-\solArg{ref}-\scaledDecoder(\redsolapproxArg{n}) -\frac{\beta_0 \Delta t}{\alpha_0} \left( \flux(\solArg{n})-\flux(\solArg{ref}+\scaledDecoder(\redsolapproxArg{n})) \right)\right)
    }_{\RN{1}} = \\
    \underbrace{
    -\frac{1}{\alpha_0}\obliqueprojector\resRedApproxArg{n}(\solArg{ref}+\scaledDecoder(\redsolapproxArg{n})) + \obliqueprojector \left( -\frac{1}{\alpha_0} \sum_{j=1}^k \alpha_j \left( \solArg{n-j} - \solArg{ref}-\scaledDecoder(\redsolapproxArg{n-j}) \right) + \frac{\Delta t}{\alpha_0} \sum_{j=1}^k \beta_j \left( \flux \left(\solArg{n-j}\right) - \flux \left(\solArg{ref}+\scaledDecoder(\redsolapproxArg{n-j})\right) \right) \right)
    }_{\RN{2}}
\end{multline}
Applying the reverse triangle inequality gives
\begin{equation}
    \| \RN{1} \|_2 \geq \left| \| \obliqueprojector (\solArg{n}-\solArg{ref}-\scaledDecoder(\redsolapproxArg{n})) \|_2 - \left\| \frac{\beta_0 \Delta t}{\alpha_0} \obliqueprojector ( \flux (\solArg{n}) - \flux (\solArg{ref}+\scaledDecoder(\redsolapproxArg{n}))) \right\|_2 \right|.
\end{equation}
Now, we use relationships
\begin{equation}
    \| \obliqueprojector (\solArg{n}-\solArg{ref}-\scaledDecoder(\redsolapproxArg{n})) \|_2 = \gamma_1 \|\obliqueprojector \|_2 \| \solArg{n}-\solArg{ref}-\scaledDecoder(\redsolapproxArg{n}) \|_2
\end{equation}
and
\begin{equation}
    \| \obliqueprojector (\flux (\solArg{n}) - \flux (\solArg{ref}+\scaledDecoder(\redsolapproxArg{n}))) \|_2 = \gamma_2 \|\obliqueprojector \|_2 \| \flux (\solArg{n}) - \flux (\solArg{ref}+\scaledDecoder(\redsolapproxArg{n})) \|_2
\end{equation}
where $0 < \gamma_1 \leq 1$ and $0 < \gamma_2 \leq 1$. Then, we have
\begin{equation}
    \begin{split}
        \| \RN{1} \|_2 &\geq \left| \gamma_1 \|\obliqueprojector \|_2 \| \solArg{n}-\solArg{ref}-\scaledDecoder(\redsolapproxArg{n}) \|_2 - \gamma_2 \|\obliqueprojector \|_2 \left\| \frac{\beta_0 \Delta t}{\alpha_0} ( \flux (\solArg{n}) - \flux (\solArg{ref}+\scaledDecoder(\redsolapproxArg{n}))) \right\|_2 \right| \\
        &= \|\obliqueprojector \|_2 \left| \gamma_1 \| \solArg{n}-\solArg{ref}-\scaledDecoder(\redsolapproxArg{n}) \|_2 - \gamma_2 \left\| \frac{\beta_0 \Delta t}{\alpha_0} ( \flux (\solArg{n}) - \flux (\solArg{ref}+\scaledDecoder(\redsolapproxArg{n}))) \right\|_2 \right|.
    \end{split}
\end{equation}
If $f$ is Lipschitz continuous with $L$ and $\Delta t$ is sufficiently small such that $\Delta t < \frac{\gamma_1 |\alpha_0|}{\gamma_2|\beta_0|L}$, we have
\begin{equation}\label{eq:I}
    \| \RN{1} \|_2 \geq \|\obliqueprojector\|_2\left(\gamma_1-\gamma_2\frac{|\beta_0|\Delta t L}{|\alpha_0|}\right)\|\solArg{n}-\solArg{ref}-\scaledDecoder(\redsolapproxArg{n})\|_2.
\end{equation}
With triangle inequality and Lipschitz continuity of $f$, we have
\begin{equation}\label{eq:II}
    \|\RN{2}\|_2 \leq \frac{1}{|\alpha_0|}\|\obliqueprojector\resRedApproxArg{n}(\solArg{ref}+\scaledDecoder(\redsolapproxArg{n}))\|_2 + \|\obliqueprojector\|_2\frac{1}{|\alpha_0|}\sum_{j=1}^k \left((|\alpha_j| + |\beta_j|\Delta t L)\|\solArg{n-j}-\solArg{ref}-\scaledDecoder(\redsolapproxArg{n-j})\|_2 \right).
\end{equation}
Combining Eq.~\eqref{eq:I} and ~\eqref{eq:II} yields
\begin{equation}\label{eq:CombIandII}
\begin{split}
    \| \solArg{n}-\solArg{ref}-\scaledDecoder(\redsolapproxArg{n}) \|_2 
    &\leq \frac{1}{\| \obliqueprojector \|_2 (\gamma_1-\gamma_2\frac{|\beta_0|\Delta t L}{|\alpha_0|})|\alpha_0|}\left \|\obliqueprojector \resRedApproxArg{n}(\solArg{ref}+\scaledDecoder(\redsolapproxArg{n}))
    \right \|_2 \\
    &+\sum_{j=1}^k \frac{|\alpha_j|+|\beta_j|\Delta t L}{(\gamma_1-\gamma_2\frac{|\beta_0|\Delta t L}{|\alpha_0|})|\alpha_0|}
    \left \|
    \solArg{n-j}-\solArg{ref}-\scaledDecoder(\redsolapproxArg{n-j})
    \right \|_2.
\end{split}
\end{equation}
The error bound for NM-Galerkin-HR Eq.~\eqref{eq:nmgalerkinerrorbound} is proved. Furthermore, noting that NM-LSPG-HR solution $\redsolapproxArg{n}$ minimizes the term $\| \obliqueprojector \resRedApproxArg{n}(\solArg{ref}+\scaledDecoder(\redsolapproxArg{n})) \|_2$ in Eq.~\eqref{eq:CombIandII} proves the error bound for NM-LSPG-HR Eq.\eqref{eq:nmlspgerrorbound}.
\hspace*{1em}\endproof}

From the error bound for NM-LSPG-HR, we know that the NM-LSPG-HR solutions
satisfy sequential minimization of the error bound. 

\section{Numerical results}\label{sec:numericresults}
We demonstrate the accuracy and speed-up of the nonlinear manifold reduced
order model for two advection-dominated problems: (i) a parameterized
1D inviscid Burgers equation in Section~\ref{sec:1dburgers} and (ii) a
parameterized 2D viscous Burgers equation with a large Reynolds number (i.e.,
the advection-dominated case) in Section~\ref{sec:2dburgers}. The ROMs are
trained with solution snapshot associated with train parameters in a parameter
space and are used to predict the solution of the parameter that is not
included in the train parameters. We refer this to the predictive case. The
accuracy of ROM solution $\solapproxFuncOnlyParam$ is assessed from its 
% mean squared state-space error:
% \begin{equation}\label{eq:relErr} 
%   \mseError  = \left.
%   \sqrt{\sum_{n=1}^\ntimedof\|\solapproxFuncArg{n}-\solFuncArg{n}\|_2^2}
%   \middle/ 
%   \sqrt{\sum_{n=1}^\ntimedof\|\solFuncArg{n}\|_2^2} \right.,
% \end{equation}
% where $\solFunc$ is the corresponding FOM solution 
% and the 
maximum relative error:
\begin{equation} \label{eq:maxErr}
  \maxError  = \left.
  \max_{n \in \nat{\ntimedof}}\left(\frac{\|\solapproxFuncArg{n}-\solFuncArg{n}\|_2}{\|\solFuncArg{n}\|_2}\right) \right.
\end{equation}
where $\solFunc$ is the corresponding FOM solution.  We also introduce the
projection errors for the lower bounds of LS-ROMs and NM-ROMs maximum relative
errors:
\begin{equation}\label{eq:linProjErr}
    \projError = \left.
    \sqrt{\sum_{n=1}^\ntimedof\|\left(\identity{}-\basismatspace\basismatspace^T\right)\left(\solFuncArg{n}-\solArg{ref}\left(\param\right)\right)\|_2^2}
    \middle/ 
    \sqrt{\sum_{n=1}^\ntimedof\|\solFuncArg{n}\|_2^2} \right.
\end{equation}
for linear subspace projection and
\begin{equation}\label{eq:nonlinProjErr}
    \projError = \left.
        \sqrt{\sum_{n=1}^\ntimedof\|\left(\solFuncArg{n}-\solArg{ref}\left(\param\right)\right)-\scaledDecoder \circ \scaledEncoder \left(\solFuncArg{n}-\solArg{ref}\left(\param\right)\right) \|_2^2}
    \middle/ 
    \sqrt{\sum_{n=1}^\ntimedof\|\solFuncArg{n}\|_2^2} \right.
\end{equation}
for nonlinear manifold projection, where $\basismatspace$ denotes a POD basis
matrix, and the scaled decoder $\scaledDecoder$ and the scaled encoder
$\scaledEncoder$ are a nonlinear manifold and its approximate inverse
function that are obtained from an autoencoder, respectively. The
computational cost is measured in terms of the CPU wall time.  Specifically,
timing is obtained by performing calculations on an Intel(R) Xeon(R) CPU
E5-2637 v3 @ 3.50 GHz and DDR4 Memory @ 1866 MT/s. The autoencoders are trained
on a NVIDIA Quadro M6000 GPU with 3072 NVIDIA CUDA Cores and 12 GB GDDR5 GPU
Memory using PyTorch \cite{paszke2019pytorch} which is the open source machine learning frame work.
% The autoencoders are trained on the GPU machine at Lawrence Livermore
% National Laboratory using PyTorch which is the open source machine learning
% frame work. The GPU machine has 4 NVIDIA Quadro M6000 GPUs, each with 3072
% NVIDIA CUDA Cores and 12 GB GDDR5 GPU Memory.

\subsection{1D inviscid Burgers' equation}\label{sec:1dburgers}
We consider a parameterized 1D inviscid Burgers' equation
\begin{align}\label{eq:1dburgers_eq} 
   \frac{\partial u(x,t;\mu)}{\partial t} &+ u(x,t;\mu)\frac{\partial u(x,t;\mu)}{\partial x} = 0,\\ 
   x &\in \spaceDomain = [0,2]\\
   t &\in[0,\totaltime],
\end{align} 
where $u\in\RR{}$ denotes a scalar-valued time dependent state variable 
with the periodic boundary condition
\begin{equation}
    u(2,t;\mu)=u(0,t;\mu)
\end{equation}
and the initial condition
\begin{equation}
    u(x,0;\mu)= 
\left\{
\begin{array}{ll}
1+\frac{\mu}{2}\left(\sin{\left(2\pi x-\frac{\pi}{2}\right)}+1\right) \quad &\text{if } 0 \leq x \leq 1 \\
1 \quad &\text{otherwise}
\end{array}
\right. 
\end{equation}
where $\mu \in \paramDomain=[0.9,1.1]$ is a parameter. Discretizing the space
domain $\spaceDomain$ into $\nx-1$ uniform mesh gives $\nx$ grid points
$x_i=(i-1)\Delta x$ where $i \in \{1,2,\cdots, \nx \}$ and $\Delta x =
\frac{2}{\nx-1}$. We denote the discrete solutions on grid points as
$u_i(t;\mu)=u(x_i,t;\mu)$, where $i \in \nat{\nx}$. Then, the backward
difference scheme $\frac{\partial u}{\partial x} \approx
\frac{u_i-u_{i-1}}{\Delta x}$ yields the semi-discretized equation which is
written by
\begin{equation}\label{eq:1dBurgersSemi}
    \frac{d\discreteU}{dt}=\flux(\discreteU)
\end{equation}
where $\discreteU=(u_1, u_2, \cdots, u_{\nx-1})^T$ and $\flux: \RR{\nx-1} \rightarrow \RR{\nx-1}$ is in the form
\begin{equation}
    \flux (\discreteU) = -\frac{1}{\Delta x} \left(\boldsymbol{M}\discreteU \odot \discreteU +  \boldsymbol{B} \discreteU \right)
\end{equation}
where 
\begin{equation}
\boldsymbol{M} = \bmat{ 1    &      &      &
\\                     -1    & 1    &      & 
\\                           &\ddots&\ddots&
\\                           &      & -1   &  1}
_{(\nx-1)\times(\nx-1)}, 
\quad
\boldsymbol{B} = \bmat{u_{\nx-1}&         &         &
\\                              & 0       &         &
\\                              &         &\ddots   &
\\                              &         &         & 0}
_{(\nx-1)\times(\nx-1)},
\end{equation} 
with $\odot$ denoting element-wise multiplication.

For a time integrator, we use the backward Euler scheme with time step size
$\Delta t = \frac{T}{\nt}$, where $T$ is final time and $\nt$ is the number of
time steps. We set $T=0.5$, $\nx=1001$, and $\nt=500$.

For the training process, we collect solution snapshots associated with the
parameter $\mu \in \paramDomain_{train}=\{0.9,1.1\}$ such that $\ntrain=2$ at
which the FOM is solved. Then, the number of train data points is
$\ntrain\cdot(\nt+1)=1002$ and $10\%$ of the train data are used for validation
purpose. We employ the Adam optimizer \cite{kingma2014adam} for SGD with
initial learning rate $0.001$ which decreases by a factor of $10$ when a
training loss stagnates for $10$ successive training epochs. We set the number
of nodes in the hidden layer of the encoder,  $M_1=2000$, and the number of
nodes in the hidden layer of the decoder, $M_2=12024$.
The weights and bias of the autoencoder are initialized via Kaiming
initialization \cite{he2015delving}. The size of the batch is $20$ and the
maximum number of epochs is $10,000$. The training process is stopped if the
loss on the validation dataset stagnates for $200$ epochs. 

\begin{figure}[!htbp]
    \centering
    \includegraphics[width=0.4\textwidth]{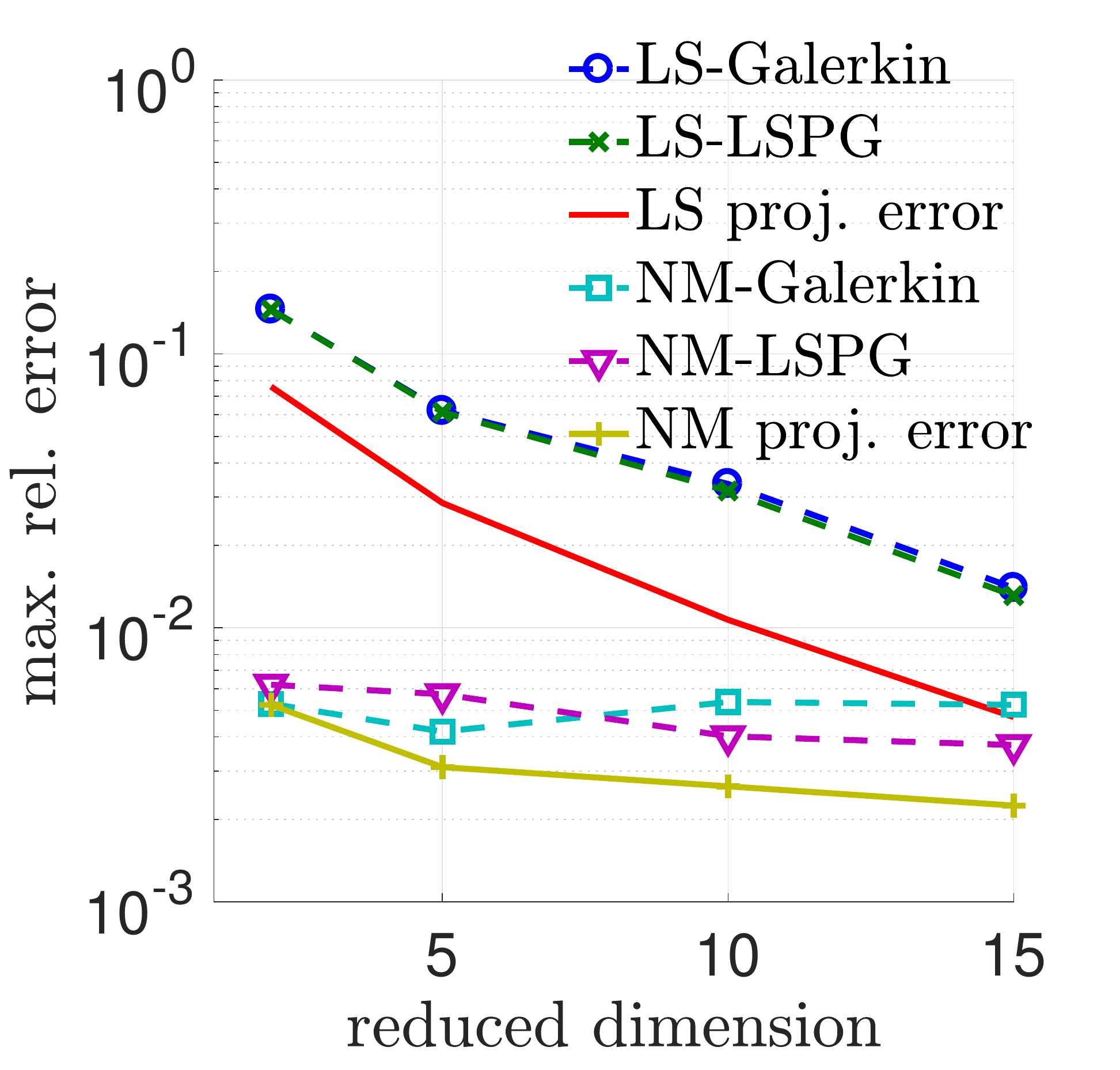}
    \caption{1D Burgers' equation. Relative errors vs reduced dimensions.}
    \label{fg:1dErrVSredDim}
\end{figure}

After the training is done, the NM-ROMs and LS-ROMs solve the
Eq.~\eqref{eq:1dburgers_eq} with the target parameter $\mu=1$ which is not
included in the train dataset for training the autoencoder and the linear
subspace.  Fig.~\ref{fg:1dErrVSredDim} shows the relative error versus the
reduced dimension $\nbasisspace$. 
It also shows the projection errors for
LS-ROMs and NM-ROMs, which are defined in \eqref{eq:linProjErr} and
\eqref{eq:nonlinProjErr}. These are the lower bounds for LS-ROMs and NM-ROMs,
respectively. As expected the relative errors for the NM-ROMs are lower than
the ones for the LS-ROMs.  We even observe that the relative errors of NM-ROMs
are even lower than the lower bounds of LS-ROMs.

\begin{figure}[!htbp]
    \centering
    \includegraphics[width=0.4\textwidth]{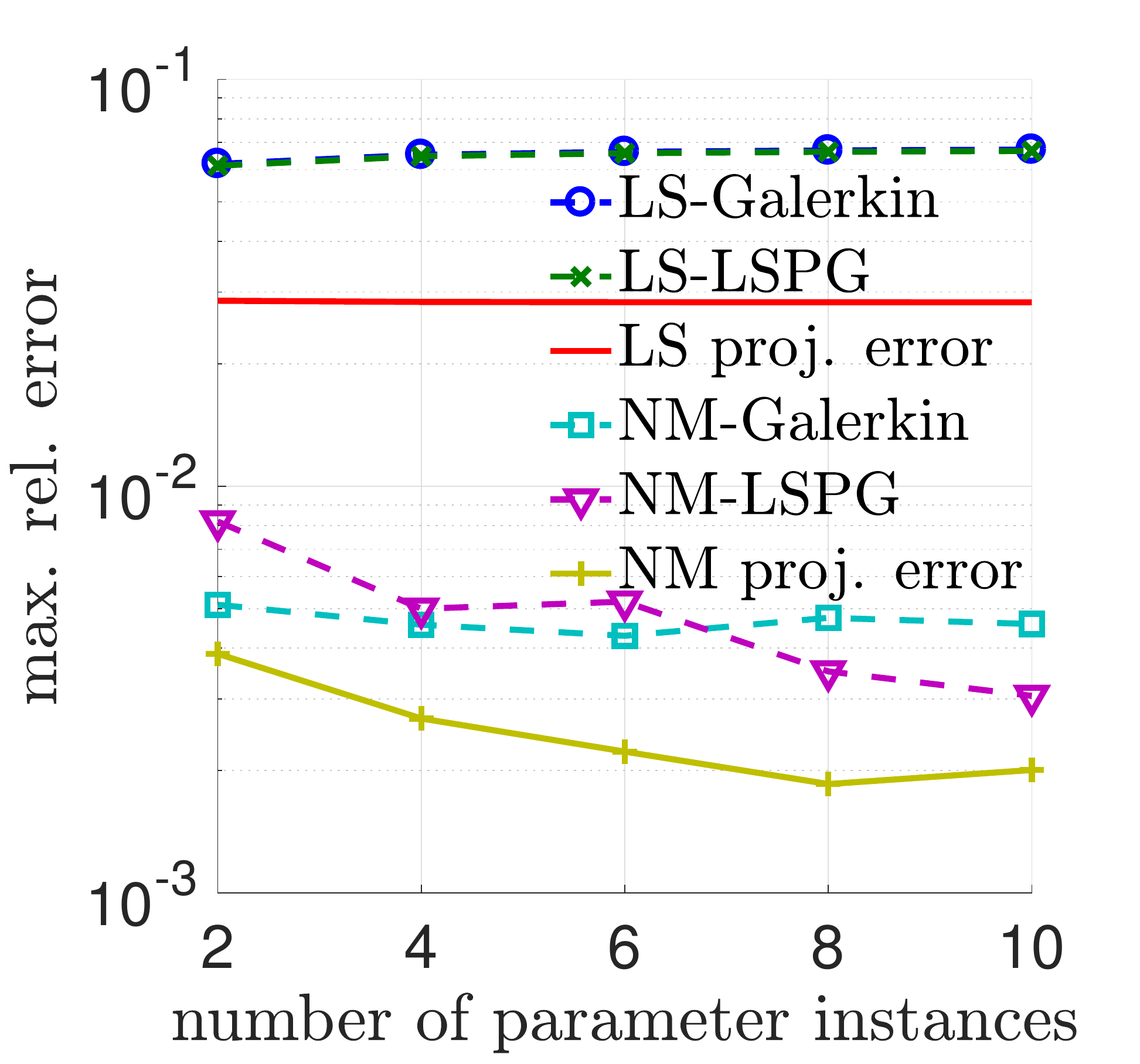}
    \caption{1D Burgers' equation. Relative errors vs the number of parameter
    instances. We use $\paramDomain_{train}=\{0.9,1.1\}$ for $\ntrain=2$,
    $\paramDomain_{train}=\{0.9,0.95,1.05,1,1\}$ for $\ntrain=4$,
    $\paramDomain_{train}=\{0.9,0.9+\frac{1}{3},0.9+\frac{2}{3},1+\frac{1}{3},1+\frac{2}{3},1.1\}$
    for $\ntrain=6$, and
    $\paramDomain_{train}=\{0.9,0.925,0.95,0.975,1.025,1.05,1.075,1.1\}$ for
    $\ntrain=8$.} 
    \label{fg:1dErrVSParam}
\end{figure}

To see the trends regarding the number of training parameter instances, we
increase the number of parameters starting from $\ntrain=2$ with the fixed reduced
dimension $\nbasisspace=5$ to achieve less than $1\%$ maximum relative error for
NM-ROMs. In Fig.~\ref{fg:1dErrVSParam}, we observe that $\ntrain=2$ is
enough.

LS-ROMs with $\nbasisspace=5$ are able to achieve speed-up, but their
accuracies are not as good as NM-ROMs. For example, LS-ROMs are about $5$ to
$6$ times faster than FOM on average and have a maximum relative error of $6$ \%.
NM-ROMs solve the problem with less than the maximum relative error of
$1$ \%. For LS-ROMs, a hyper-reduction improves speed-up (e.g., $9$ to $10$
times faster than FOM) but accuracy doesn't get better. On the other hand,
NM-ROMs without a hyper-reduction do not achieve any speed-up with respect to
the corresponding FOM simulation. For example, the FOM simulation takes $1.30$
seconds, while the NM-Galerkin and NM-LSPG with $\nbasisspace=5$ takes $1.67$
and $1.35$ seconds, respectively. Therefore, the hyper-reduction is essential
to achieve a speed-up with a reasonable accuracy for the NM-ROMs. Now, we
compute the maximum relative error and wall-clock time for the hyper-reduced
ROMs, i.e., NM-LSPG-HR and LS-LSPG-HR, by varying the number of residual basis
and residual samples with the fixed number of training parameter instances
$\ntrain=2$ and the reduced dimension $\nbasisspace=5$ and show the results in
Table~\ref{tb:1dDEIMtest}.  Although the LS-LSPG-HR can achieve a better
speed-up than the NM-LSPG-HR, the relative error of the LS-LSPG-HR is too large, e.g., the relative errors of around $6 \%$. On the other hand, the
NM-LSPG-HR achieves much better accuracy, i.e., a relative error of around $1
\%$, than the LS-LSPG-HR with a speedup of higher than $2$.
\begin{table}[!htbp]
\caption{The top 6 maximum relative errors and wall-clock times at different
  numbers of residual basis and samples which range from $30$ to
  $50$.}\label{tb:1dDEIMtest}
\centering
\resizebox{\textwidth}{!}{\begin{tabular}{|c|c|c|c|c|c|c|c|c|c|c|c|c|}
\hline
 & \multicolumn{6}{|c|}{NM-LSPG-HR} & \multicolumn{6}{|c|}{LS-LSPG-HR}\\
\hline
Residual basis & 31 & 33 & 36 & 32 & 40 & 32 & 30 & 30 & 30 & 31 & 41 & 41\\
\hline
Residual samples & 47 & 49 & 40 & 47 & 42 & 46 & 47 & 48 & 49 & 49 & 49 & 48\\
\hline
Max. rel. error (\%) & 1.03 & 1.07 & 1.18 & 1.23 & 1.23 & 1.25 & 6.07 & 6.08 & 6.08 & 6.09 & 6.11 & 6.11\\
\hline
Wall-clock time (sec) & 0.63 & 0.51 & 0.49 & 0.50 & 0.51 & 0.50 & 0.14 & 0.13 & 0.13 & 0.23 & 0.14 & 0.13\\
\hline
Speed-up & 2.07 & 2.53 & 2.63 & 2.62 & 2.56 & 2.62 & 9.29 & 9.80 & 9.71 & 5.65 & 9.63 & 9.82 \\
\hline
\end{tabular}}
\end{table}

\begin{figure}[!htbp]
  \centering
  \subfigure[FOM]{
  \includegraphics[width=0.3\textwidth]{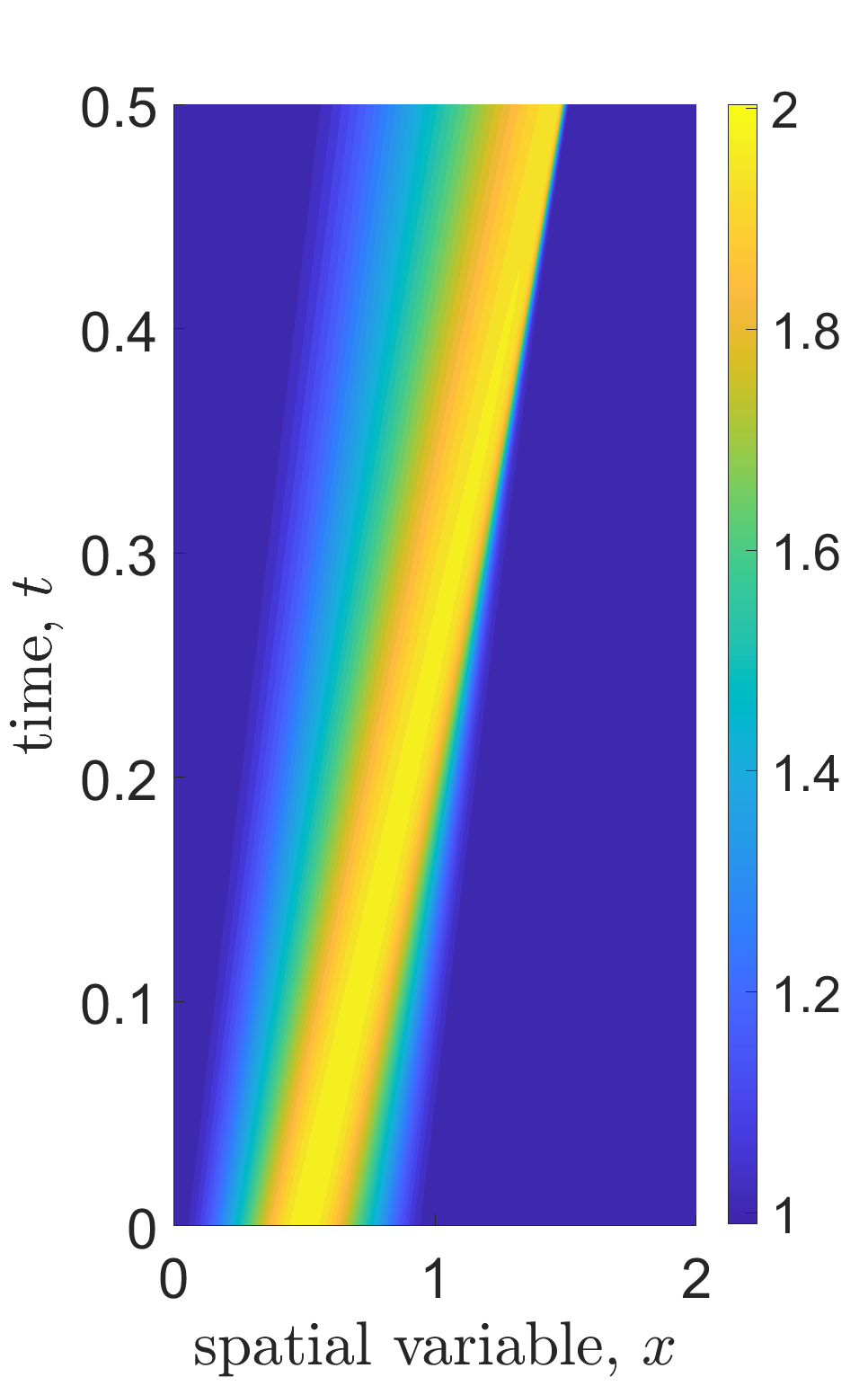}
}~~~~~~~
\subfigure[NM-LSPG-HR]{
  \includegraphics[width=0.3\textwidth]{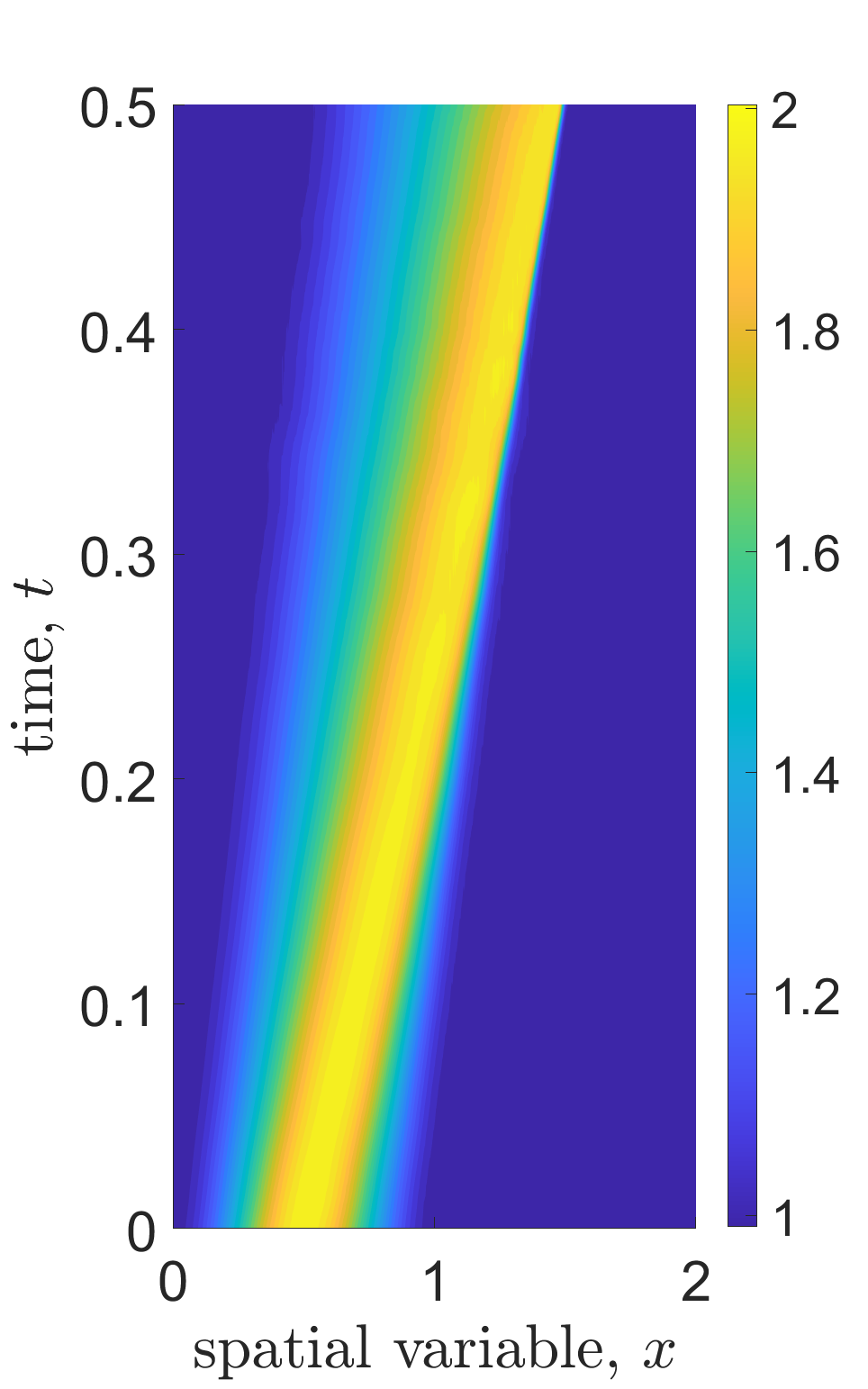}
}~~~~~~~
\subfigure[LS-LSPG-HR]{
  \includegraphics[width=0.3\textwidth]{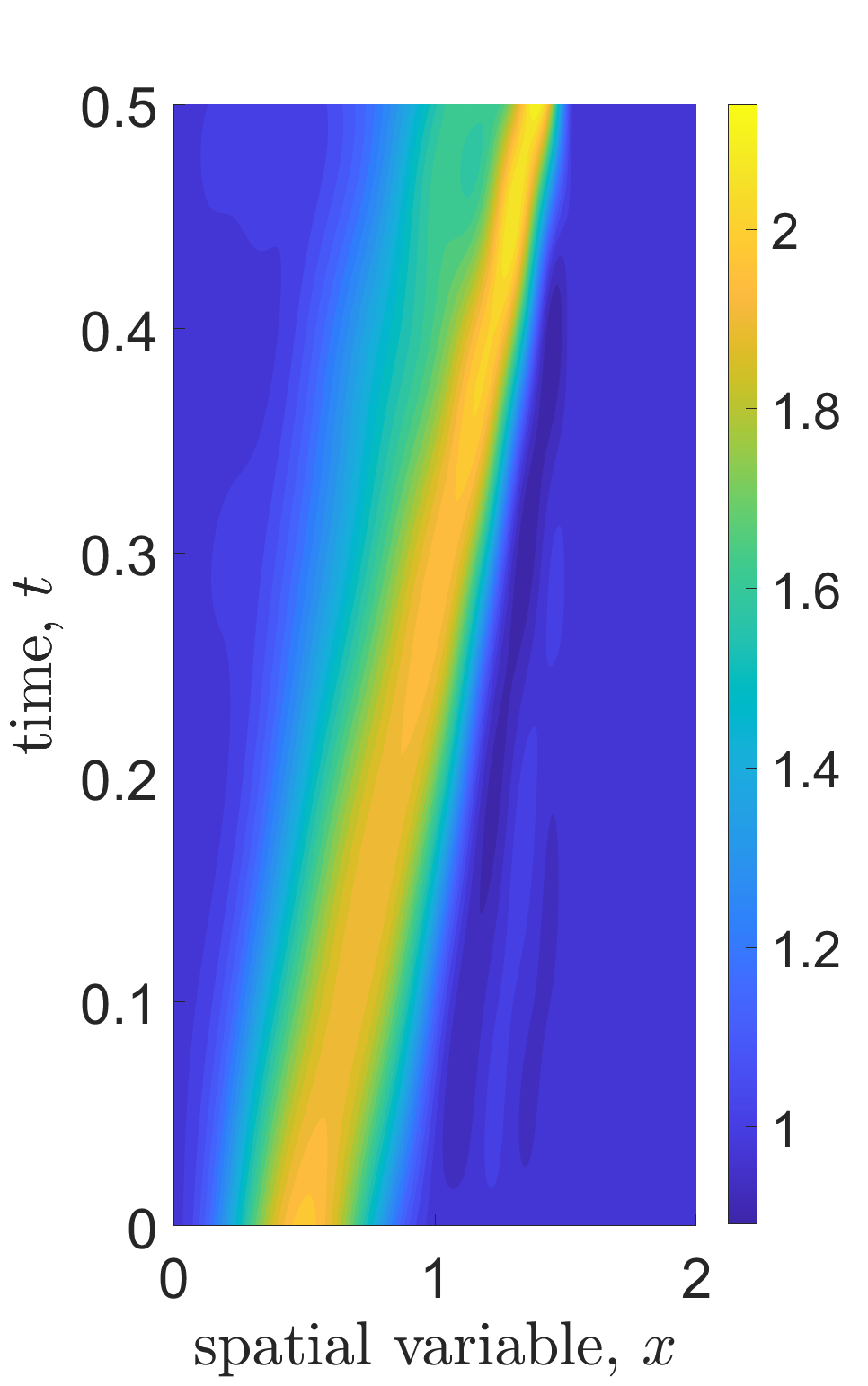}
}~~~~~~~
  \caption{A space--time solution instances of FOM and ROMs for 1D Burgers'
  equation.}
  \label{fg:1dburgersFOMROMS}
\end{figure}

Fig.~\ref{fg:1dburgersFOMROMS} shows solutions in both space and time domain of
FOM, NM-LSPG-HR, and LS-LSPG-HR with the reduced dimension being
$\nbasisspace=5$. For NM-LSPG-HR, $31$ residual basis and $47$ residual samples are used and for LS-LSPG-HR, $30$ residual basis and $47$ residual samples are used. In
fact, the NM-LSPG-HR is able to achieve an accuracy as good as the NM-LSPG for
some combinations of the small number of residual basis and residual samples.

\begin{figure}[!htbp]
  \centering
  \includegraphics[width=0.5\textwidth]{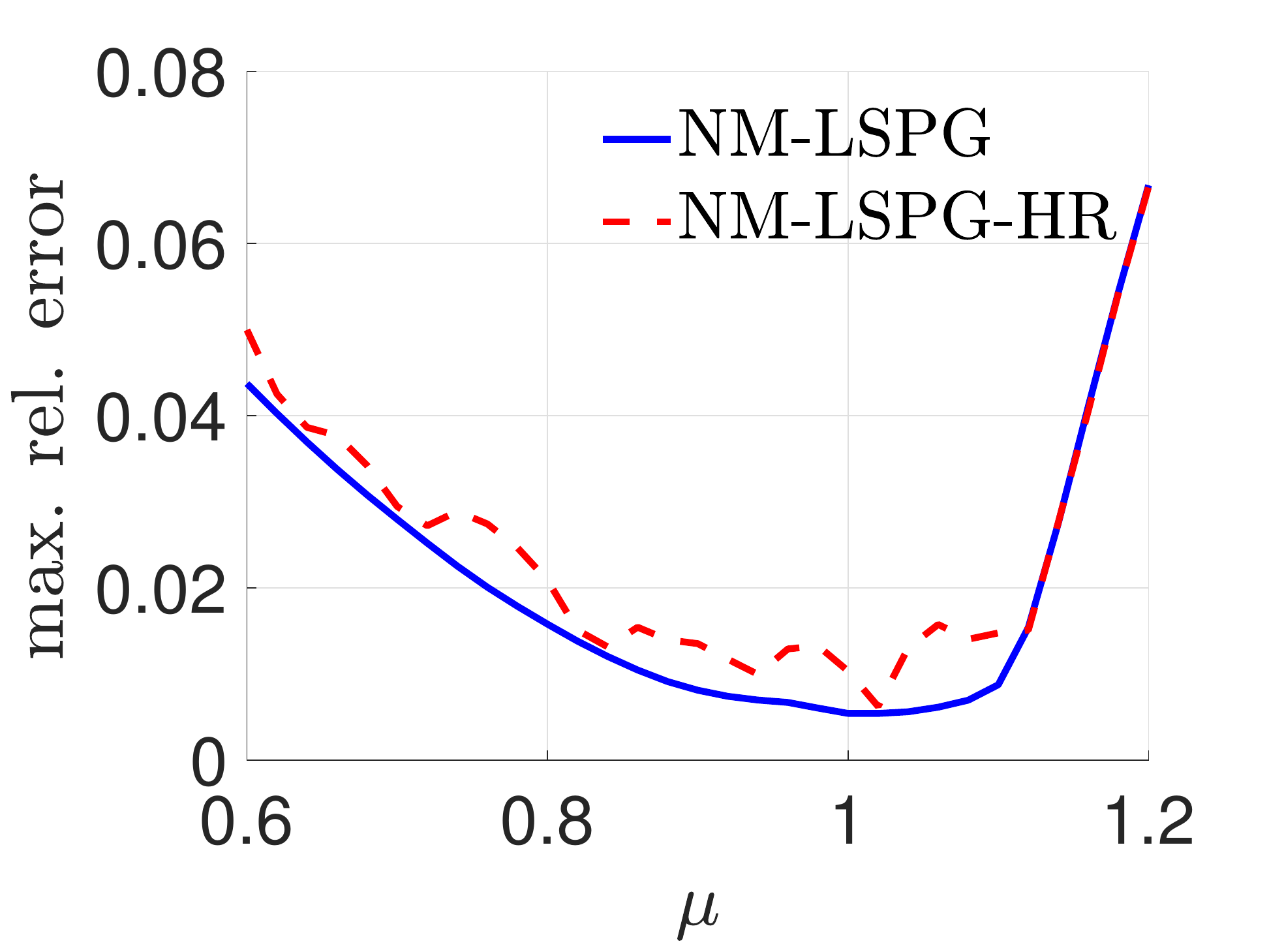}
  \caption{The comparison of the NM-LSPG-HR and NM-LSPG on the maximum relative
  error vs $\mu$}
  \label{fg:1dPredictionTest}
\end{figure}

We look into the numerical tests to see the generalization capability of the
NM-LSPG and NM-LSPG-HR, i.e., the robustness of the NM-LSPG and NM-LSPG-HR
outside of the trained domain. The training sample point set,
$\mu\in\paramDomain_{train}=\{0.9,1.1\}$, is used to train a NM-LSPG-HR. Then
the trained NM-LSPG-HR model is used to predict the following parameter points,
$\mu\in\paramDomain_{test} = \{\mu | \mu = 0.6+0.02i, i = 0,1,\cdots,30\}$.
The residual basis dimension and the number of residual
samples for each test case are given in Table~\ref{tb:1dPredictionTest}.
Fig.~\ref{fg:1dPredictionTest} shows the maximum relative error over the test
range of the parameter points.  Note that the NM-LSPG and NM-LSPG-HR are the
most accurate within the range of the training points, i.e., $[0.9, 1.1]$.  As
the parameter points go beyond the training parameter domain, the accuracy of
the NM-LSPG and NM-LSPG-HR start to deteriorate gradually. This implies that
the NM-LSPG and NM-LSPG-HR have a trust region. Its trust region should be
determined by an application. For example, if the application is okay with the
maximum relative error of $10$ \%, then the trust region of this particular
NM-LSPG-HR is $[0.6, 1.2]$. However, if the application requires a higher
accuracy, e.g., the maximum relative error of $2$ \%, then the trust region of
the NM-LSPG-HR is $[0.82, 1.12]$.
Note that the average speed-up of the NM-LSPG-HR for all the test cases is
$2.72$ (see Table~\ref{tb:1dPredictionTest}).

\begin{table}[!htbp]
\caption{The residual basis dimension and the number of residual samples for
  each test parameter $\mu$. The wall-clock time and the speed-up of the NM-LSPG-HR with respect to
  the corresponding FOM are also reported.}\label{tb:1dPredictionTest}
\centering
\resizebox{0.6\textwidth}{!}{\begin{tabular}{|c|c|c|c|c|}
\hline
$\mu$ & Residual basis & Residual samples & Wall-clock time (sec) & Speed-up \\ \hline
0.60	&	46	&	48	&	0.54	&	2.41	\\	\hline
0.62	&	37	&	39	&	0.46	&	2.83	\\	\hline
0.64	&	37	&	39	&	0.47	&	2.77	\\	\hline
0.66	&	44	&	46	&	0.52	&	2.50	\\	\hline
0.68	&	42	&	44	&	0.48	&	2.71	\\	\hline
0.70	&	42	&	44	&	0.48	&	2.71	\\	\hline
0.72	&	42	&	44	&	0.48	&	2.71	\\	\hline
0.74	&	42	&	44	&	0.48	&	2.71	\\	\hline
0.76	&	43	&	45	&	0.53	&	2.45	\\	\hline
0.78	&	36	&	45	&	0.46	&	2.83	\\	\hline
0.80	&	38	&	47	&	0.47	&	2.77	\\	\hline
0.82	&	38	&	47	&	0.47	&	2.77	\\	\hline
0.84	&	38	&	47	&	0.47	&	2.77	\\	\hline
0.86	&	38	&	47	&	0.47	&	2.77	\\	\hline
0.88	&	37	&	46	&	0.46	&	2.83	\\	\hline
0.90	&	33	&	33	&	0.45	&	2.89	\\	\hline
0.92	&	33	&	33	&	0.46	&	2.83	\\	\hline
0.94	&	33	&	33	&	0.46	&	2.83	\\	\hline
0.96	&	33	&	33	&	0.45	&	2.89	\\	\hline
0.98	&	31	&	47	&	0.45	&	2.89	\\	\hline
1.00	&	31	&	47	&	0.45	&	2.89	\\	\hline
1.02	&	33	&	49	&	0.48	&	2.71	\\	\hline
1.04	&	31	&	48	&	0.46	&	2.83	\\	\hline
1.06	&	30	&	48	&	0.46	&	2.83	\\	\hline
1.08	&	33	&	39	&	0.48	&	2.71	\\	\hline
1.10	&	33	&	40	&	0.48	&	2.71	\\	\hline
1.12	&	33	&	42	&	0.48	&	2.71	\\	\hline
1.14	&	44	&	49	&	0.54	&	2.41	\\	\hline
1.16	&	31	&	48	&	0.48	&	2.71	\\	\hline
1.18	&	31	&	48	&	0.47	&	2.77	\\	\hline
1.20	&	44	&	48	&	0.57	&	2.28	\\	\hline
\end{tabular}}
\end{table}

\subsection{2D Burgers' equation}\label{sec:2dburgers}
\begin{figure}[htbp]
  \centering
  \subfigure[$Re=100$]{
  \includegraphics[width=0.45\textwidth]{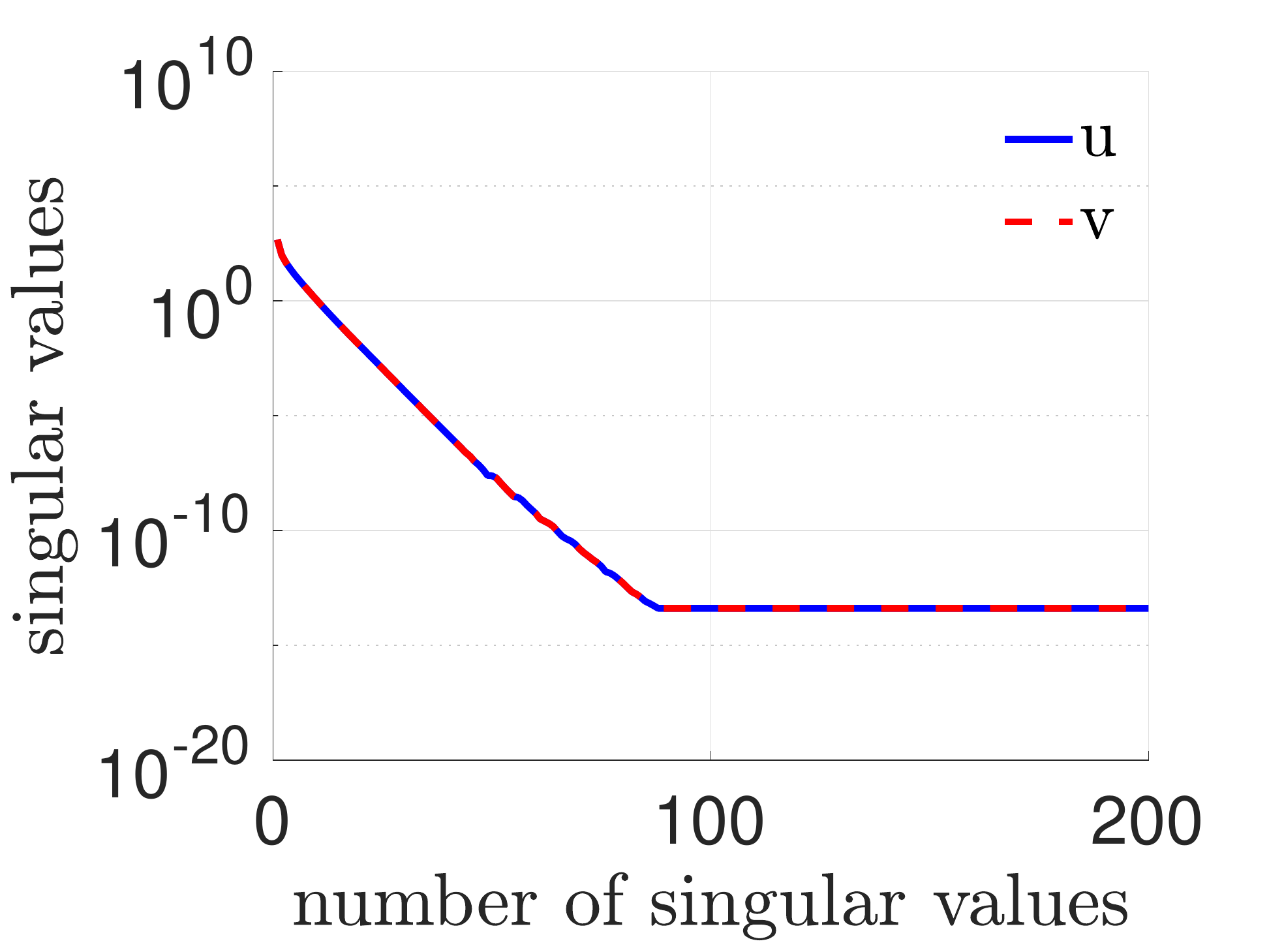}
  }
\subfigure[$Re=10000$]{
  \includegraphics[width=0.45\textwidth]{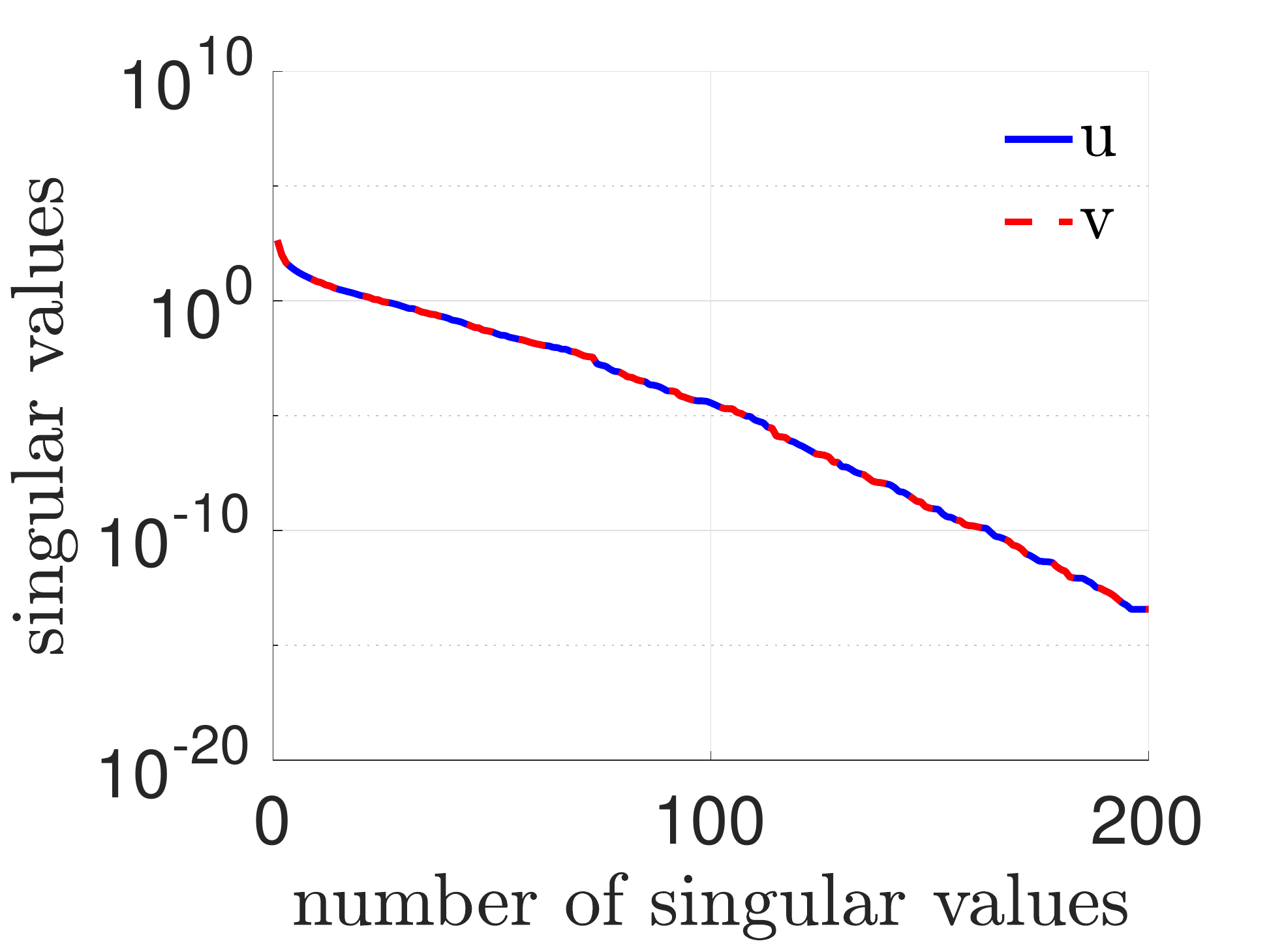}
  }
  \caption{The effect of Reynolds number on the singular value decay.}
  \label{fg:2dburgersSingular}
\end{figure}

\begin{figure}[!htbp]
  \centering
  \subfigure[FOM, $u$]{
  \includegraphics[width=0.3\textwidth]{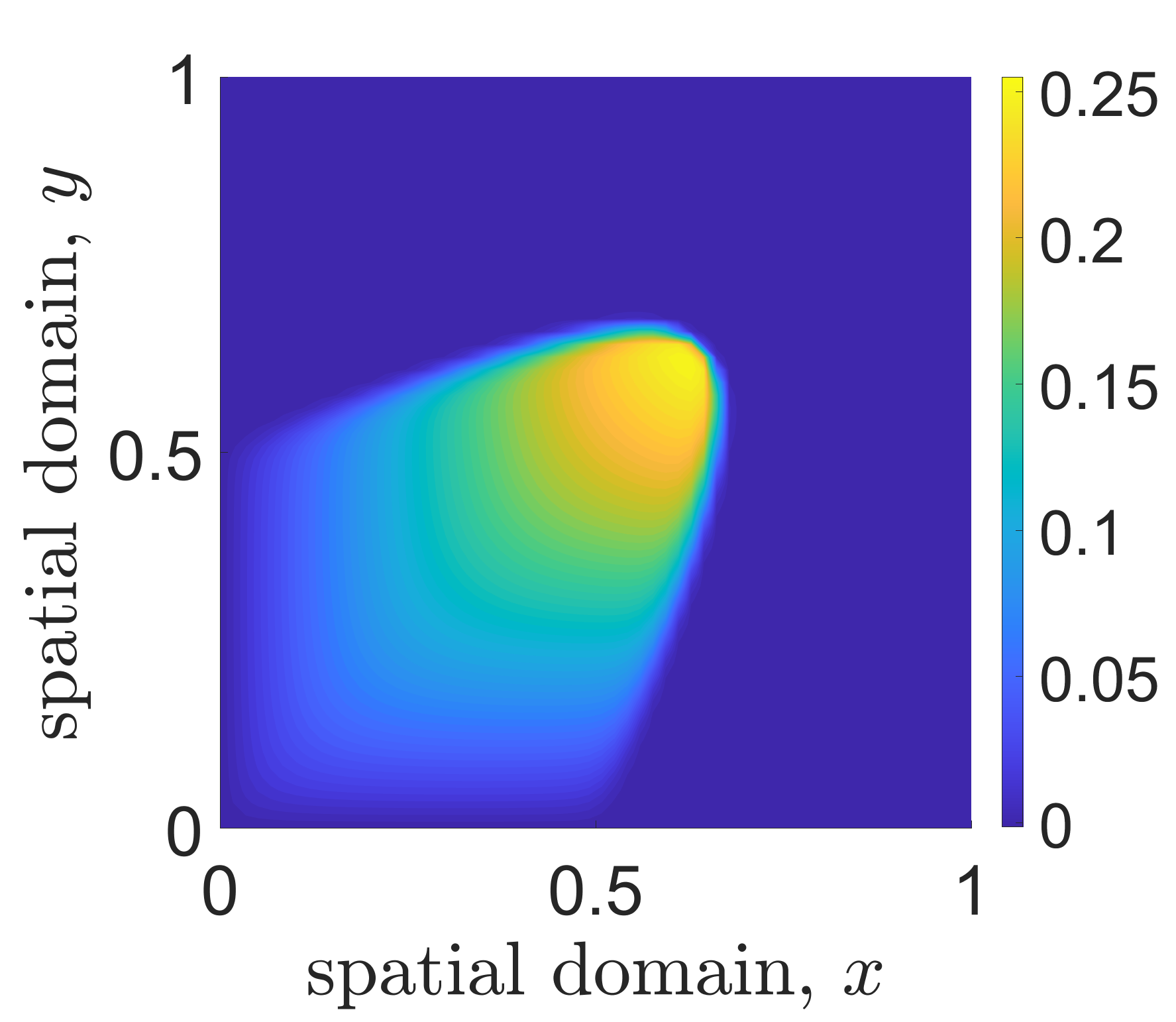}
}~~~~~~~
\subfigure[NM-LSPG-HR, $u$]{
  \includegraphics[width=0.3\textwidth]{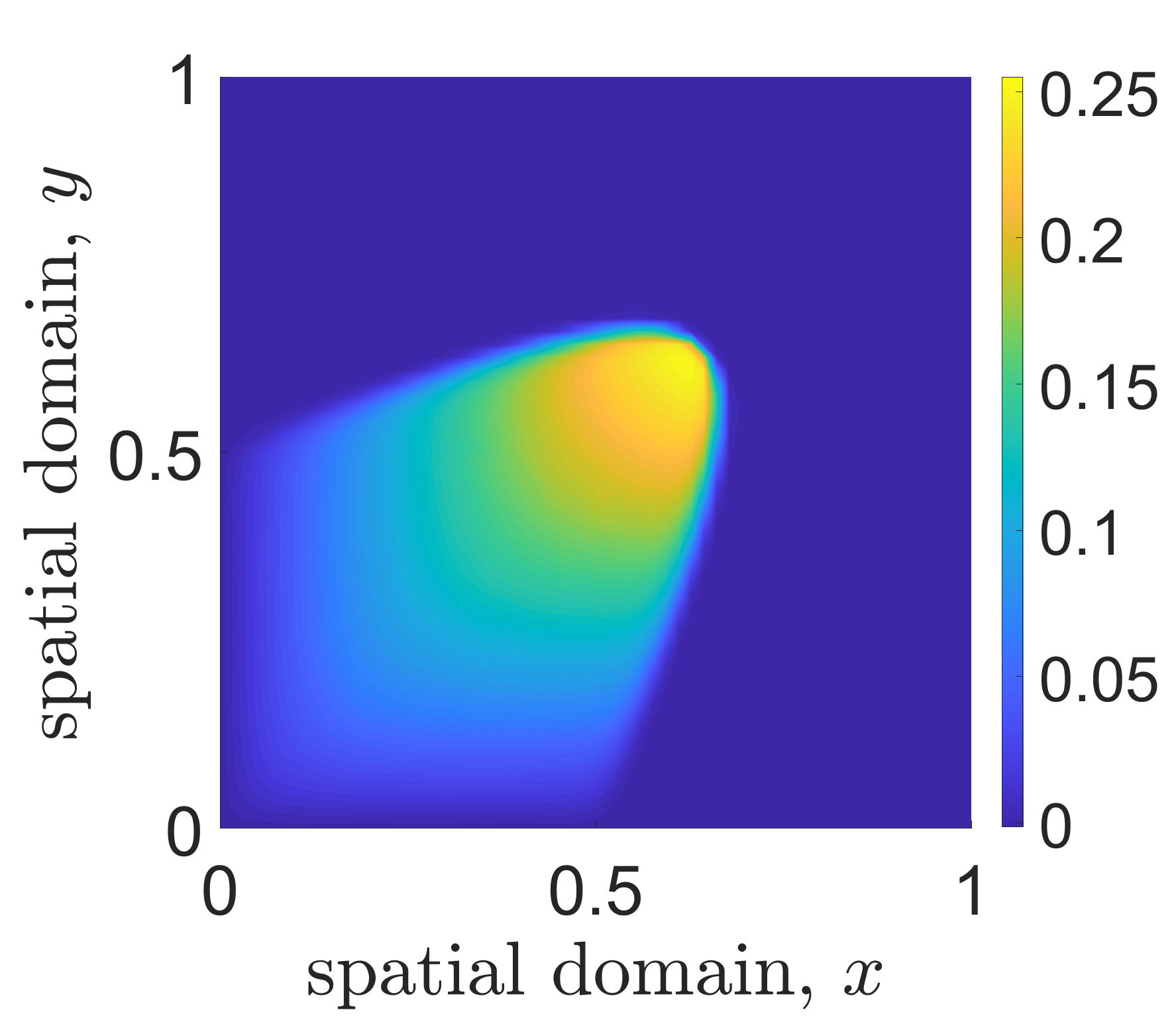}
}~~~~~~~
\subfigure[LS-LSPG-HR, $u$]{
  \includegraphics[width=0.3\textwidth]{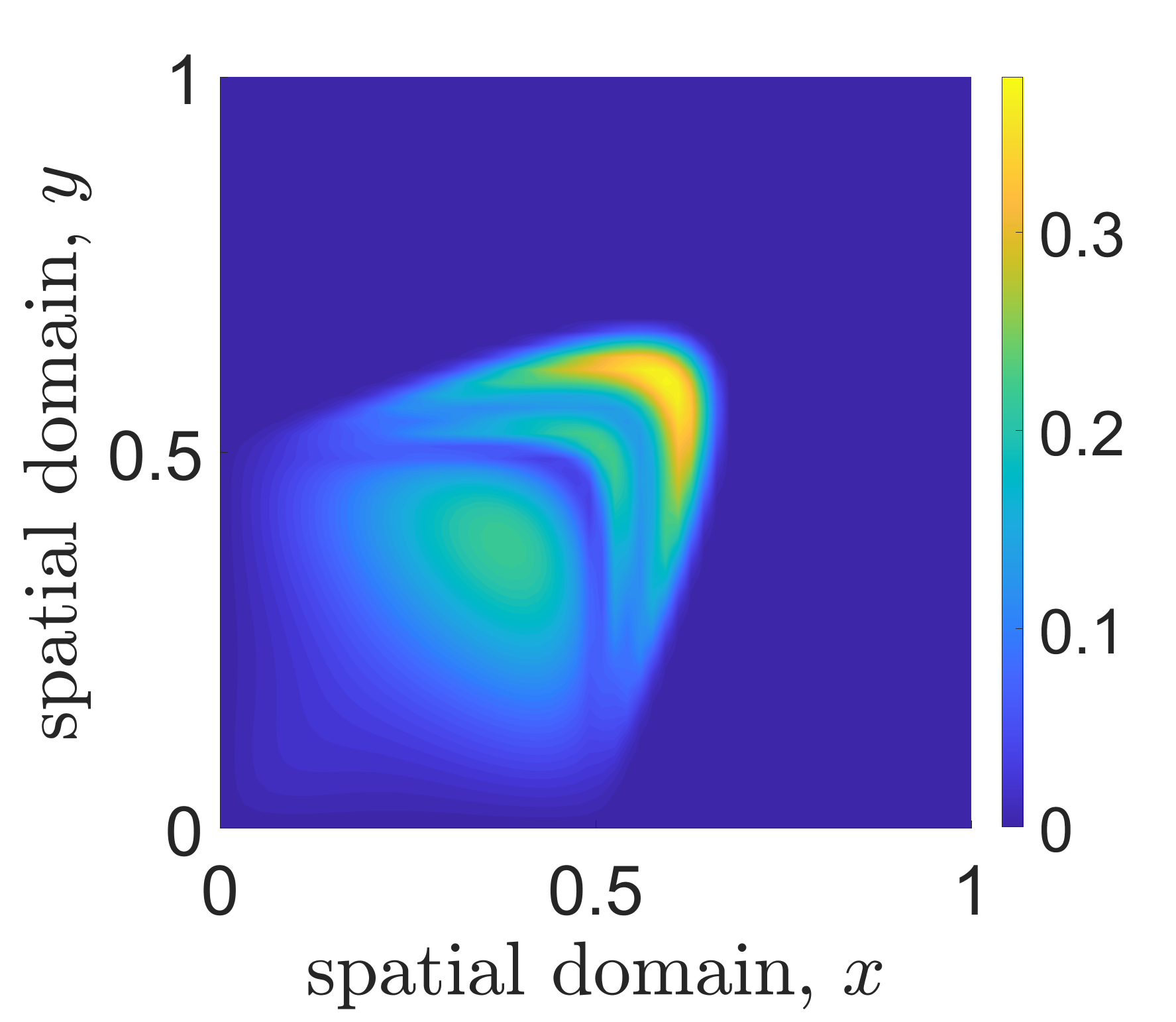}
}~~~~~~~\\
  \subfigure[FOM, $v$]{
  \includegraphics[width=0.3\textwidth]{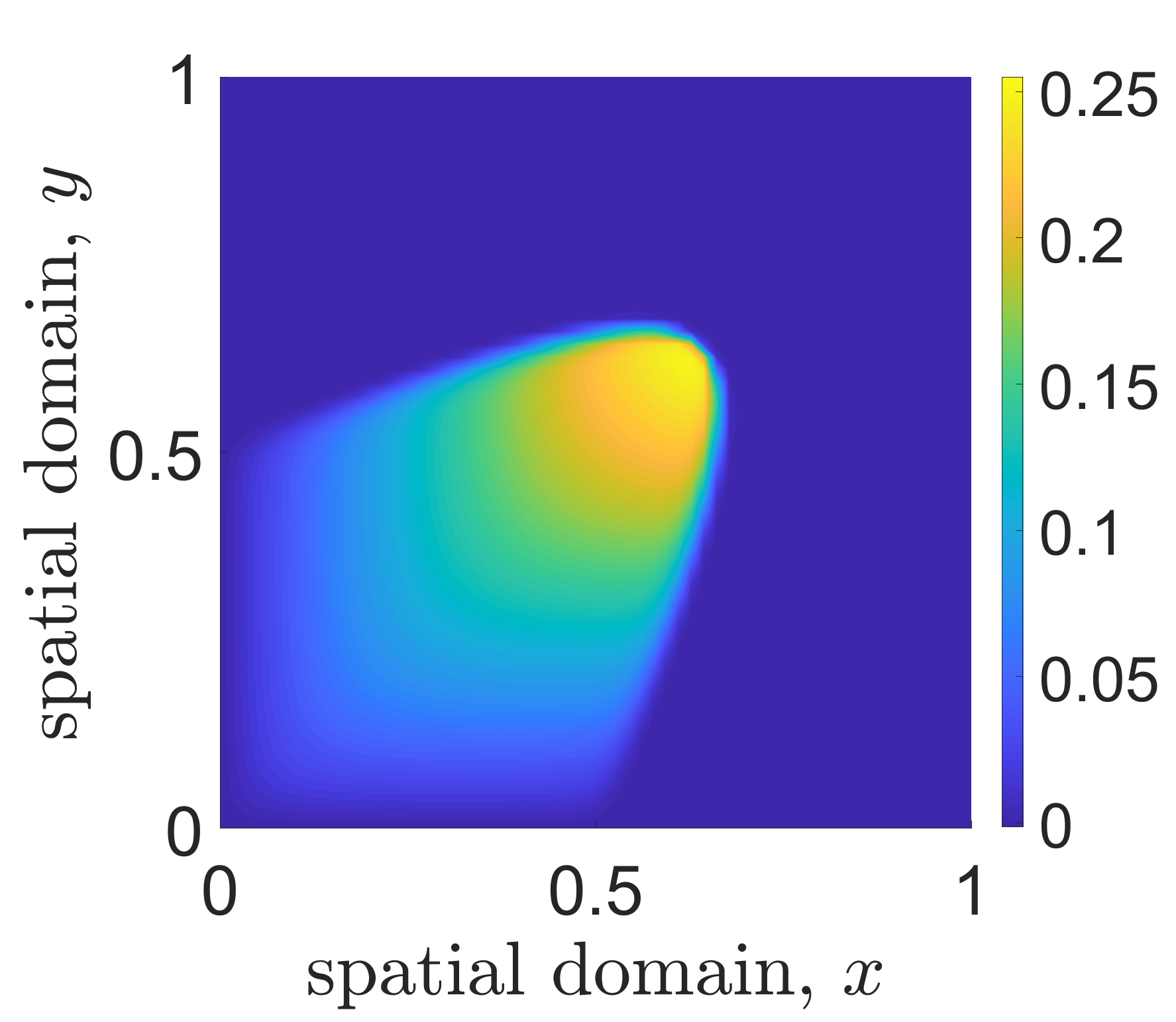}
}~~~~~~~
\subfigure[NM-LSPG-HR, $v$]{
  \includegraphics[width=0.3\textwidth]{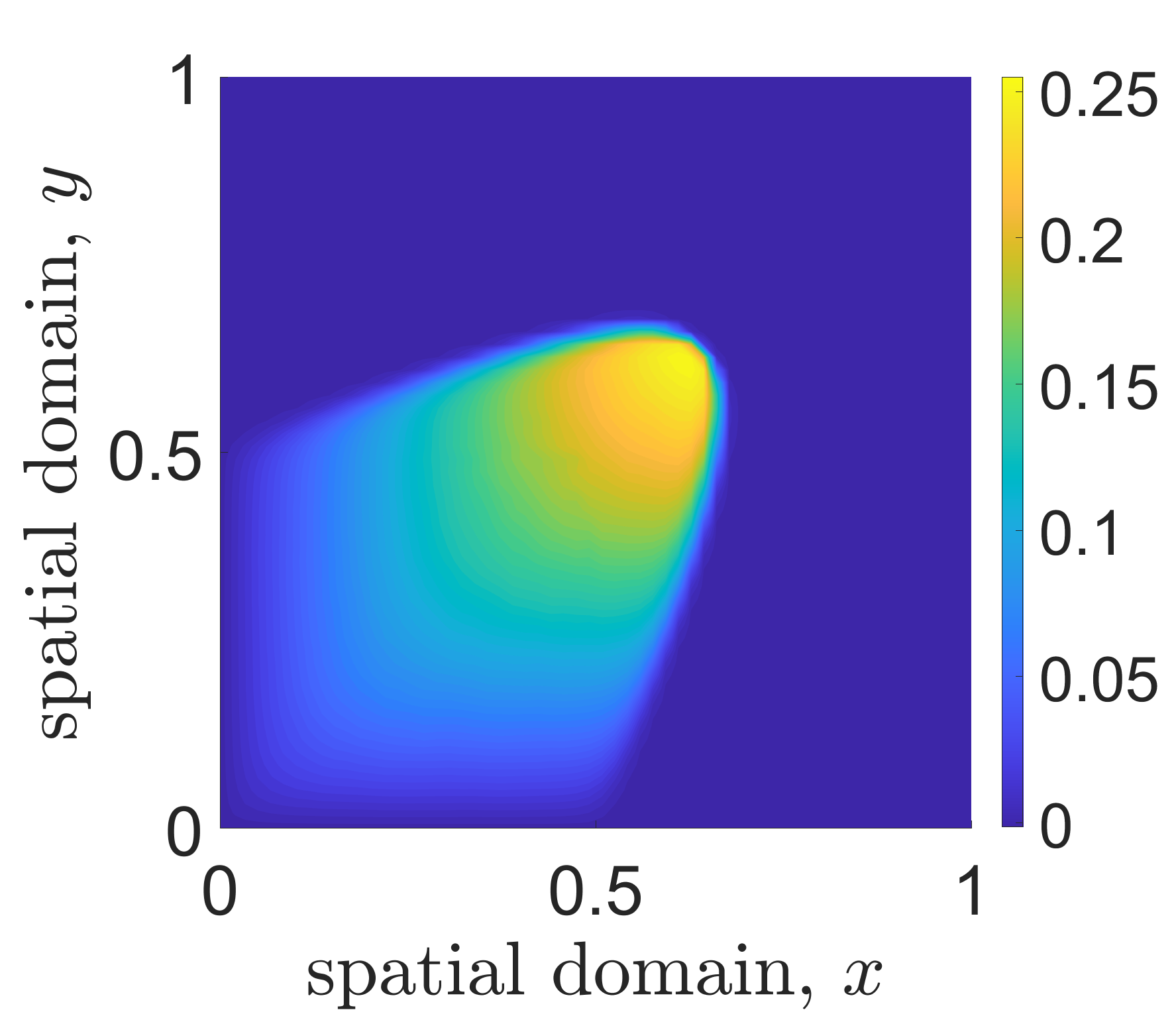}
}~~~~~~~
\subfigure[LS-LSPG-HR, $v$]{
  \includegraphics[width=0.3\textwidth]{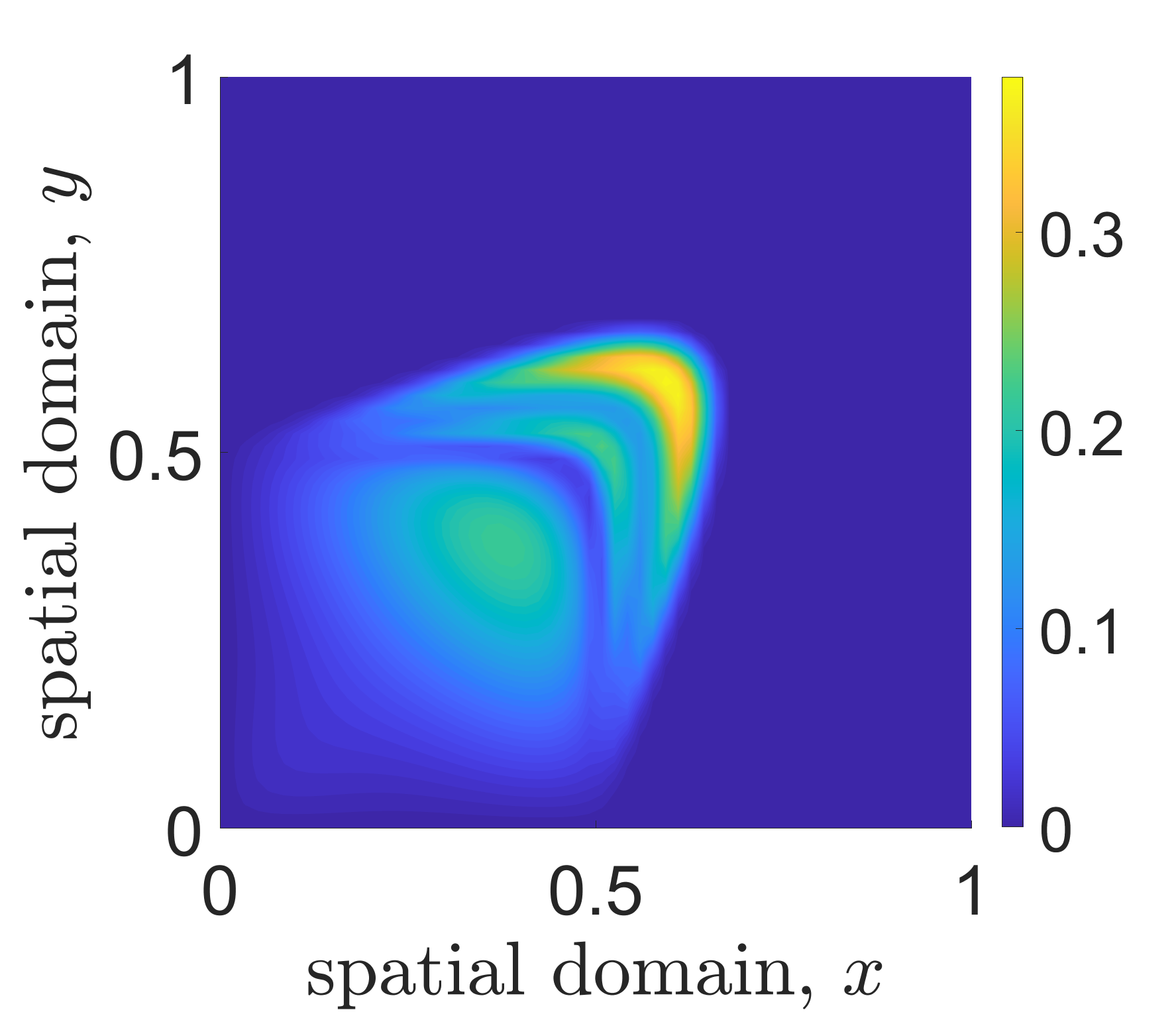}
}~~~~~~~
  \caption{Solution snapshots of FOM, NM-LSPG-HR, and LS-LSPG-HR at $t=2$}
  \label{fg:2dburgersFOMROMS}
\end{figure}

We now consider a parameterized 2D viscous Burgers' equation
\begin{align}\label{eq:2dburgers_eq} 
\frac{\partial u}{\partial t} + u\frac{\partial u}{\partial x} + v\frac{\partial u}{\partial y} &= \frac{1}{Re}\left(\frac{\partial^2 u}{\partial x^2}+\frac{\partial^2 u}{\partial y^2}\right) \\ 
\frac{\partial v}{\partial t} + u\frac{\partial v}{\partial x} + v\frac{\partial v}{\partial y} &= \frac{1}{Re}\left(\frac{\partial^2 v}{\partial x^2}+\frac{\partial^2 v}{\partial y^2}\right) \\ 
(x,y) &\in \spaceDomain = [0,1] \times [0,1] \\
t &\in[0,2],
\end{align}
with the boundary condition
\begin{equation}
    u(x,y,t;\mu)=v(x,y,t;\mu)=0 \quad \text{on} \quad \Gamma = \left\{(x,y)|x\in\{0,1\},y\in\{0,1\}\right\}
\end{equation}
and the initial condition
\begin{align}
    u(x,y,0;\mu) = \left\{
                                  \begin{array}{ll}
                                    \mu \sin{(2\pi x)}\cdot \sin{(2\pi y)} \quad &\text{if } (x,y) \in [0,0.5]\times [0,0.5] \\
                                    0 \quad &\text{otherwise}
                                  \end{array}
                                  \right.  \\
    v(x,y,0;\mu) = \left\{
                                  \begin{array}{ll}
                                    \mu \sin{(2\pi x)}\cdot \sin{(2\pi y)} \quad &\text{if } (x,y) \in [0,0.5]\times [0,0.5] \\
                                    0 \quad &\text{otherwise}
                                  \end{array}
                                  \right.
\end{align}
where $\mu\in \paramDomain=[0.9,1.1]$ 
% $\param=(\mu_1,\mu_2)\in \paramDomain=[0.9,1.1]\times[0.9,1.1]$ 
is a parameter and $u(x,y,t;\mu)$ and $v(x,y,t;\mu)$ denote the $x$ and $y$
directional velocities, respectively, with $u:\spaceDomain \times [0,2] \times
\paramDomain \rightarrow \RR{}$ and $v:\spaceDomain \times [0,2] \times
\paramDomain \rightarrow \RR{}$ defined as the solutions to
Eq.~\eqref{eq:2dburgers_eq}, and $Re$ is a Reynolds number which is set
$Re=10000$. In the case of $Re=10000$ (an advection-dominated case), the FOM
solution snapshot shows slowly decaying singular values compared to the case of
$Re=100$ as shown in Fig.~\ref{fg:2dburgersSingular} and we observe that a
sharp gradient, i.e., a shock, appears in Figs.~\ref{fg:2dburgersFOMROMS}(a)
and \ref{fg:2dburgersFOMROMS}(d).

Discretizing the space domain $\spaceDomain$ into $\nx-1$ and $\ny-1$ uniform
meshes in $x$ and $y$ directions, respectively, gives $\nx \times \ny$ grid
points $(x_i,y_j)$. $x_i$ is given by $x_i=(i-1)\Delta x$ where $i \in
\{1,2,\cdots, \nx \}$ and $\Delta x = \frac{1}{\nx-1}$ and $y_j$ is given by
$y_j=(j-1)\Delta y$ where $j \in \{1,2,\cdots, \ny \}$ and $\Delta y =
\frac{1}{\ny-1}$. We denote the discrete solutions on grid points as
$u_{i,j}(t;\mu)=u(x_i,y_j,t;\mu)$ and
$v_{i,j}(t;\mu)=v(x_i,y_j,t;\mu)$, where $i\in \nat{\nx}$ and $j\in
\nat{\ny}$. Let $n_{xy}=(\nx-2)\times(\ny-2)$. Then, the backward difference
scheme for the first spatial derivative terms
\begin{align}
    \frac{\partial (\cdot)}{\partial x} &\approx \frac{(\cdot)_{i,j}-(\cdot)_{i-1,j}}{\Delta x},\\
    \frac{\partial (\cdot)}{\partial y} &\approx \frac{(\cdot)_{i,j}-(\cdot)_{i,j-1}}{\Delta y}
\end{align}
and the central difference scheme for the second spatial derivative terms
\begin{align}
    \frac{\partial^2 (\cdot)}{\partial x^2} &\approx \frac{(\cdot)_{i+1,j}-2(\cdot)_{i,j}+(\cdot)_{i-1,j}}{\Delta x^2},\\
    \frac{\partial^2 (\cdot)}{\partial y^2} &\approx \frac{(\cdot)_{i,j+1}-2(\cdot)_{i,j}+(\cdot)_{i,j-1}}{\Delta y^2}
\end{align}
yield the semi-discretized equation which is written by
\begin{align}\label{eq:2dBurgersSemi}
    \frac{d\discreteU}{dt}&=\flux_u(\discreteU,\discreteV),\\
    \frac{d\discreteV}{dt}&=\flux_v(\discreteU,\discreteV)
\end{align}
where $\discreteU=(u_{2,2}, u_{3,2}, \cdots, u_{\nx-2,2},u_{2,3}, u_{3,3},
\cdots, u_{\nx-2,3},\cdots u_{2,\ny-2}, u_{3,\ny-2}, \cdots,
u_{\nx-2,\ny-2})^T$ and $\discreteV=(v_{2,2}, v_{3,2}, \cdots,
v_{\nx-2,2},v_{2,3}, v_{3,3}, \cdots, v_{\nx-2,3},\cdots v_{2,\ny-2},
v_{3,\ny-2}, \cdots, v_{\nx-2,\ny-2})^T$ with superscript $T$ standing for
the transpose operation and $\flux_u: \RR{\nxy} \times \RR{\nxy} \rightarrow
\RR{\nxy}$ and $\flux_v: \RR{\nxy} \times \RR{\nxy} \rightarrow \RR{\nxy}$ are
in the form
\begin{align}
\flux_u(\discreteU,\discreteV) &= 
    -\frac{1}{\Delta x} \discreteU \odot \left(\boldsymbol{M}\discreteU - \boldsymbol{b}_{ux1} \right) 
    -\frac{1}{\Delta y} \discreteV \odot \left(\boldsymbol{N}\discreteU - \boldsymbol{b}_{uy1} \right) \\
    &+\frac{1}{Re \Delta x^2}\left(\boldsymbol{D}_x \discreteU + \boldsymbol{b}_{ux2} \right)
    +\frac{1}{Re \Delta y^2}\left(\boldsymbol{D}_y \discreteU + \boldsymbol{b}_{uy2} \right) \\
\flux_v(\discreteU,\discreteV) &= 
    -\frac{1}{\Delta x} \discreteU \odot \left(\boldsymbol{M}\discreteV - \boldsymbol{b}_{vx1} \right) 
    -\frac{1}{\Delta y} \discreteV \odot \left(\boldsymbol{N}\discreteV - \boldsymbol{b}_{vy1} \right) \\
    &+\frac{1}{Re \Delta x^2}\left(\boldsymbol{D}_x \discreteV + \boldsymbol{b}_{vx2} \right)
    +\frac{1}{Re \Delta y^2}\left(\boldsymbol{D}_y \discreteV + \boldsymbol{b}_{vy2} \right)
\end{align}
where 
\begin{align}
    \boldsymbol{M} &= \bmat{
    \boldsymbol{M}_b & & \\
     & \ddots & \\ 
     & & \boldsymbol{M}_b}_{\nxy \times \nxy},
    \quad
    \boldsymbol{M}_b = \begin{bmatrix}
    1 & & &\\
    -1 & 1 & & \\
     & \ddots & \ddots & \\
     & & -1 & 1
    \end{bmatrix}_{(\nx-2)\times(\nx-2)}, \\
    \boldsymbol{N} &= \begin{bmatrix}
    \boldsymbol{N}_b & & & \\
    -\boldsymbol{N}_b & \boldsymbol{N}_b & & \\
    & \ddots & \ddots & \\
    & & -\boldsymbol{N}_b & \boldsymbol{N}_b
    \end{bmatrix}_{\nxy \times \nxy}, 
    \quad
    \boldsymbol{N}_b  = \begin{bmatrix}
    1 & & \\ 
    & \ddots & \\
    & & 1
    \end{bmatrix}_{(\nx-2)\times(\nx-2)}, \\
    \boldsymbol{b}_{ux1} &=\left(\left(u_{1,2},u_{1,3},\cdots,u_{1,\ny-1}\right)_{1\times(\ny-2)} \otimes \left(1,0,\cdots,0\right)_{1\times(\nx-2)}\right)^T, \\
    \boldsymbol{b}_{uy1} &=\left(\left(1,0,\cdots,0\right)_{1\times(\ny-2)} \otimes \left(u_{2,1},u_{3,1},\cdots,u_{\nx-1,1}\right)_{1\times(\nx-2)}\right)^T, \\
    \boldsymbol{b}_{vx1} &=\left(\left(v_{1,2},v_{1,3},\cdots,v_{1,\ny-1}\right)_{1\times(\ny-2)} \otimes \left(1,0,\cdots,0\right)_{1\times(\nx-2)}\right)^T, \\
    \boldsymbol{b}_{vy1} &=\left(\left(1,0,\cdots,0\right)_{1\times(\ny-2)} \otimes \left(v_{2,1},v_{3,1},\cdots,v_{\nx-1,1}\right)_{1\times(\nx-2)}\right)^T, \\
    \boldsymbol{D}_{x} &= \begin{bmatrix}
    \boldsymbol{D}_{xb} & & \\
    & \ddots & \\
    & & \boldsymbol{D}_{xb}
    \end{bmatrix}_{\nxy \times \nxy},
    \quad
    \boldsymbol{D}_{xb}=\begin{bmatrix}
    -2 & 1 & \\
    1 & \ddots & 1 \\
    & 1 & -2
    \end{bmatrix}_{(\nx-2)\times(\nx-2)},\\
    \boldsymbol{D}_{y} &= \begin{bmatrix}
    -2\boldsymbol{D}_{yb} & \boldsymbol{D}_{yb} & \\
    \boldsymbol{D}_{yb} & \ddots & \boldsymbol{D}_{yb} \\
    & \boldsymbol{D}_{yb} & -2\boldsymbol{D}_{yb}
    \end{bmatrix}_{\nxy \times \nxy},
    \quad
    \boldsymbol{D}_{yb}=\begin{bmatrix}
    1 &  & \\
    & \ddots & \\
    & & 1
    \end{bmatrix}_{(\nx-2)\times(\nx-2)},\\
    \boldsymbol{b}_{ux2}&=(u_{1,2},0,\cdots,0,u_{\nx,2},u_{1,3},0,\cdots,0,u_{\nx,3},\cdots,u_{1,\ny-1},0,\cdots,0,u_{\nx,\ny-1})^T,\\
    \boldsymbol{b}_{uy2}&=(u_{2,1},u_{3,1},\cdots,u_{\nx-1,1},0,\cdots,0,u_{2,\ny},u_{3,\ny},\cdots,u_{\nx-1,\ny})^T,\\
    \boldsymbol{b}_{vx2}&=(v_{1,2},0,\cdots,0,v_{\nx,2},v_{1,3},0,\cdots,0,v_{\nx,3},\cdots,v_{1,\ny-1},0,\cdots,0,v_{\nx,\ny-1})^T,\\
    \boldsymbol{b}_{vy2}&=(v_{2,1},v_{3,1},\cdots,v_{\nx-1,1},0,\cdots,0,v_{2,\ny},v_{3,\ny},\cdots,v_{\nx-1,\ny})^T
\end{align} 
with $\odot$ and $\otimes$ denoting the element-wise multiplication and the
Kronecker product, respectively.

For a time integrator, we use the backward Euler scheme with time step size
$\Delta t = \frac{2}{\nt}$, where $\nt$ is the number of time steps. We set
$\nx=60$, $\ny=60$, and $\nt=1500$.

For the training process, we collect solution snapshots associated with the
parameter $\mu \in \paramDomain_{train}=\{0.9,0.95,1.05,1.1\}$
% $\param = (\mu_1,\mu_2) \in \paramDomain_{train}=\{(0.9,0.9),(0.95,0.95),(1.05,1.05),(1.1,1.1)\}$
such that $\ntrain=4$ at which the FOM is solved. Then, the number of train data points
is $\ntrain\cdot(\nt+1)=6004$ and $10\%$ of the train data are used for
validation purpose. We employ the Adam optimizer \cite{kingma2014adam} with the
SGD and the initial learning rate of $0.001$, which decreases by a factor of
$10$ when a training loss stagnates for $10$ successive training epochs. 
Here, we have two autoencoders. One for $\discreteU$ vector and the other for
$\discreteV$ vector.  The reason why we have such two autoencoders is that we
can use less memory for training process compared to one autoencoder for
$(\discreteU^T,\discreteV^T)^T$ vector and train both of them at the same time.
We set the number of nodes in hidden layer in the encoder, $M_1=6728$, and the
number of nodes in hidden layer in the decoder, $M_2=33730$.  The weights and
bias of the autoencoder are initialized via Kaiming initialization
\cite{he2015delving} as in the first numerical example. The size of the batch is
$240$ and the maximum number of epochs is $10,000$. The training process is
stopped if the loss on the validation dataset stagnates for $200$ epochs.

\begin{figure}[!htbp]
    \centering
    \subfigure[State variable, $u$]{
    \includegraphics[width=0.48\textwidth]{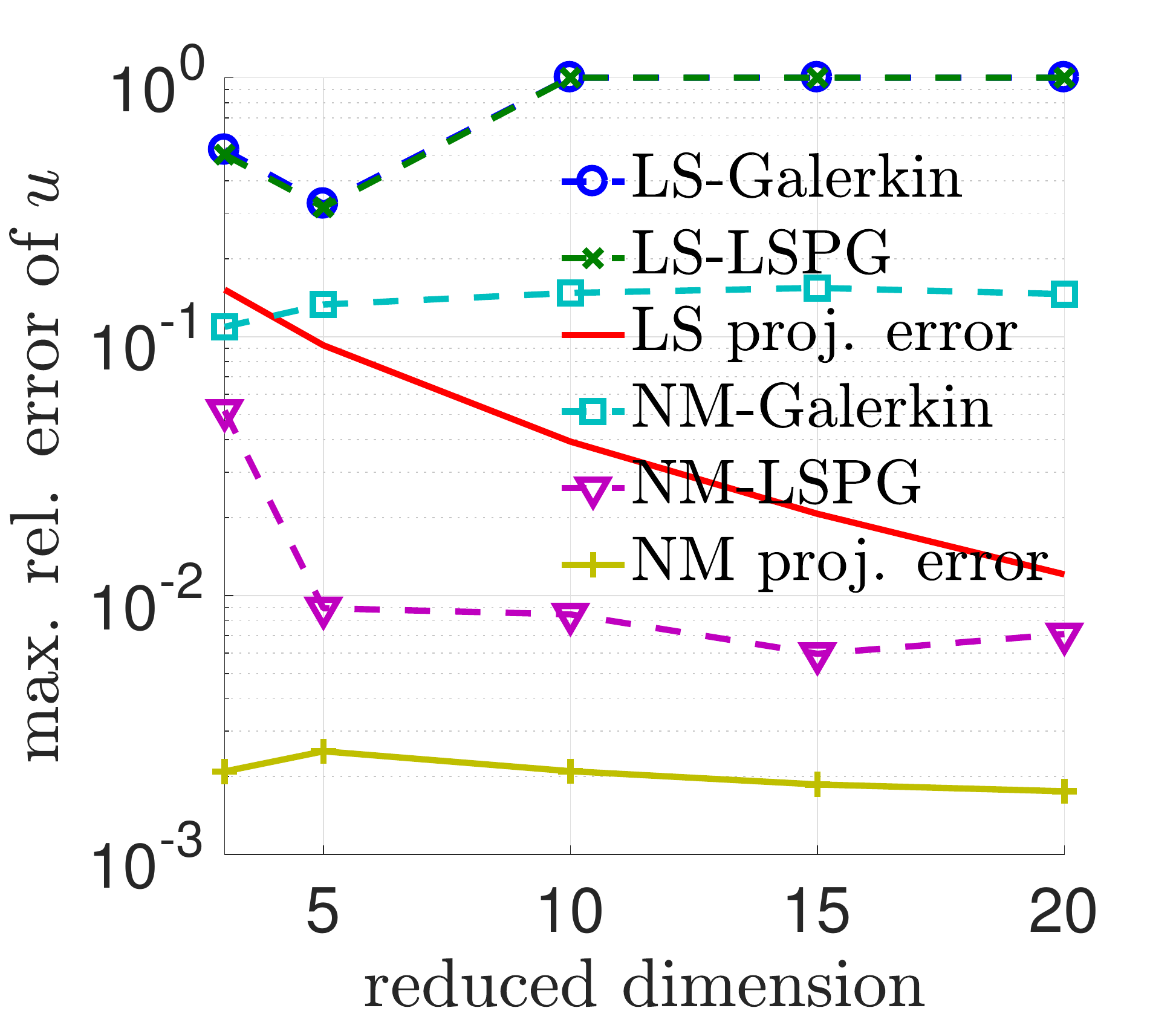}
    }
    \subfigure[State variable, $v$]{
    \includegraphics[width=0.48\textwidth]{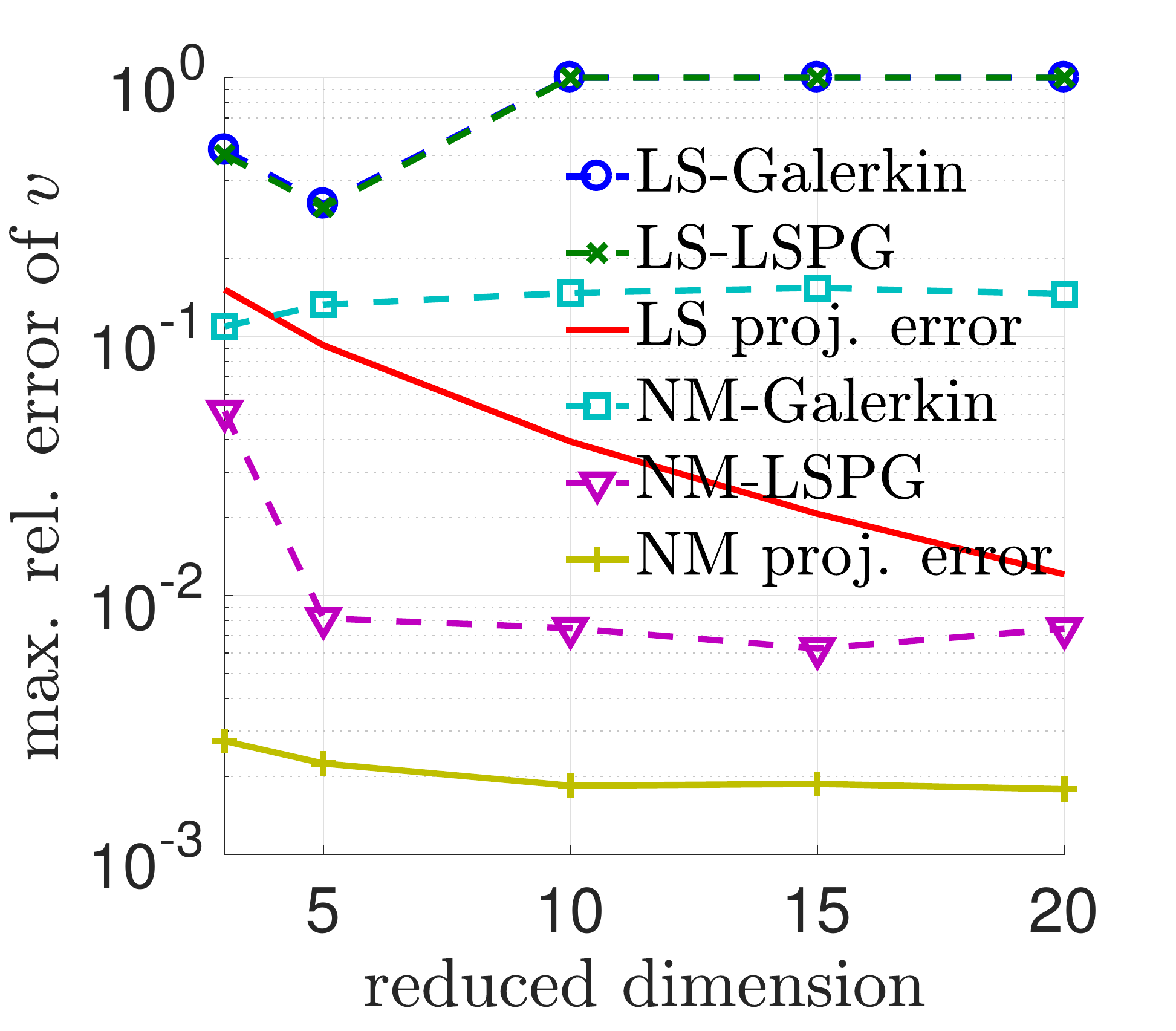}
    }
    \caption{2D Burgers' equation. Relative errors vs reduced dimensions. A
    maximum relative error that is $1$ means the ROM failed to solve the problem.}
    \label{fg:2dErrVSredDim}
\end{figure}

After the training is done, the NM-ROMs and LS-ROMs solve the
Eq.~\eqref{eq:2dburgers_eq} with the target parameter $\mu=1$
% $(\mu_1,\mu_2)=(1,1)$
,which is not included in the train dataset for training the autoencoder and
the linear subspace.  Fig.~\ref{fg:2dErrVSredDim} shows the relative error
versus the reduced dimension $\nbasisspace$ for both NM-ROMs and LS-ROMs.  It
also shows the projection errors for LS-ROMs and NM-ROMs, which are defined in
\eqref{eq:linProjErr} and \eqref{eq:nonlinProjErr}. These are the lower bounds
for LS-ROMs and NM-ROMs, respectively. As expected the relative errors for the
NM-ROMs are lower than the ones for the LS-ROMs.  We even observe that the
relative errors of NM-LSPG are even lower than the lower bounds of LS-ROMs.
One notable observation is that NM-Galerkin is not able to achieve a good
accuracy, while the NM-LSPG does. Another observation is that LS-ROMs struggle
more for this problem than the 1D invisid Burgers' equations, e.g., some
LS-ROMs fail to converge.

\begin{figure}[!htbp]
    \centering
    \subfigure[State variable, $u$]{
    \includegraphics[width=0.48\textwidth]{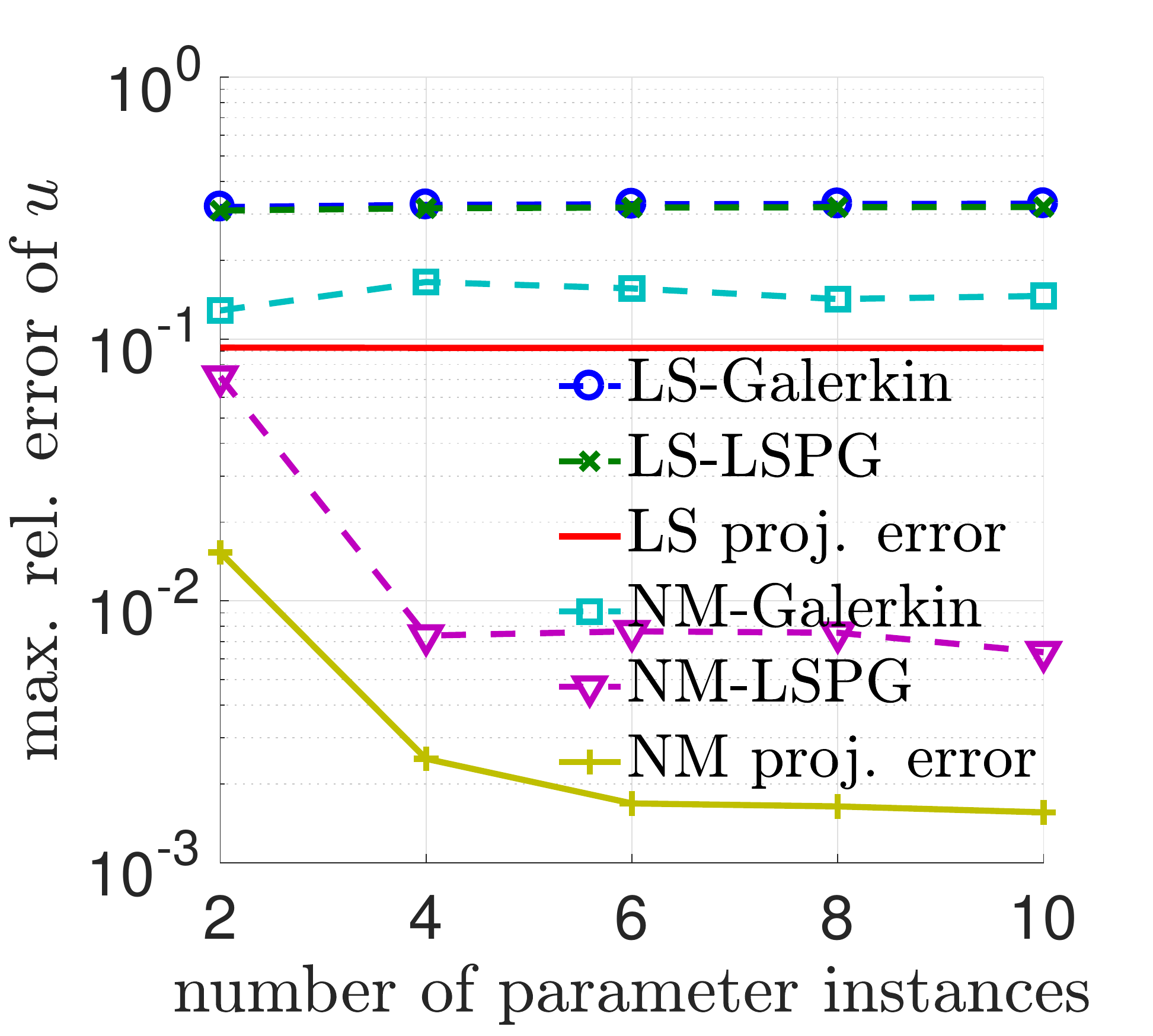}
    }
    \subfigure[State variable, $v$]{
    \includegraphics[width=0.48\textwidth]{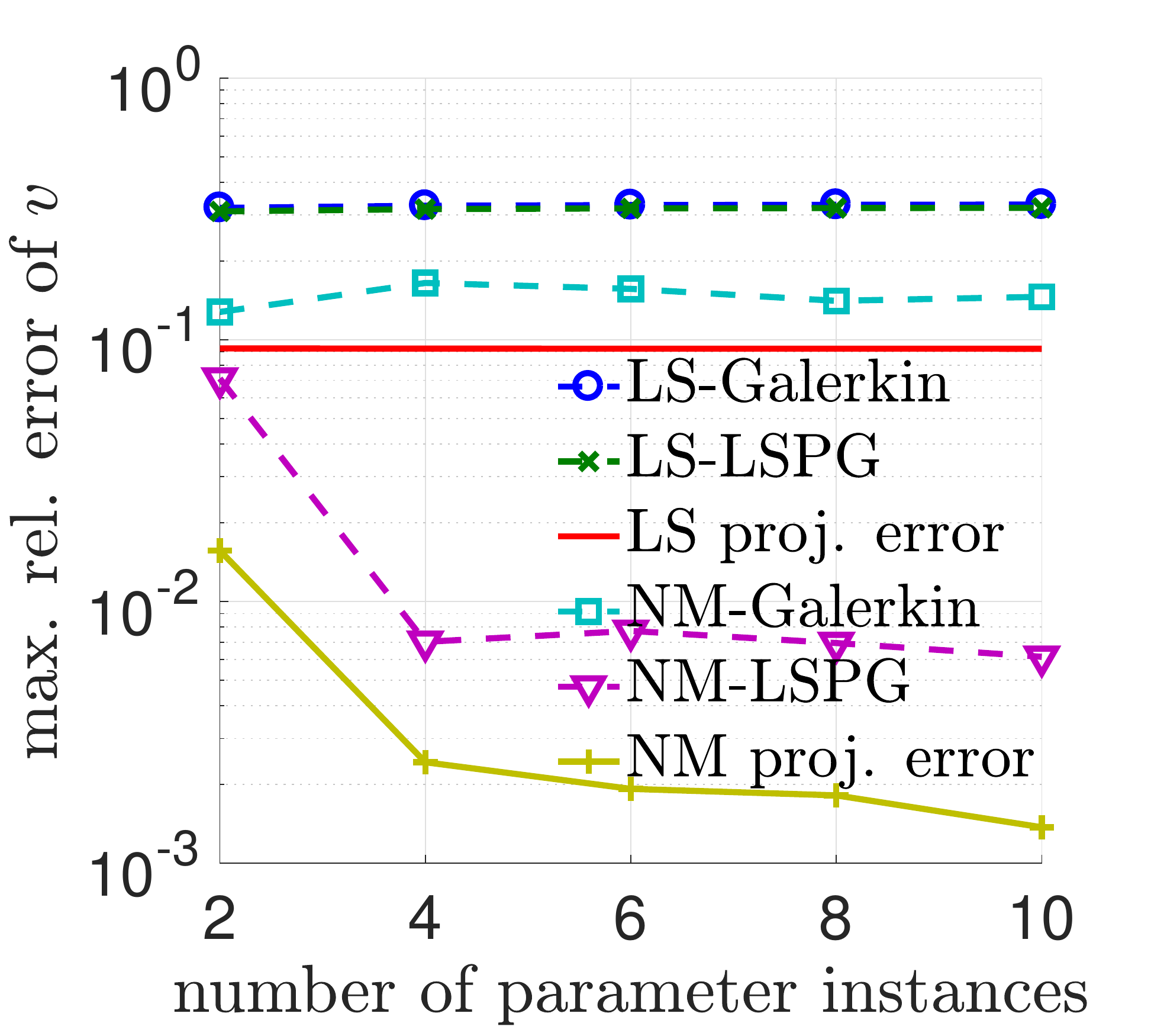}
    }
    \caption{2D Burgers' equation. Relative errors vs the number of parameter
    instances. The following parameter sets are used:
    $\paramDomain_{train}=\{0.9,1.1\}$ for $\ntrain=2$,
    $\paramDomain_{train}=\{0.9,0.95,1.05,1.)\}$ for $\ntrain=4$,
    $\paramDomain_{train}=\{0.9,0.9+\frac{1}{3},0.9+\frac{2}{3},1+\frac{1}{3},1+\frac{2}{3},1.1\}$
    for $\ntrain=6$, and
    $\paramDomain_{train}=\{0.9,0.925,0.95,0.975,1.025,1.05,1.075,1.1\}$ for
    $\ntrain=8$.}
    % \caption{2D Burgers equation. Relative errors vs the number of parameter
    % instances. $\paramDomain_{train}=\{(0.9,0.9),(1.1,1.1)\}$ for $\ntrain=2$,
    % $\paramDomain_{train}=\{(0.9,0.9),(0.95,0.95),(1.05,1.05),(1.1,1.1)\}$ for
    % $\ntrain=4$,
    % $\paramDomain_{train}=\{(0.9,0.9),(0.9+\frac{1}{3},0.9+\frac{1}{3}),(0.9+\frac{2}{3},0.9+\frac{2}{3}),(1+\frac{1}{3},1+\frac{1}{3}),(1+\frac{2}{3},1+\frac{2}{3}),(1.1,1.1)\}$
    % for $\ntrain=6$, and
    % $\paramDomain_{train}=\{(0.9,0.9),(0.925,0.925),(0.95,0.95),(0.975,0.975),(1.025,1.025),(1.05,1.05),(1.075,1.075),(1.1,1.1)\}$
    % for $\ntrain=8$. }
    \label{fg:2dErrVSParam}
\end{figure}

To see the trends regarding the number of training parameter instances, we
increase the number of parameters starting from $\ntrain=2$ with the fixed
reduced dimension $\nbasisspace=5$ to achieve less than $1\%$ maximum relative
error for NM-ROMs.  In Fig.~\ref{fg:2dErrVSParam}, we observe that $\ntrain=4$
is enough.  

Both NM-Galerkin and LS-ROMs without a hyper-reduction do not achieve any
speed-up with respect to the corresponding FOM simulation. For example, the
FOM simulation takes $140.67$ seconds, while the NM-Galerkin, NM-LSPG, LS-Galerkin, and
LS-LSPG  with $\nbasisspace=5$ takes $143.41$, $78.22$, $519.12$ and $2193.70$ seconds,
respectively. Although NM-LSPG is able to achieve a speed-up, it is not
considerable. Therefore, the hyper-reduction is essential to achieve a considerable
speed-up with a reasonable accuracy.  

We compute the maximum relative
error by choosing the larger of the two errors between the maximum relative error of $u$ and the maximum relative error of
$v$. We vary the number of residual basis and
residual samples, with the fixed number of training parameter instances $\ntrain=4$ and the reduced dimension $\nbasisspace=5$, and measure the wall-clock time. The results are shown in
Table~\ref{tb:2dDEIMtest}. Although the LS-LSPG-HR can achieve better speedup
than the NM-LSPG-HR, the relative error of the LS-LSPG-HR is too large to be
reasonable, e.g., the relative errors of around $37 \%$. On the other hand, the
NM-LSPG-HR achieves much better accuracy, i.e., a relative error of around $1
\%$, than the LS-LSPG-HR with a good speedup, i.e., a speedup of higher than
$11$.
\begin{table}[!htbp]
\caption{The top 6 maximum relative errors and wall-clock times at different
  numbers of residual basis and samples which range from $40$ to
  $60$.}\label{tb:2dDEIMtest}
\centering
\resizebox{\textwidth}{!}{\begin{tabular}{|c|c|c|c|c|c|c|c|c|c|c|c|c|}
\hline
 & \multicolumn{6}{|c|}{NM-LSPG-HR} & \multicolumn{6}{|c|}{LS-LSPG-HR}\\
\hline
Residual basis & 55 & 56 & 51 & 53 & 54 & 44 & 59 & 53 & 53 & 53 & 53 & 53\\
\hline
Residual samples & 58 & 59 & 54 & 56 & 57 & 47 & 59 & 58 & 59 & 56 & 55 & 53\\
\hline
Max. rel. error (\%) & 0.93 & 0.94 & 0.95 & 0.97 & 0.97 & 0.98 & 34.38 & 37.73 & 37.84 & 37.95 & 37.96 & 37.97\\
\hline
Wall-clock time (sec) & 12.15 & 12.35 & 12.09 & 12.14 & 12.29 & 12.01 & 5.26 & 5.02 & 4.86 & 5.05 & 4.75 & 7.18\\
\hline
Speed-up & 11.58 & 11.39 & 11.63 & 11.58 & 11.44 & 11.71 & 26.76 & 28.02 & 28.95 & 27.83 & 29.61 & 19.58 \\
\hline
\end{tabular}}
\end{table}

Fig.~\ref{fg:2dburgersFOMROMS} shows solutions at the last time step of
FOM, NM-LSPG-HR, and LS-LSPG-HR with the reduced dimension being
$\nbasisspace=5$. For NM-LSPG-HR, $55$ residual basis dimension and $58$ residual samples are used and for LS-LSPG-HR, $59$ residual basis dimension and $59$ residual samples are used. Both
FOM and NM-LSPG-HR show good agreement in their solutions, while the LS-LSPG-HR
is not able to achieve a good accuracy.
In fact, the NM-LSPG-HR is able to achieve an accuracy as good
as the NM-LSPG for some combinations of the small number of residual basis and
residual samples as in Section~\ref{sec:1dburgers}. 

\begin{figure}[!htbp]
  \centering
  \subfigure[State variable, $u$]{
  \includegraphics[width=0.48\textwidth]{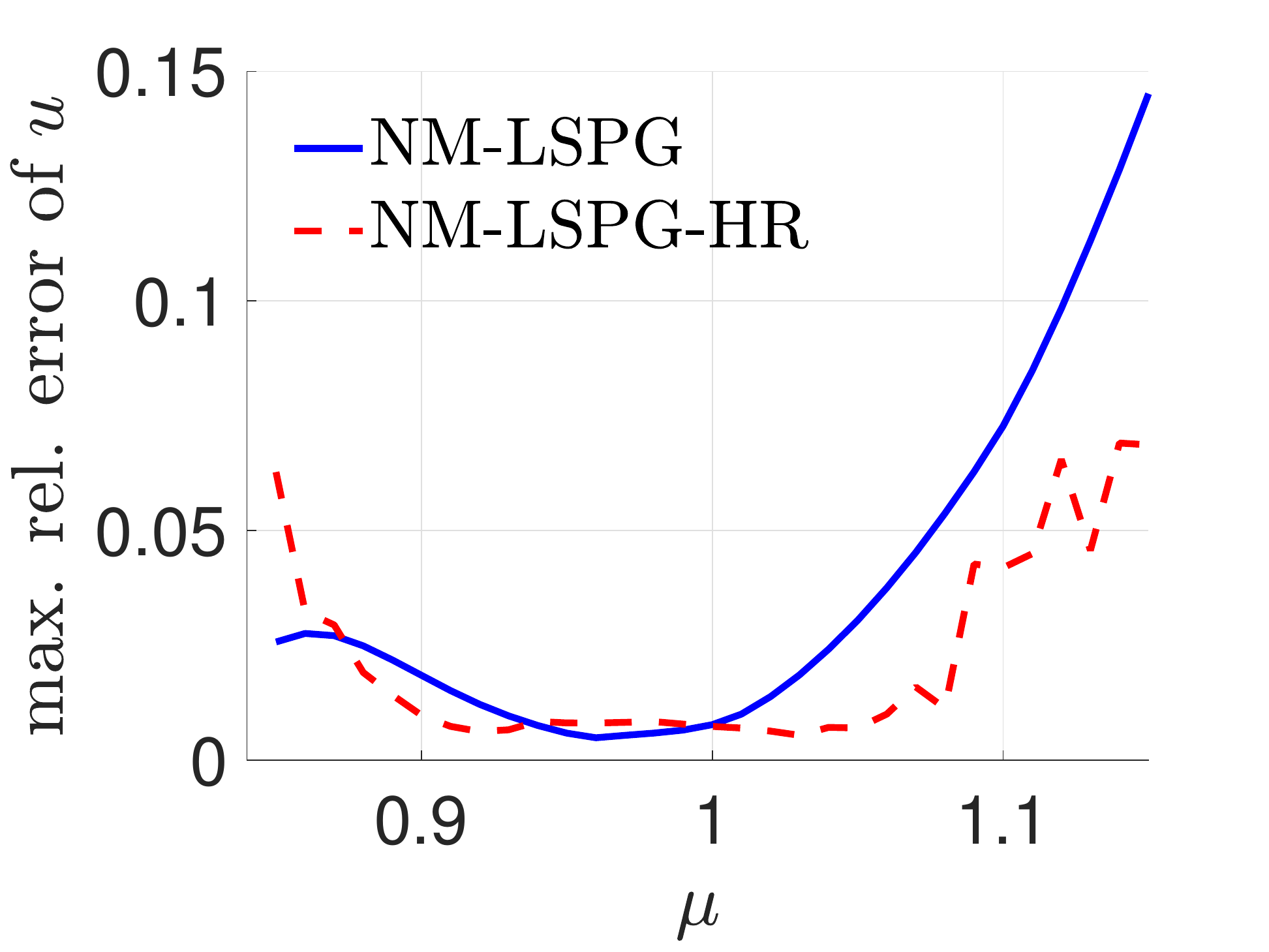}
  }
  \subfigure[State variable, $v$]{
  \includegraphics[width=0.48\textwidth]{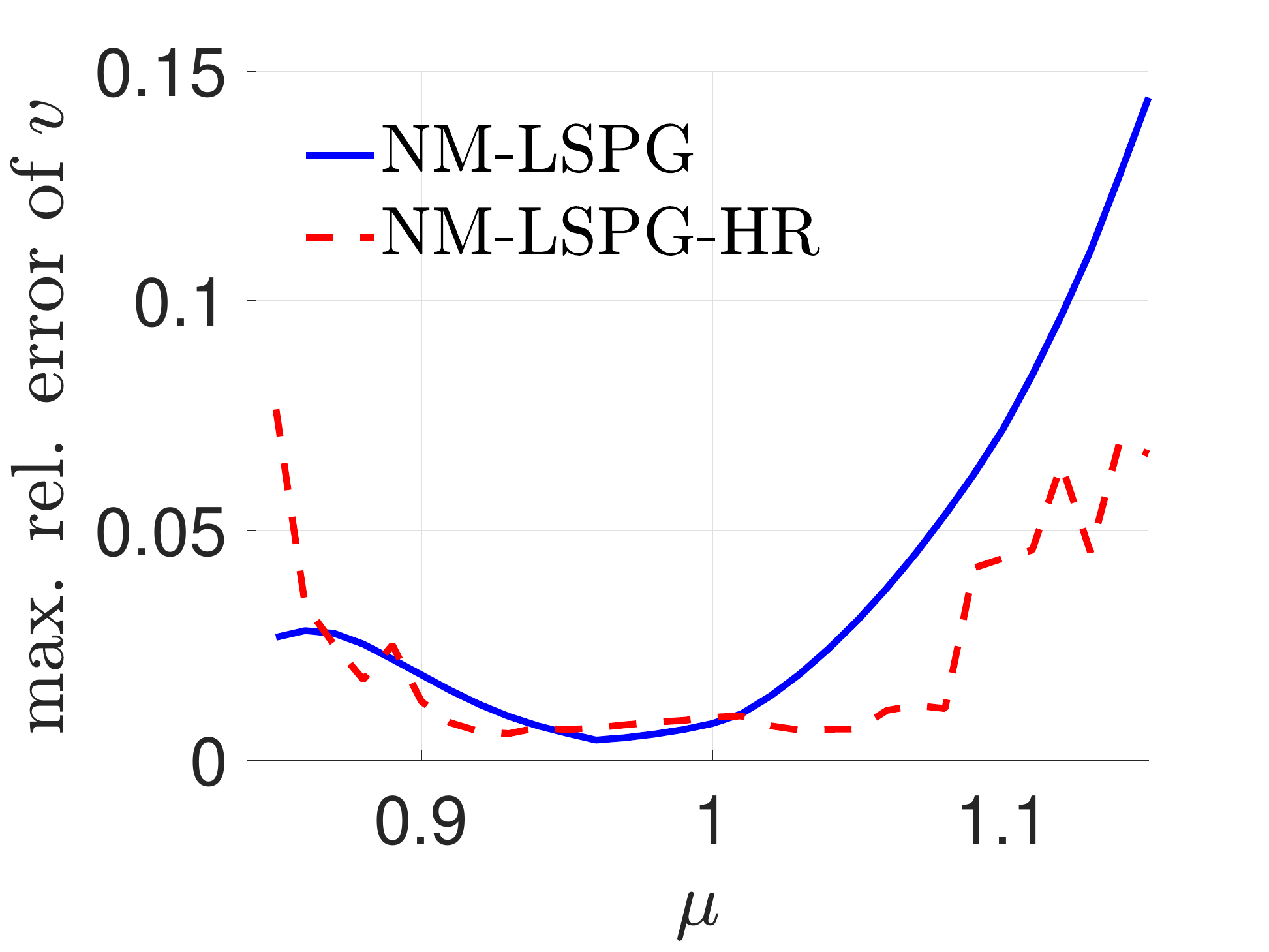}
  }
  \caption{The comparison of the NM-LSPG-HR and NM-LSPG on the maximum relative
  error vs $\mu$}
  \label{fg:2dPredictionTest}
\end{figure}

We look into the numerical tests to see the generalization capability of the
NM-LSPG and NM-LSPG-HR, i.e., the robustness of the NM-LSPG and NM-LSPG-HR
outside of the trained domain. The training sample point set,
$\mu\in\paramDomain_{train}=\{0.9,0.95,1.05,1.1\}$, is used to train a
NM-LSPG-HR. Then the trained NM-LSPG-HR model is used to predict the following
parameter points, $\mu\in\paramDomain_{test} = \{\mu | \mu = 0.85 + 0.01i,
i=0,1, \cdots,30\}$.  The residual basis dimension and the number of residual
samples for each test case are given in Table~\ref{tb:2dPredictionTest}.
Fig.~\ref{fg:2dPredictionTest} shows the maximum relative error over the test
range of the parameter points.  Note that the NM-LSPG and NM-LSPG-HR are the
most accurate within the range of the training points, i.e., $[0.9, 1.1]$.  As
the parameter points go beyond the training parameter domain, the accuracy of
the NM-LSPG and NM-LSPG-HR start to deteriorate gradually. This implies that
the NM-LSPG and NM-LSPG-HR have a trust region. Its trust region should be
determined by an application. For example, if the application is okay with the
maximum relative error of $10$ \%, then the trust region of this particular
NM-LSPG-HR is $[0.85, 1.15]$. However, if the application requires a higher
accuracy, e.g., the maximum relative error of $2$ \%, then the trust region of
the NM-LSPG-HR is $[0.87, 1.08]$.
Note that the average speed-up of the NM-LSPG-HR for all the test cases is
$10.61$ (see Table~\ref{tb:2dPredictionTest}).

\begin{table}[!htbp]
\caption{The residual basis dimension and the number of residual samples for
  each test parameter $\mu$. The wall-clock time and the speed-up of the
  NM-LSPG-HR with respect to the corresponding FOM are also
  reported.}\label{tb:2dPredictionTest}
\centering
\resizebox{0.6\textwidth}{!}{\begin{tabular}{|c|c|c|c|c|}
\hline
$\mu$ & Residual basis & Residual samples & Wall-clock time (sec) & Speed-up \\ \hline
0.85	&	47	&	59	&	13.65	&	10.31	\\	\hline
0.86	&	50	&	50	&	13.19	&	10.66	\\	\hline
0.87	&	45	&	45	&	12.61	&	11.16	\\	\hline
0.88	&	49	&	50	&	12.69	&	11.08	\\	\hline
0.89	&	52	&	52	&	13.41	&	10.49	\\	\hline
0.90	&	53	&	57	&	13.35	&	10.54	\\	\hline
0.91	&	59	&	59	&	13.60	&	10.34	\\	\hline
0.92	&	55	&	58	&	13.41	&	10.49	\\	\hline
0.93	&	51	&	54	&	13.17	&	10.68	\\	\hline
0.94	&	54	&	57	&	13.32	&	10.56	\\	\hline
0.95	&	55	&	58	&	13.52	&	10.40	\\	\hline
0.96	&	55	&	58	&	13.54	&	10.39	\\	\hline
0.97	&	54	&	57	&	13.39	&	10.51	\\	\hline
0.98	&	52	&	55	&	13.20	&	10.66	\\	\hline
0.99	&	52	&	55	&	13.18	&	10.67	\\	\hline
1.00	&	55	&	58	&	13.38	&	10.51	\\	\hline
1.01	&	46	&	49	&	12.80	&	10.99	\\	\hline
1.02	&	50	&	53	&	13.35	&	10.54	\\	\hline
1.03	&	50	&	53	&	13.40	&	10.50	\\	\hline
1.04	&	52	&	53	&	13.40	&	10.50	\\	\hline
1.05	&	46	&	58	&	13.21	&	10.65	\\	\hline
1.06	&	54	&	57	&	13.58	&	10.36	\\	\hline
1.07	&	45	&	57	&	13.20	&	10.66	\\	\hline
1.08	&	45	&	57	&	13.23	&	10.63	\\	\hline
1.09	&	43	&	55	&	13.27	&	10.60	\\	\hline
1.10	&	44	&	48	&	13.31	&	10.57	\\	\hline
1.11	&	40	&	43	&	12.79	&	11.00	\\	\hline
1.12	&	48	&	59	&	13.66	&	10.30	\\	\hline
1.13	&	42	&	51	&	13.25	&	10.62	\\	\hline
1.14	&	46	&	49	&	13.10	&	10.74	\\	\hline
1.15	&	40	&	50	&	13.11	&	10.73	\\	\hline
\end{tabular}}
\end{table}

\begin{figure}[!htbp]
  \centering
  \includegraphics[width=0.7\textwidth]{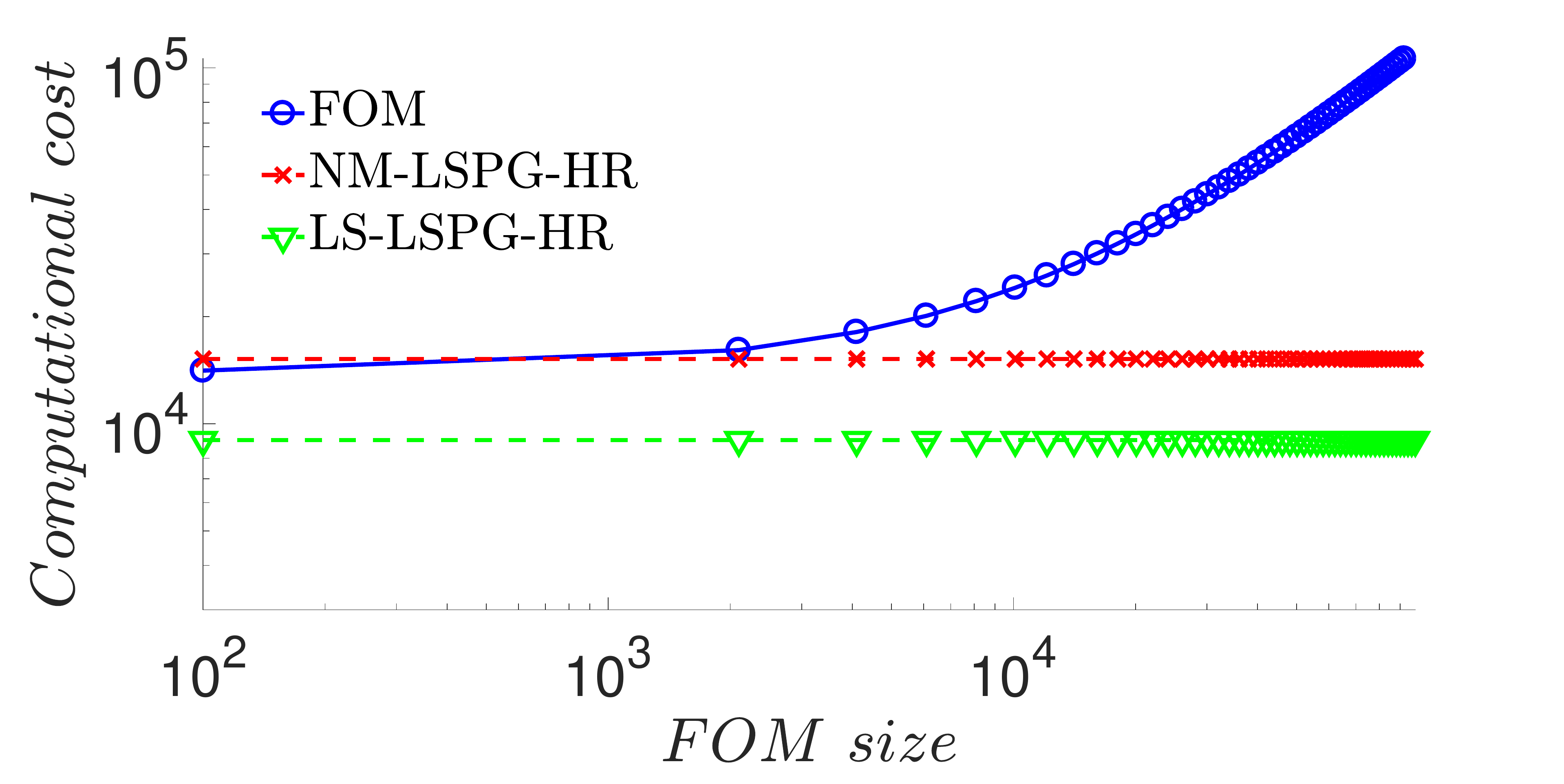}
  \caption{Computational cost vs FOM size. The figure shows that the higher the
  speed up will be achieved, the larger the underlying FOM problem is. The graph
  is generated based on the computational cost analysis done in
  Appendix~\ref{sec:appendixComputationalCosts}}.
  \label{fg:costVSsize}
\end{figure}

\section{Discussion \& conclusion}\label{sec:discussion-conclusion}
In this work, we have successfully developed an accurate and efficient
nonlinear manifold based reduced order model. We demonstrated that the linear
subspace based reduced order model is not able to represent advection-dominated
or sharp gradient solutions of 1D inviscid Burgers' equation and 2D viscous
Burgers' equation with a high Reynolds number. However, our new approach,
NM-LSPG-HR, solves such problems accurately and efficiently. For 1D case, the
NM-LSPG-HR method achieves a good accuracy i.e., the maximum relative error of
around $1\%$ with the speed-up of higher than $2$. For 2D case, the NM-LSPG-HR
method was able to achieve even better accuracy, i.e., the maximum relative
error of less than $1\%$, with even better speed-up of around $12$ than the 1D
case. We also presented \textit{a posteriori} error bounds for NM-Galerkin-HR
and NM-LSPG-HR.  The speed-up of the NM-LSPG-HR is achieved by choosing the
sparse shallow decoder as the nonlinear manifold and applying the efficient
hyper-reduction computation, which can be done by constructing a subnet.
Furthermore, we expect more speed-up as the FOM size increases because the
difference in the computational cost between the FOM and NM-LSPG-HR increases as
shown in Fig.~\ref{fg:costVSsize}.

Compared with the deep neural networks for computer vision and natural language
processing applications, our neural networks are shallow with a small number of
parameters. However, these networks were able to capture the variation in our 1D
and 2D Burgers' simulations. A main future work for transferring this work to
more complex simulations, will be to find the right balance between a shallow
network that is large enough to capture the data variance and yet small enough
to run faster than the FOM. Another future work will be to find an efficient way
of determining the proper size of the residual basis and the number of sample
points \textit{a priori}. To find the optimal size of residual basis and the
number of sample points for hyper-reduced ROMs, we relied on test results. This
issue is not just for NM-LSPG-HR but also for LS-LSPG-HR.

\section*{Acknowledgments}
%The authors acknowledge the helpful comments provided by the anonymous
%reviewers.  
This work was performed at Lawrence Livermore National Laboratory and was
supported by the LDRD program (project 20-FS-007).  Youngkyu was also supported
for this work through generous funding from DTRA. Lawrence Livermore National
Laboratory is operated by Lawrence Livermore National Security, LLC, for the
U.S. Department of Energy, National Nuclear Security Administration under
Contract DE-AC52-07NA27344 and LLNL-JRNL-814844.

\appendix
\section{Time integrators}\label{sec:appendixTimeIntegrators}
\subsection{The linear multistep methods}\label{sec:LM}
Applying a linear $k$-step method to numerically solve Eq.~\eqref{eq:fom} yields an
O$\Delta$E characterized by
 the following system of nonlinear algebraic residual function that needs to be
 solved for the numerical solution
 $\solArg{n}\in\RR{\nspacedof}$ at each time instance:
 \begin{align} \label{eq:LMMresidual}
    \resn(\solArg{n};\solArg{n-1}\ldots,\solArg{n-\kn},\param) &\defeq
    \sum_{j=0}^{\kn}\alpha_j^n\solArg{n-j}
    -\dt\sum_{j=0}^{\kn}\beta_j^n\fluxArg{n-j} \\ 
    &= \zerobold,\quad n\innat{\ntimedof},
 \end{align} 
where coefficients $\alpha_j^n, \beta_j^n\in\RR{}$, $j=0,\ldots,\kn$ define
 a particular linear multistep scheme. It is necessary for consistency to have
 $\alpha_0^n \neq 0$ and
 $\sum_{j=0}^{\kn}\alpha_j^n = 0$. 
Here, $\kn(\leq n)$ denotes the number of steps used by the linear multistep
method at time instance $n$.  The linear multistep methods include the one-step
Euler methods, the implicit Adams--Moulton methods, the explicit
Adams--Bashforth methods, and the Backward Differentiation Formulas (BDFs). 

The second order Adams--Bashforth (AB) method numerically solves 
Eq.~\eqref{eq:fom}, by solving the following nonlinear system of
equations for $\solArg{n}$ at $n$-th time step:
\begin{equation} \label{eq:AF2}
\solArg{n} - \solArg{n-1} = \dt\left (
\frac{3}{2}\fluxArg{n-1}-\frac{1}{2}\fluxArg{n-2}\right ),
\end{equation} 
The residual function of the second AB method is defined as
\begin{align}\label{eq:residual_AF2} 
\begin{split}
\resn_{\ADB}(\solArg{n};\solArg{n-1},\param) &\defeq 
\solArg{n} - \solArg{n-1}
-\dt\left ( \frac{3}{2}\fluxArg{n-1}-\frac{1}{2}\fluxArg{n-2}\right ).
\end{split}
\end{align} 

The second order Adams--Moulton (AM) method numerically solves 
Eq.~\eqref{eq:fom}, by solving the following nonlinear system of
equations for $\solArg{n}$ at $n$-th time step:
\begin{equation} \label{eq:AM2}
\solArg{n} - \solArg{n-1} = \frac{1}{2}\dt(
\fluxArg{n}+\fluxArg{n-1}),
\end{equation} 
The residual function of the second AM method is defined as
\begin{align}\label{eq:residual_AM2} 
\begin{split}
\resn_{\ADM}(\solArg{n};\solArg{n-1},\param) &\defeq 
\solArg{n} - \solArg{n-1}
-\dt\frac{1}{2}(\fluxArg{n}+\fluxArg{n-1}).
\end{split}
\end{align} 
      
The 2-step BDF numerically solves Eq.~\eqref{eq:fom}, by solving the following nonlinear system of equations for $\solArg{n}$ at n-th time step:
\begin{equation} \label{eq:BDF2}
\solArg{n} - \frac{4}{3}\solArg{n-1} + \frac{1}{3}\solArg{n-2}
= \frac{2}{3}\dt\fluxArg{n},
\end{equation}
The residual function of the two-step BDF method is defined as
\begin{align}\label{eq:residual_BDF2}
\begin{split}
\resn_{\BDF}(\solArg{n};\solArg{n-1},\solArg{n-2},\param) &\defeq
\solArg{n} - \frac{4}{3}\solArg{n-1} + \frac{1}{3}\solArg{n-2}
-\frac{2}{3}\dt\fluxArg{n}.
\end{split}
\end{align}

\subsection{The midpoint Runge--Kutta method}\label{sec:RK}
    The midpoint method, a $2$-stage Runge--Kutta method, 
    takes the following two stages to advance at n-th time step of 
    Eq.~\eqref{eq:fom}:
    \begin{equation}\label{eq:RK2}
      \begin{aligned}
        \solArg{n-\frac{1}{2}} &= \solArg{n-1} +
          \frac{\dt}{2}\fluxArg{n-1} \\
        \solArg{n} &= \solArg{n-1} + \dt\fluxArg{n-\frac{1}{2}}.
      \end{aligned}
    \end{equation}

\section{Computational costs}\label{sec:appendixComputationalCosts}
Let's denote the size of FOM and ROM as $m$ and $f$, respectively. Because of
mathematical models that require local information, we need not only the indices
selected from the hyper-reduction, but also their neighbors. We denote the total
number of indices as $z$ and assume $z\ll m$, $z^2 < m$ and $z>f$ (e.g.,
$z=10f$). For simplicity, we assume the mask matrix for the sparse decoder has
the same structure as the mask matrix for 1D Burgers equation as depicted in
Section~\ref{sec:NN}. To generate the mask matrix, two variables $b$ and $\delta
b$ are used, where $b$ denotes the number of nodes in the hidden layer to
compute single output element and $\delta b$ denotes the amount by which the
block of $b$ nodes shifts. Then, the number of nodes in the hidden
layer can be computed as $M_2=b+(m-1) \delta b$.

\subsection{Computational costs of NM-LSPG}\label{sec:ComputationalCostNMLSPG}
Since the decoder is a single hidden layer neural network, the cost of the
decoder and its Jacobian evaluation is $\mathcal{O}(mb)+\mathcal{O}(M_2f)$ and
$\mathcal{O}(fM_2)+\mathcal{O}(mbf)$, respectively. Computing residual,
$\tilde{r}$, includes only element-wise vector calcualtion, resulting in
$\mathcal{O}(m)$. Jacobian of the residual, $\tilde{J}$, can be computed using
row-wise multiplication of matrix and vector because of local connectivity of
mathematical model (e.g., discrete 1D and 2D Burgers equation) in
$\mathcal{O}(fm)$. For the Gauss--Newton method, we need to construct
$\hat{r}=\tilde{J}^T\tilde{r}$ and $\hat{J}=\tilde{J}^T\tilde{J}$, which
requires $\mathcal{O}(fm)$ and $\mathcal{O}(f^2m)$, respectively. It also takes
$\mathcal{O}(f^2)$ to compute each update, $\delta u = -\hat{J}^{-1}\hat{r}$,
iteratively. Assuming the number of Gauss--Newton iterations is in the same
order for the given governing equation, we can factor out the number of
iterations.  Thus, the total computational costs of NM-LSPG for each time step
is $\mathcal{O}(fM_2)+\mathcal{O}(mbf)+\mathcal{O}(f^2m)$. With the assumption
of $M_2 \approx m\delta b$, $f<b$, and $\delta b < b$, we have
$\mathcal{O}(mbf)$.

\subsection{Computational costs of NM-LSPG-HR}\label{sec:ComputationalCostNMLSPGHR}
The size of the weight matrix connecting the hidden layer and the output layer
is reduced to $z$ by $\beta M_2$, where $\beta =\frac{z}{m}$ for the best case
($z$ successive points are selected) and $\beta =1$ for the worst case ($z$
uniformly separate points are selected). Note that when $z$ is small, it is
possible to have $\beta <1$ even for the worst case. Then, replacing $m$ with
$z$ and $M_2$ with $\beta M_2$ in the decoder and its Jacobian evaluation gives
us $\mathcal{O}(zb)+\mathcal{O}(\beta M_2f)$ and $\mathcal{O}(f\beta
M_2)+\mathcal{O}(zbf)$, respectively. Costs of computing residual,
$\tilde{r}_{HR}=Z^T\tilde{r}$ and its Jacobian, $\tilde{J}_{HR}=Z^T\tilde{J}$
for NM-LSPG-HR are $\mathcal{O}(z)$ and $\mathcal{O}(fz)$, respectively because
the sampling matrix $Z^T$ selects $z$ elements of the residual and $z$ rows of
its Jacobian without constructing the sampling matrix. For the Gauss--Newton
method, we need to construct $\hat{r}=\tilde{J}^T_{HR}\obliqueprojector
\tilde{r}_{HR}$ and $\hat{J}=\tilde{J}^T_{HR}\obliqueprojector \tilde{J}_{HR}$,
where $\obliqueprojector$ is the pre-computed $z\times z$ matrix, which
require $\mathcal{O}(fz)+\mathcal{O}(z^2)$ and
$\mathcal{O}(fz^2)+\mathcal{O}(f^2z)$, respectively. It also takes
$\mathcal{O}(f^2)$ to compute each update, $\delta u = -\hat{J}^{-1}\hat{r}$,
iteratively. Assuming the number of Gauss--Newton iterations is in the same
order for the given governing equation, we can factor out the number of
iterations.  Thus, the total computational costs of NM-LSPG-HR for each time
step is $\mathcal{O}(f\beta M_2)+\mathcal{O}(zbf)+\mathcal{O}(fz^2)$. With the
assumption of $M_2 \approx \delta b m$, we have $\mathcal{O}(f\beta \delta b
m)+\mathcal{O}(zbf)+\mathcal{O}(fz^2)$. For the best case, $\beta =
\frac{z}{m}$, the computational costs is $\mathcal{O}(fz\delta
b)+\mathcal{O}(zbf)+\mathcal{O}(fz^2)$. Assuming $\delta b < b$, we have
$\mathcal{O}(zbf)+\mathcal{O}(fz^2)$. For the worst case, $\beta =1$, we have
$\mathcal{O}(f\delta b m)+\mathcal{O}(zbf)+\mathcal{O}(fz^2)$.

\subsection{Computational costs of LS-LSPG}\label{sec:ComputationalCostLSLSPG}
The decoder $\scaledDecoder(\redsolapprox)$ and its Jacobian
$J_g(\redsolapprox)$ are replaced with $\basismatspace \redsolapprox$ and
$\basismatspace$, respectively. Thus, the cost of $\basismatspace \redsolapprox$
is $\mathcal{O}(mf)$ and the cost of its Jacobian evaluation is zero. The costs
of computing residual and its Jacobian are the same as for NM-LSPG . Also, the
costs of the Gauss--Newton method is the same as in
Section~\ref{sec:ComputationalCostNMLSPG}. Assuming the number of Gauss--Newton
iterations is in the same order for the given governing equation, we can factor
out the number of iterations. Thus, the total computational costs of LS-LSPG for
each time step is $\mathcal{O}(f^2m)$.

\subsection{Computational costs of LS-LSPG-HR}\label{sec:ComputationalCostNLSSPGHR}
For LS-LSPG-HR, we construct reduced model with the size of basis matrix
$\basismatspace_{HR}$ being $z$ by $f$, where
$\basismatspace_{HR}:=\samplemat\basismatspace$. Thus, the costs of $\basismatspace_{HR}
\redsolapprox$ is $\mathcal{O}(zf)$ and the cost of its Jacobian evaluation is
zero. The costs of computing residual and its Jacobian are the same as for
NM-LSPG-HR. Also, the costs of the Gauss--Newton method is the same as in
Section~\ref{sec:ComputationalCostNMLSPGHR}. Assuming the number of
Gauss--Newton iterations is in the same order for the given governing equation,
we can factor out the number of iterations. Thus, the total computational costs
of LS-LSPG-HR for each time step is $\mathcal{O}(f^2z)+\mathcal{O}(fz^2)$.

\bibliographystyle{plain}
\bibliography{references}

\end{document}